\documentclass[a4paper,12pt]{article}
\usepackage{fullpage}

\usepackage[]{graphicx}\usepackage[]{color}

\makeatletter
\def\maxwidth{ %
  \ifdim\Gin@nat@width>\linewidth
    \linewidth
  \else
    \Gin@nat@width
  \fi
}
\makeatother

\usepackage{subfigure}
\usepackage{booktabs} 
\usepackage{natbib}
\usepackage{hyperref}
\usepackage{algorithm}
\usepackage{algpseudocode}
\usepackage{algorithmicx}
\usepackage{thumbpdf} %recommanded by jss

\usepackage{amsmath}
\usepackage{amssymb}
\usepackage{mathtools}
\usepackage{amsthm}

\usepackage[capitalize,noabbrev]{cleveref}

\usepackage{autobreak}
\usepackage{setspace}

\newcommand{\x}{\mathbf{x}}
\newcommand{\y}{\mathbf{y}}

\newcommand{\A}{\mathbf{A}}
\newcommand{\C}{\mathbf{C}}
\newcommand{\Id}{\mathbf{I}}
\newcommand{\U}{\mathbf{U}}
\newcommand{\W}{\mathbf{W}}
\newcommand{\X}{\mathcal{X}}
\newcommand{\Xset}{\mathbf{X}}
\newcommand{\XC}{\mathbf{X}_C}
\newcommand{\XE}{\mathbf{X}_E}
\newcommand{\XR}{\mathbf{X}_R}

\newcommand{\R}{\mathbb{R}}
\newcommand{\veck}{\mathbf{k}}

\newcommand{\K}{\mathbf{K}}

\newcommand{\Cov}[1]{\mathrm{Cov}\left[{#1}\right]}
\newcommand{\KC}{\mathbf{K}_C}
\newcommand{\KE}{\mathbf{K}_E}
\newcommand{\kC}{\mathbf{k}_C}
\newcommand{\kE}{\mathbf{k}_E}
\newcommand{\Wij}{\mathbf{W}_{i,j}}

\title{Combining additivity and active subspaces for high-dimensional Gaussian process modeling}

\author{Micka\"{e}l Binois\thanks{Corresponding author: Inria, Université Côte d’Azur, CNRS, LJAD, Sophia Antipolis, France \href{mailto:mickael.binois@inria.fr}{\tt mickael.binois@inria.fr}}
\and
Victor Picheny\thanks{Secondmind, Cambridge, UK}
}

\begin{document}

\maketitle

\begin{abstract}
Gaussian processes are a widely embraced technique for regression due to their
good prediction accuracy, analytical tractability, and built-in capabilities
for uncertainty quantification. However, they suffer from the curse of
dimensionality whenever the number of variables increases. This challenge is
generally addressed by assuming additional structure in the problem, the
preferred options being either additivity or low intrinsic dimensionality.
After discussing their relative merits, our contribution for high-dimensional
Gaussian process modeling is to combine them with a multi-fidelity strategy.
We detail the corresponding construction and showcase the advantages through
experiments on synthetic functions and datasets.
\end{abstract}

\onehalfspacing

\section{Introduction}

As a surrogate modeling option, Gaussian processes (GPs), also known as
kriging, enjoy widespread use across applied scientific domains, including
engineering, machine learning and physics \citep[see,
e.g.,][] {williams2006gaussian,Gramacy2020}. Appreciated for their efficiency
on small datasets, GPs offer a full predictive distribution in closed form
and are readily accessible from numerous software packages. Nonetheless, the
case of many input variables remains one of the most challenging topic in GP
modeling, and particularly in its use within Bayesian
optimization(BO) \citep [see, e.g.,][]{garnett2023bayesian}. The root of the
issue is that the typical covariance kernels employed in applications, e.g.,
the squared exponential or Matérn ones, hinge on Euclidean and absolute
distances between data points, as discussed, e.g., by \cite{wilson16}. As the
dimension increases, so do the distances between designs. Hence, in
high-dimensional spaces, GP predictions predominantly operate in the
extrapolation regime, where the selection of a trend (or mean function) is
critical, see, e.g., \cite{Journel1989}.

Among possible structural assumptions to scale with respect to the number of
variables, three main categories have emerged. The first one is to identify
the most important variables and then reduce the problem dimensionality. This
adaptive variable selection is applied for instance
by \cite {cao2022scalable}, while remaining efficient in both the number of
variables and dimensions through the use of the Vecchia approximation, which
restricts the conditioning set of each observations to only a few designs. A
second, more general approach, entails assuming that the problem has a low
intrinsic dimensionality, meaning that the variation of the function is
concentrated on a few directions only. This is also referred to as using
linear embeddings \cite{Wang2013,Garnett2013} or active
subspaces \citep[AS,][] {constantine2014active,Eriksson2019}. Hence, both
approach avoid dealing directly with high-dimensional distances. A third
perspective is to consider additive decompositions of the function, where
components involve only a few variables, thereby limiting the degree of
interaction between variables, see
e.g., \cite{Duvenaud2011,Durrande2010,rolland2018high,Lu2022a}.
In \citep {Ginsbourger2016}, a regular product Gaussian kernel is decomposed
into ANOVA terms, allowing the resulting GP to undergo a similar
decomposition. The originality is to separate elements into additive and
ortho-additive components (i.e., those that capture all the non-additive
parts). For a more comprehensive review of high-dimensional GPs and BO, we
refer to \cite {malu2021bayesian,binois2022survey}.

Here we focus on these two promising avenues that have been explored
separately: the use of linear embeddings and of additive models. Each
approach comes with its own set of strengths and limitations. Linear
embeddings offer great scalability and allow to capture complex interactions,
but only as long as the intrinsic dimension is low. Besides, the intrinsic
dimension is generally not known beforehand and the method is very sensitive
to its estimation, as it is to the choice of the embedding. Additive models
typically capture very well the main trends from high-dimensional data.
However, as the number of possible interaction terms makes the full inference
combinatorially intractable, with too many hyperparameters, it is common
practice to limit the number of interactions to a manageable few, which
severely limits the expressivity of those models. Additional concerns may
arise in BO due to the vanishing variance at unexplored
locations \cite {Durrande2010}.

Given the distinct relative advantages of each approach -- namely the
scalability and interpretability offered by additive GPs and the ability to
capture high-order interactions through AS-based GPs -- we aim to propose a
hybrid model that incorporates the strengths of both. As each can capture the
features of the other, with high-order interactions or large intrinsic
dimensionality, combining them is not trivial. To address identifiability
issues, we introduce orthogonality by adopting a multi-fidelity approach,
typically used when inexpensive yet coarse approximations of the target
function are available, see,
e.g., \cite {kennedy2000predicting,brevault2020overview}. In a nutshell, our
model will consider two fidelity levels: a coarse level corresponding to a
first-order additive model and a high-fidelity level corresponding to a GP on
an active subspace whose dimension is learned.

Our contributions are the following:
\begin{itemize}
    \item We concentrate on two efficient structural assumptions for
    high-dimensional GP modeling, namely additivity and linear embeddings. We
    further detail their respective strengths and weaknesses, plus the
    complexity of a naive combination.
    \item We develop a multi-fidelity approach designed to efficiently perform
     this combination, capturing additive and linear embedding
     contributions.
    \item For AS-based GPs, we discuss the inference of the intrinsic dimension
    within one- or two-stage methods.
    \item We conduct a thorough comparison of our multi-fidelity method
     against baselines within an extensive benchmark comprising test
     functions and datasets, focusing on the moderate-dimensional, small data
     regime. Our findings confirm that the multi-fidelity approach improves
     over a standard GP when either additivity or active subspaces are
     present. Importantly, the performance remains robust when such
     structures are absent except sometimes on the smallest data regime
     examples.
\end{itemize}

\section{Background}

We want to fit a Gaussian process model of
$f : \Xset \subset \R^d \rightarrow \R$ when $d$ is relatively large.
What \emph{large} means depends on the dataset size and the complexity of the
problem at hand. In the derivative free black-box context, it is generally
considered that ten variables problems already fall within the realm of high
dimensionality. We briefly outline Gaussian process regression, the additive
and linear embedding versions, before introducing the multi-fidelity model.

\subsection{Gaussian Processes Regression}

Given $n \in \mathbb{N}$ input designs $\x^{(i)} \in \Xset$ with corresponding
observations $f(\x^{(i)}) = y_i$ (possibly noisy), GPs are a form of spatial
modeling that only depends on a mean and a covariance function. Typically,
the mean function is taken to be zero, and all the modeling effort is placed
on the covariance function $k$. From this GP prior, the posterior
distribution is another GP and the prediction at any $\x$ follows: $Y(\x)
| (\x, y_i)_{1
\leq i \leq n} \sim
\mathcal{N} \left(m_n(\x), s_n^{2}(\x)\right)$ where, see, e.g.,
\cite{williams2006gaussian,Gramacy2020}:
\begin{align*}
m_n(\x) &= \veck(\x)^\top \K^{-1} \y,\\
s_n^{2}(\x) &= k(\x, \x) - \veck(\x)^\top \K^{-1} \veck(\x)
\end{align*}
with $\y:= (y_1, \dots,
y_n)$, $\veck(\x) := (k(\x, \x^{(i)}))_{1 \leq i \leq n}$, $\K:= (k(\x^{(i)},
\x^{(j)}) + \tau^2 \mathbf{1}_{i=j})_{1 \leq i,j \leq n}$. $\tau^2$ is the noise
hyperparameter, when assuming $y_i = f(\x^{(i)}) + \varepsilon_i$, with
$\varepsilon_i \sim \mathcal{N}(0, \tau^2)$.

The covariance kernel function must be a positive definite function. In
practice, parameterized families such as the Gaussian and Matérn covariances
are employed, see, e.g., \cite{williams2006gaussian}. As an example, the
Matérn 5/2 kernel in product form writes: $k(\x, \x') = \sigma^2
\prod \limits_{i=1}^d k_i(x_i, x'_i)$ with $k_i(x_i, x'_i) = \left(1 +
\sqrt{5} h_i / \theta_i + 5 h_i^2/(3\theta_i^2) \right) \times \allowbreak
\exp{\left(-\sqrt{5}h_i/\theta_i\right)}$, $h_i = |x_i - x'_i|$. For inferring
 the hyperparameters, we rely on the log-likelihood: $\log L := - n/2 \log
 (2\pi) - 1/2 \log|\K| - 1/2 \y^\top \K^{-1} \y$. When the variance parameter
 $\sigma^2$ can be factorized, i.e., $\K = \sigma^2 \mathbf{R}$, with
 $\mathbf{R}$ the correlation matrix, its estimator is available in
 closed-form: $\hat{\sigma}_n^2 := n ^{-1} \y^\top  \mathbf{R}^
 {-1} \y$, while the other hyperparameters are obtained by maximizing the
 concentrated log-likelihood: $\log \tilde{L} := -n/2 \log(2 \pi) - n/2
\log(\hat{\sigma}_n^2) - 1/2
\log|\mathbf{R}| - n/2$.

A typical extension of GPs to offer much better scalability with data size is
to follow the Sparse Variational GP (SVGP)
framework \citep {titsias2009,hensman2013}. While not considered in this
paper, our model would naturally apply to this framework.

\subsection{Additive Model}

Product kernels can quickly decrease to zero as two points move apart,
reducing the modelling generalization. However sum kernels do not suffer this
problem. This latter form of covariance amounts to considering additive
models, that is, decompositions of the original function into several
components. The general model writes $f(\x) \approx \mu + \sum \limits_{i =
1}^M g_i(\x_{A_i})$ with component functions $g_i$ acting on subsets of
variables $A_i$, plus a constant term $\mu$. These subsets can simply be the
original variables
\cite{neal1997monte,plate1999accuracy,Duvenaud2011,Durrande2010}, disjoint
groups \cite{Kandasamy2015,Gardner2017} or more general subsets of variables
\cite{rolland2018high}.

The sum form of the covariance, i.e., $k_A(\x, \x') = \sum \limits_{i = 1}^d
k_i(x_i, x'_i)$, translates in the model, where the mean becomes the sum of
component-wise means $m_{n,A}(\x) = \veck_A(\x)^\top
\K_A^{-1} \y = \sum \limits_{i = 1}^d \veck_i(x_i) \K_A^ {-1} \y
 = \sum \limits_{i = 1}^d m_{n,i} (x_i)$ with $\veck_i(x_i):= (k_i (x_i,
 x^{(j)}_i))_{1 \leq j \leq n}$. It becomes useful for visualization and
 interpretation, e.g., with main effect plots. As for the predictive
 variance, it does not have a similar decomposition but it can be zero at
 unobserved locations, unless noise is present, see, e.g.,
\cite{Durrande2010}.

Inference can involve learning variance and scale parameters for every
component kernel, plus eventually selecting interaction order, with up to
$2^d$ components. To help inference, centering the various terms is usually
preferred to avoid non-identifiability \cite {Durrande2010,Lu2022a}.
Orthogonality constraints can be further added between the terms, leading to
functional ANOVA decomposition of the original function, see, e.g.,
\cite{Muehlenstaedt2012,Durrande2013,Ginsbourger2016}.

\subsection{Active Subspace Methods}

By not imposing the variables to match the original variables in dimension
reduction, one can rather attempt to learn the most important directions of
variation of $f$: $f(\x) \approx g(\A^\top \x)$ with $\A$ a $d \times r$
matrix, $1 \leq r \leq d$ and preferably $r \ll d$. Learning this linear
embedding encoded in $\A$ is possible with different strategies. Elements of
$\A$ can be treated as regular hyperparameters
\cite{Garnett2013,Tripathy2016,letham2020re}, or they can be random
\cite{Wang2013,nayebi2019framework}, relying on the stability of the random
projection for the $L_2$ norm.

When looking at directions where $f$ varies the most, the so-called active
subspace \cite{Constantine2015}, $\A$ is defined (up to a rotation) as the
largest $r$ eigen-vectors of the matrix $\mathbf{C} := \int_\Xset
\nabla(f(\x))^\top
\nabla(f(\x)) \lambda(d\x)$ where $\lambda$ is usually the Lebesgue measure on
hypercubic domains. Without the gradient of $f$, $\A$ may be estimated via
compressed sensing, partial least squares, principal component analysis, see
e.g., \cite{carpentier2012bandit,Djolonga2013,Bouhlel2016,Raponi2020}. For a
GP, its AS matrix $\C^{(n)}$ can be directly computed, as shown by
\cite{Wycoff2021} and detailed in Appendix \ref{ap:asgp}. These AS approaches
usually involve first learning a high-dimensional GP to estimate
$\A$, before fitting a low dimensional GP in the reduced space, see, e.g.,
\cite{Tripathy2016}. In Appendix \ref{ap:asgp}, we also show how to learn
directly the low dimensional GP.

\subsection{Multi-fidelity}
\label{sec:armufi}

Be it a number of Monte Carlo iterations, a mesh or a training set size, the
accuracy of a simulator experiment is often tunable. Accordingly, GP models
have been adapted to take into account these various levels of fidelity, see,
e.g.,
\cite{kennedy2000predicting,Forrester2008,le2014recursive,tighineanu2022transfer}.
We only review the two levels case here, coarse (resp.\ fine) level denoted by
$C$ (resp.\ $E$), with the auto-regressive (AR) model: $f_E(\x) = \rho f_C(\x) +
\delta(\x)$, $\delta(\cdot) \perp f_C(\cdot)$. This model, proposed by
\cite{kennedy2000predicting}, assumes that $\forall \x \neq \x'$,
$\Cov{Y_E(\x), Y_C(\x') | Y_C(\x)} = 0$, i.e., that nothing more can be learned
for $Y_E(\x)$ from the coarse model once $Y_C(\x)$ is known.

Denote the $n_C$ observations $\y_C$ (resp.\ $\y_E$) at $\XC: =
\left(\x_C^{(1)}, \dots, \x_C^{(n_C)}\right)$ (resp.\ $\XE$). Given the
following covariances:
\begin{align*}
\Cov{Y_C(\x), Y_C(\x')} &= k_C(\x, \x'),\\
\Cov{Y_E(\x), Y_C(\x')} &= \rho k_C(\x, \x'),\\
\Cov{Y_E(\x), Y_E(\x')} &= \rho^2 k_C(\x, \x') + k_E(\x, \x'),
\end{align*}
the corresponding predictive equations for the zero mean version are given by
(see, e.g., \cite{kennedy2000predicting,Forrester2008} for the derivation):
\begin{align}
\begin{split}
m_{n,E}(\x) &= \tilde{\veck}(\x)^\top \tilde{\K}^{-1} \tilde{\y},\\
s_{n,E}^{2}(\x) &= \rho^2 k_C(\x, \x) + k_E(\x, \x) - \tilde{\veck}(\x)^\top \tilde{\K}^{-1} \tilde{\veck}(\x)
\label{eq:armufi}
\end{split}
\end{align}
with $\tilde{\veck}(\x)^\top = [\rho k_C(\XC, \x), \rho^2 k_C(\XE, \x) + k_E(\XE, \x)]$,\\
$\tilde{\K} = \begin{bmatrix}
k_C(\XC, \XC) & \rho  k_C(\XC, \XE) \\
 \rho  k_C(\XE, \XC) & \rho^2  k_C(\XE, \XE) + k_E(\XE, \XE)
\end{bmatrix}$
and $\tilde{\y}^\top = [\y_C, \y_E]$.

For inference, the low fidelity model is independently trained first, then the
fine level hyperparameters (including $\rho$) are obtained based on
$\mathbf {d} := \y_E - \rho \y_C(\XE)$. If the designs of experiments between
fidelity levels are nested, i.e., $\XE \subseteq \XC$, the difference between
levels can be directly evaluated. Otherwise, the difference can be computed
based on the predictive mean \cite{Forrester2008,sacher2021non}, replacing
$\y_C (\XE)$ by its prediction. Further details are given in Appendix
\ref{ap:mufi}.

A recursive formulation of Eq.\ (\ref{eq:armufi}) is available to reduce the
computational effort, see, e.g., \cite{le2014recursive}, but equivalence
holds only in the noiseless setting, see Appendix \ref{ap:mufi}. Subsequently
we introduce our proposed multi-fidelity combination of additive and linear
embedding models, tailored to tackle high-dimensional problems.

\section{Multi-fidelity for High-dimensional Modeling}

\begin{figure*}[htpb]
\centering
\includegraphics[width=0.3333\textwidth, trim= 20 35 20 40, clip]{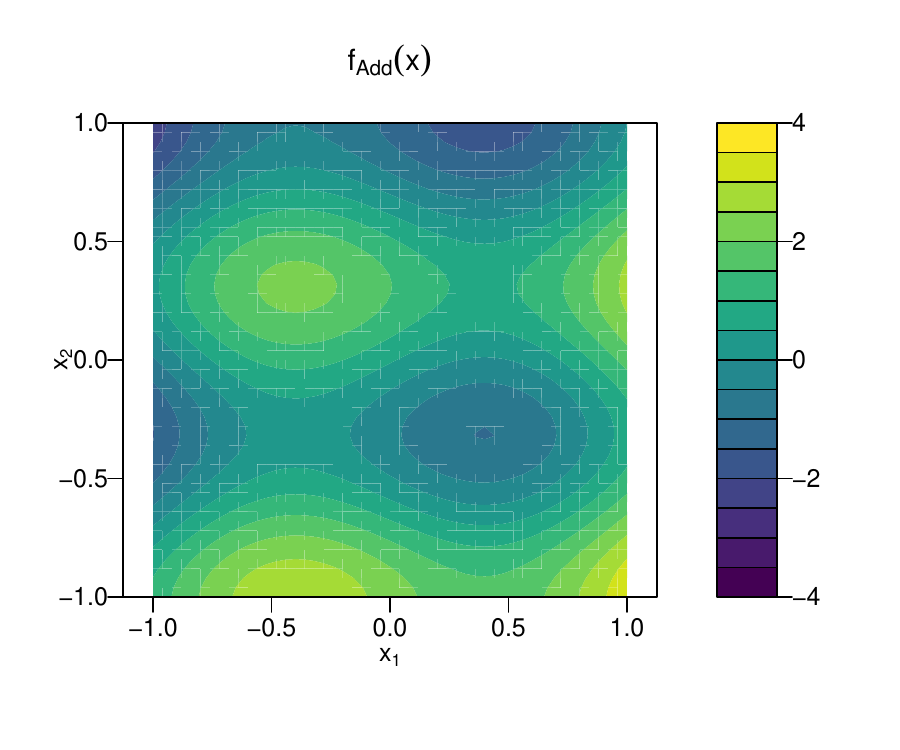}%
\includegraphics[width=0.3333\textwidth, trim= 20 35 20 40, clip]{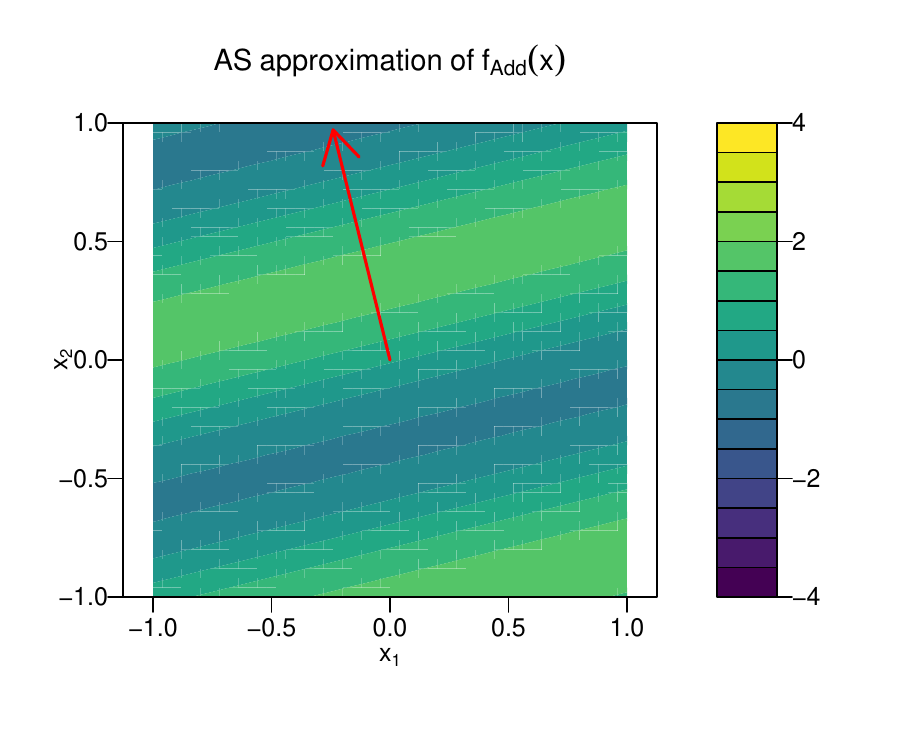}%
\includegraphics[width=0.3333\textwidth, trim= 20 35 20 40, clip]{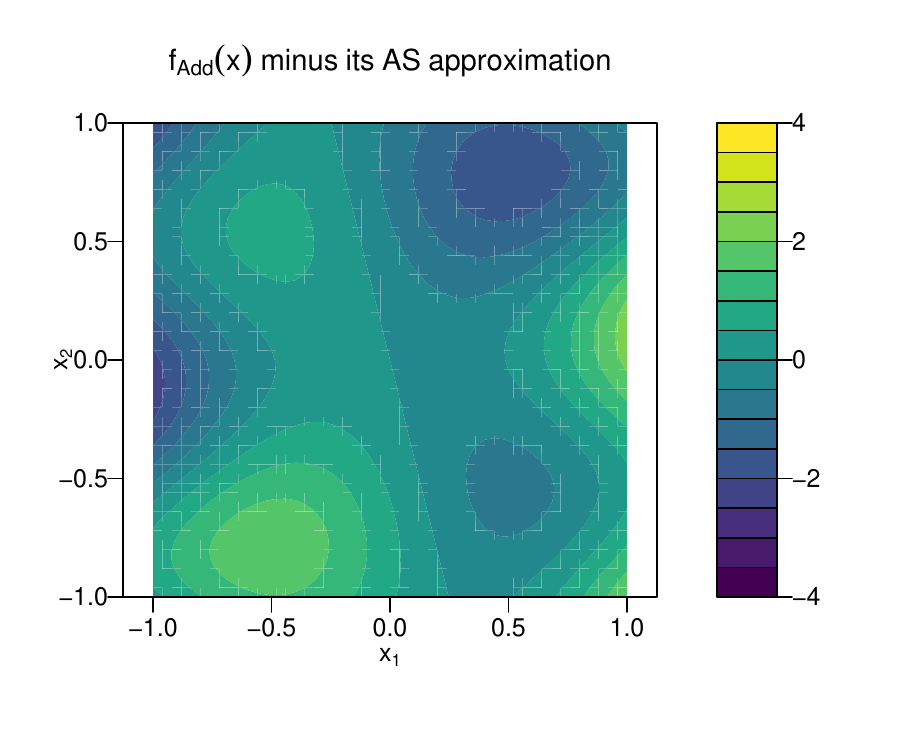}\\
\includegraphics[width=0.3333\textwidth, trim= 20 35 20 40, clip]{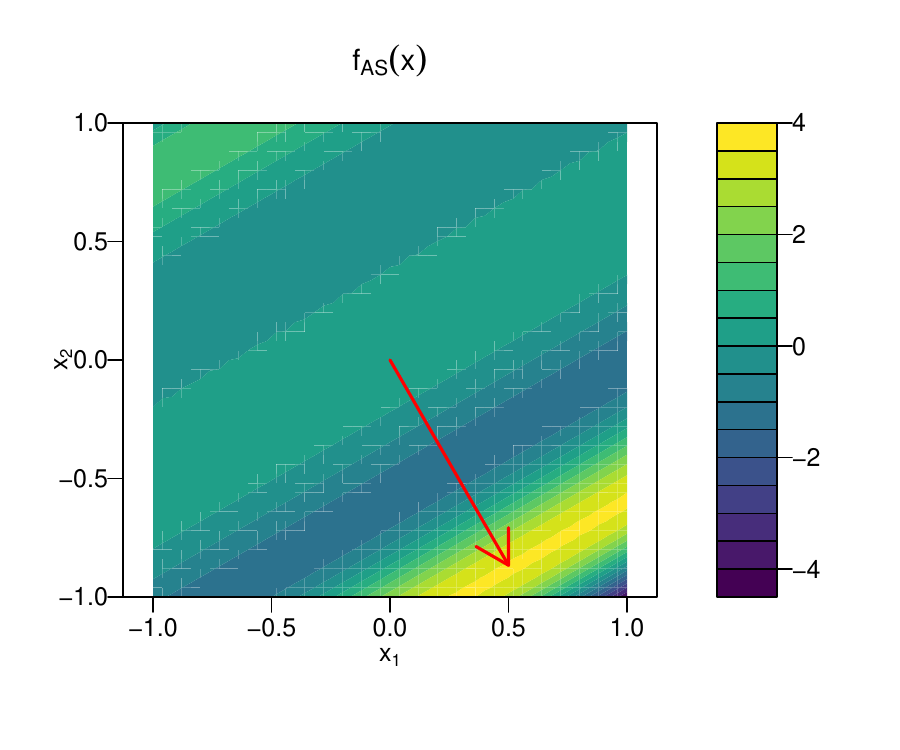}%
\includegraphics[width=0.3333\textwidth, trim= 20 35 20 40, clip]{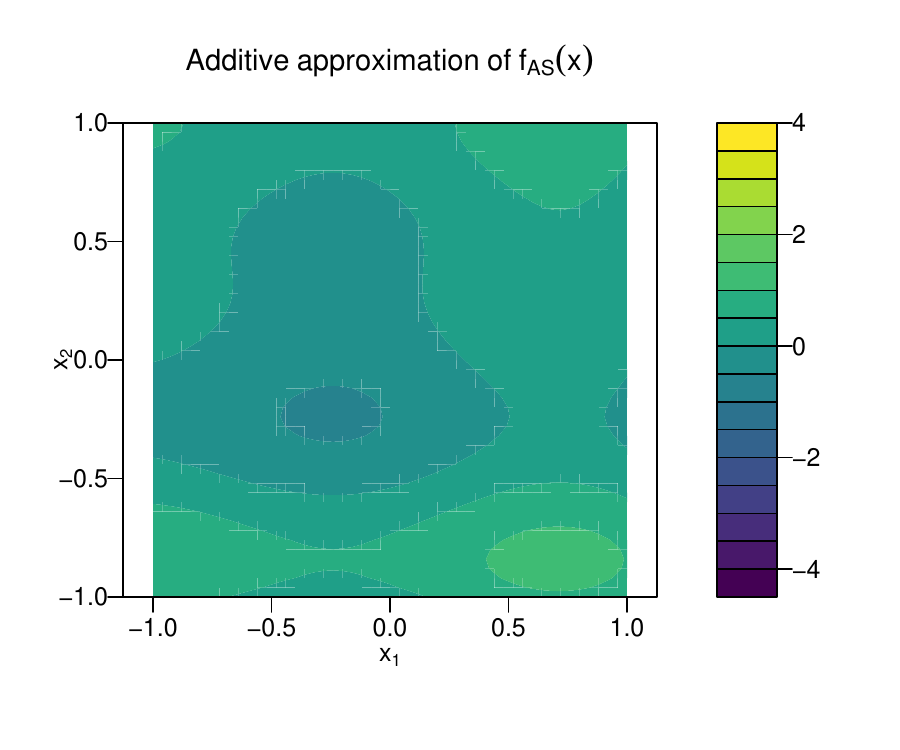}%
\includegraphics[width=0.3333\textwidth, trim= 20 35 20 40, clip]{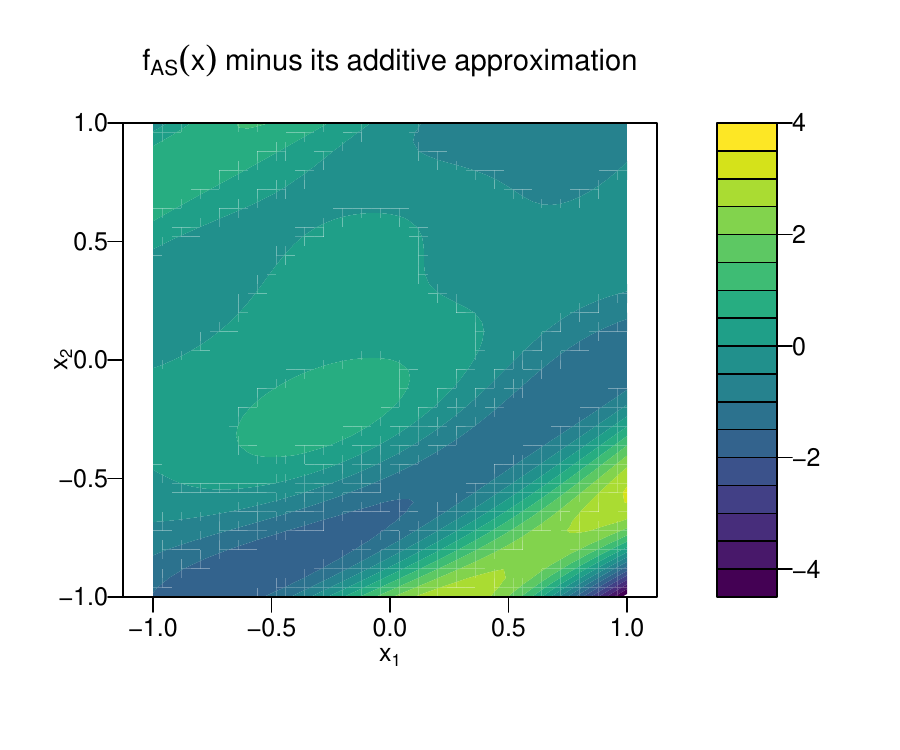}
\caption{Illustration of the overlap between the features learned by additive
 and AS models. Top: approximation of an additive function by a 1d active
 subspace. Bottom: approximation of a function with a 1d active subspace by
 an additive function. The red arrow indicates the active direction.}
\label{fig:asadd}
\end{figure*}

Our goal is to combine the advantages of both additive and linear embedding
models, without further complexifying inference. Ideally, each component
would capture distinct features of the high-dimensional black-box, thereby
enhancing the overall model performance. But this is not so simple in
practice. We illustrate this with Figure \ref{fig:asadd}. For the additive
function, the active subspace estimation gives $\C \approx \begin{bmatrix}
20672 & -5268\\
-5268 & 40684,
\end{bmatrix}$ whose eigen-vectors are not axis-aligned. Consequently, the
 remainder of the additive function after its AS approximation is not
 additive anymore. Conversely, the additive approximation of a ridge function
 leads to a remainder that is no longer ridge.

\subsection{Combination Options}

There are presumably many options to combine models. A straightforward idea
would be to simply sum the two types of models, but this raises
identifiability issues. Without the AS assumption, this is the model proposed
by \cite{plate1999accuracy} for visualization, by gradually modifying the
degree of additivity. An additive model has also been constructed on
directions learned by projection pursuit as in \cite{Gilboa2013} or directly
on random directions as in \cite{Delbridge2020}. Another such idea is first
to apply a rotation with AS as a preprocessing step, e.g., as
in \cite{wycoff2022sensitivity}, before applying an additive model on the
inactive directions. The drawbacks here are the loss of interpretability of
the additive model in the original variables and a potential lack of
interpolation if the additivity assumption does not hold. The converse is to
learn a linear embedding directly on the residuals of an additive model,
which results in a challenging inference problem if done in one step. One
workaround to include orthogonality conditions would be to
follow \cite{lenz:hal-01063741,Ginsbourger2016} with an orthogonal
decomposition between additive and ortho-additive components, before applying
AS on the latter. While appealing, this decomposition remains based on the
projection of a single high-dimensional tensor product kernel. Furthermore,
the independent integral derivations of ortho-additive and AS components do
not seamlessly carry over when combined.

\subsection{A Multi-fidelity Approach}

To maintain an orthogonality condition for identifiability, we propose to rely
on the one enjoyed by the multi-fidelity model. Indeed, this model regresses
the coarse model when no fine data is available, in order to improve
extrapolation--the predominant prediction regime characteristic of
high-dimensional problems. We opt for a first-order additive model as the
coarse level and a linear embedding model as the finer one. This choice is
the most natural since (first-order) additivity is a restrictive yet
data-efficient assumption. The linear embedding then is able to learn the
remaining high-orders of interaction, allowing flexibility in the choice of
the embedding dimension $r$. In the case $r = d$, a regular GP is fit on the
residuals between the additive model prediction and the data, still on the
rotated initial input space, thus maintaining interpolation in the noiseless
case. This rotation is shown to be helpful as a pre-processing
step \cite{wycoff2022sensitivity}. If the additive model is a good
approximation, the remaining variance of the GP on the residuals will be
smaller. This may thus alleviate the inflation of the variance in high
dimension stemming from the behavior of distances, causing over-exploration
as often observed in high-dimensional BO, see, e.g., \cite{Eriksson2019}. The
converse, taking an AS model as the coarse model is less relevant as it will
also capture additive components while a finer additive model may not be able
to interpolate deterministic data. Plus, again, the additive model would be
on a rotated input space, losing interpretability and convenience. Our
proposed model thus writes:
\begin{equation}
\left\{
    \begin{array}{ll}
       Y_E(\x) = \rho Y_C(\x) + \delta(\A \x)\\
       Y_C(\x) \perp \delta(\A \x)
    \end{array}
\right.
\label{eq:mainmodel}
\end{equation}

\begin{figure*}[htpb]
\centering
\includegraphics[width=0.45\textwidth, trim= 200 270 20 300, clip]{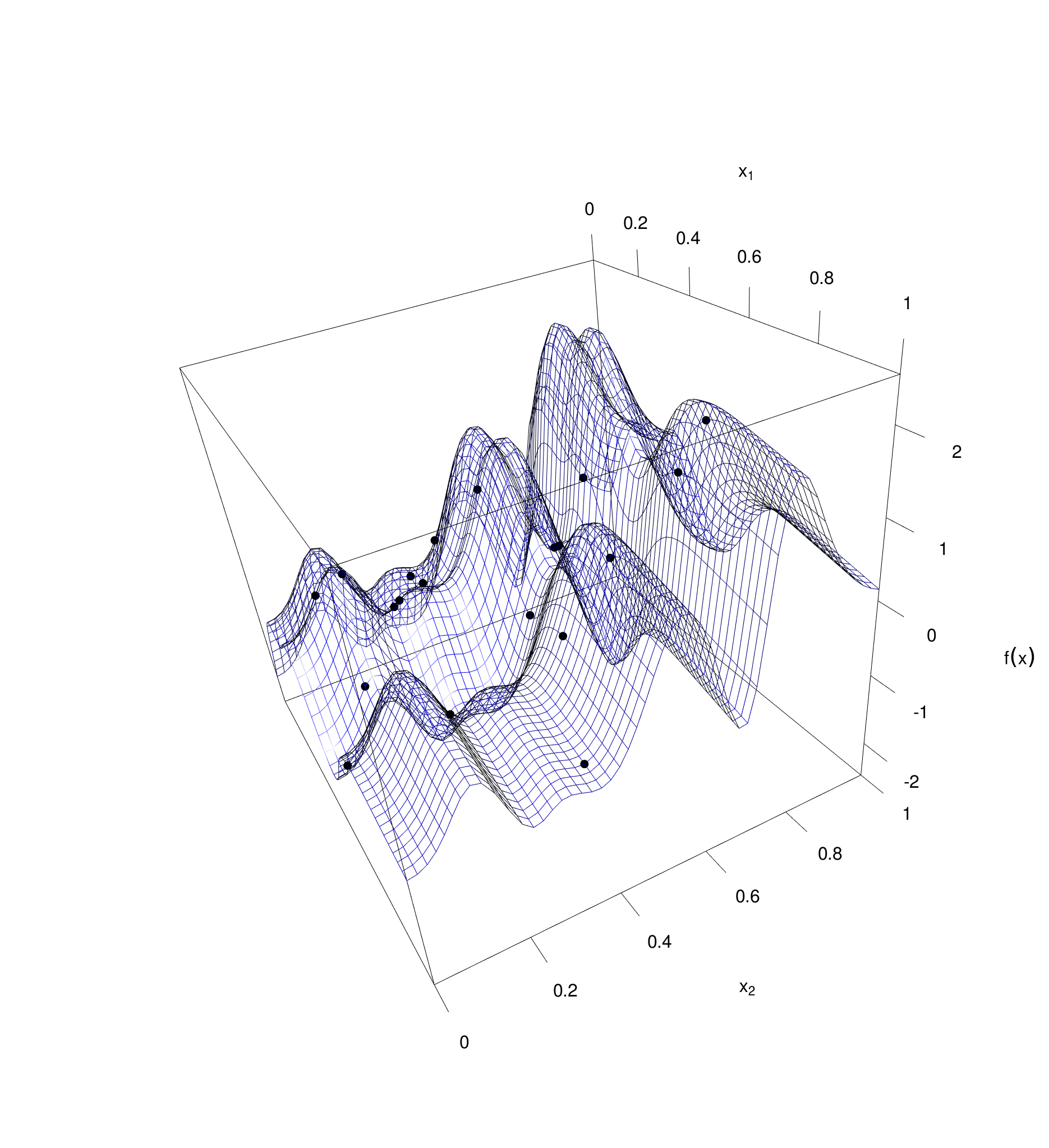}%
\includegraphics[width=0.45\textwidth, trim= 200 270 20 300, clip]{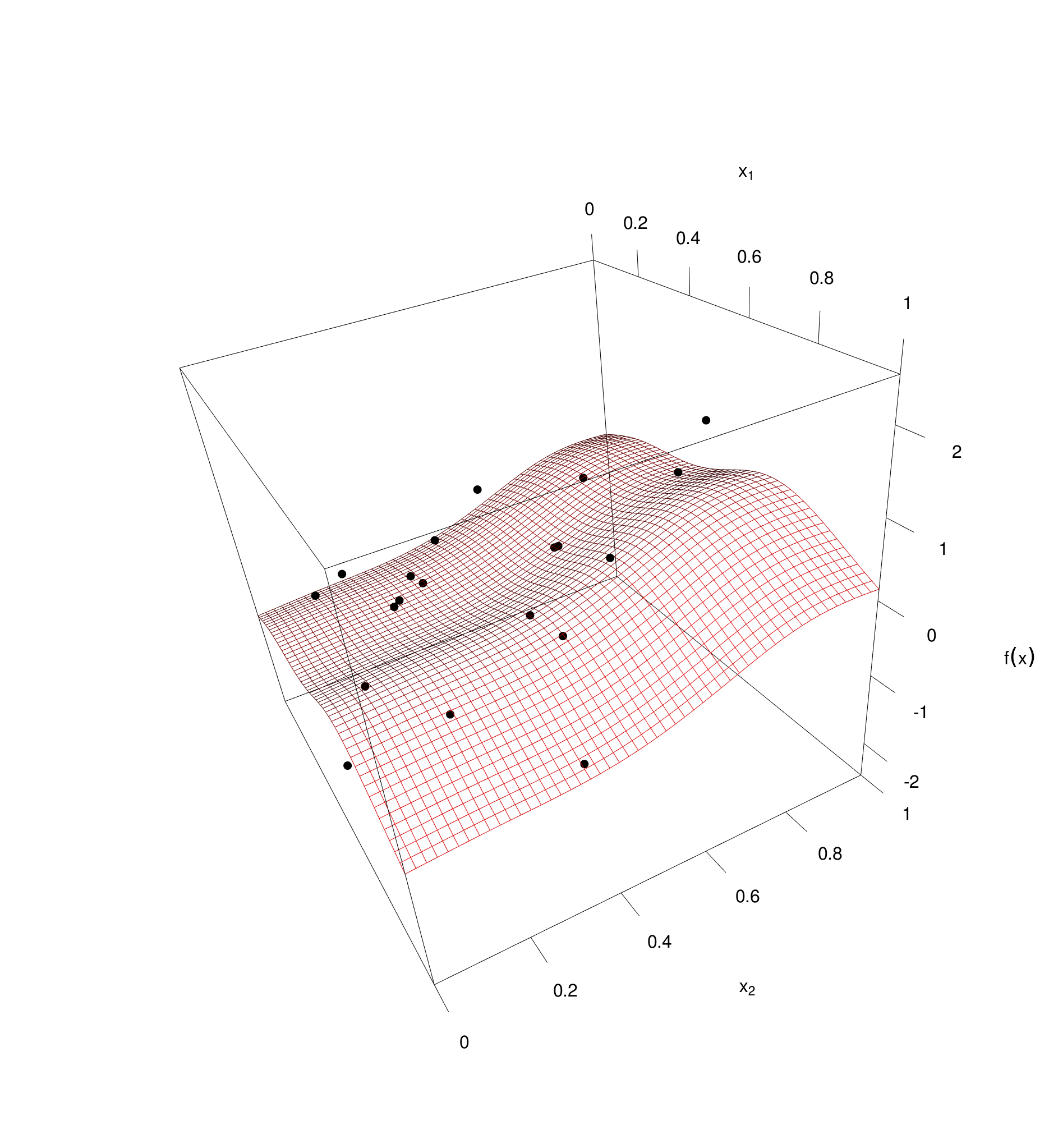}
\caption{Additive model mean prediction surfaces on 20 points from the Branin
function, interpolating (left) or approximating (right).}
\label{fig:addBran}
\end{figure*}

In practice, there is solely one set of (fine level) evaluations $\y$. To
obtain coarse level values and $\mathbf{d}$, we simply take the predictions
provided by the additive model: $\y^{(C)} = m_n^{(C)}(\XC)$. An issue that
may arise with this scheme is when an additive model fits the data perfectly.
This scenario is more likely to occur when the training dataset is small, or
the dimension is large. Figure \ref{fig:addBran} illustrates an example using
the non-additive Branin function. There are several options to cope with this
issue, which can be detected by comparing the noise variance to the process
variance. One is to restrict the range of the lengthscale (e.g., using prior
knowledge). Another one is to withhold some designs from the additive model.
That way, if the prediction by the additive model is not accurate, then it
will be corrected at the fine level. Lastly, it remains to build the finer
model on the residuals of the prediction by the coarse level. When fitting a
linear embedding model, one key question is the selection of the dimension.
We choose to rely on the likelihood to do so.

\subsection{Which AS matrix in the MuFi model?}

Concerning the fit on the residuals, it remains to be specified how the active
subspace is actually computed. Ideally, the corresponding AS matrix could be
learned when fitting the hyperparameters of the multi-fidelity GP model,
following, for instance, Appendix \ref{ap:asgp}. But this is rather
computationally intensive, since this GP fit is on $Y_d = Y_E - \rho Y_C$,
hence it depends on $\rho$. By developing the corresponding active subspace
matrix, $\C_{\mathbf{d}} := \int \limits_\Xset \nabla(f_E(\x) - \rho
f_C (\x))^\top \nabla(f_E(\x) - \rho f_C(\x))$, we get $\C_{\mathbf {d}}
= \C_E + \rho^2 \C_C - \rho (\C_{E,C} + \C_{C,E})$, where $\C_E$ and $\C_C$
are the AS matrices of $f_E$ and $f_C$ while $\C_{E,C}$ and $\C_ {C,E}$ are
their cross-term counterparts. Hence, it is possible to compute $\C_E$,
$\C_C$, $\C_{CE}$ and $\C_{EC} = \C_{EC}^\top$ only once.

\subsection{Proposed Instantiation}

\begin{algorithm}[tb]
\caption{Pseudo-code for the multi-fidelity high-dimensional GP}\label{alg:hdGP}
\begin{algorithmic}[1]
\State {\bfseries Input:}  $\XE = \XC$, $\y$, $p$ (e.g., $0.8$)
\State Train an additive model $Y_C$ on $(\XC, \y)$.
\If{$\tau^2_C \leq 0.01 \times \sum \limits_1^d \alpha_i$} 
  \State Sample $n_0 = p \times n$ data points from $\x_{1:n}$, $\y$ and
  remove the rest from $\XC$ and $\y^{(C)}$.
  \State Train an additive model $Y_C$ on $(\XC, \y^{(C)})$.
\EndIf
\State Compute the AS matrix of the additive model, $\C_C$
\State Train a regular GP $Y_E$ on $(\XE, \y)$ to compute its AS matrix $\C_E$.
\State Compute the cross AS matrices $\C_{E,C}$ and $\C_{C, E}$.
\State Predict the response of $Y_C$ at $\XE$: $m_n^{(C)}(\XE)$.
\State Train an AS multi-fidelity GP from the residual data: $\mathbf{d} = \y - \rho m_n^{(C)}(\XE)$, varying the number of dimensions kept $r$.  
\State {\bfseries Output:} Trained multi-fidelity model.
\end{algorithmic}
\end{algorithm}

 Following the modularization principle \cite{liu2009modularization}, that is,
 separating inference of different modules, we chose to perform a two-stage
 approach, rather than the single-stage one described in
 Appendix \ref {ap:asgp}. We detail the construction of this
 model (\ref{eq:mainmodel}) in Algorithm \ref{alg:hdGP}. In Step 2, we
 restrict ourselves to using a first-order additive model, avoiding
 higher-order terms selection. From Steps 3 to 6, in case the noise variance
 of the additive model $\tau^2_C$ is less than one percent of the additive
 process variance, $\sigma^2_ {n,C} = \sum \limits_{i=1}^d \alpha_i$, then,
 to avoid overfitting, the additive model is replaced by one trained on a
 fraction $p$ of the data. Note that this corresponds to the opposite of the
 typical nesting in multi-fidelity. In Step 7 to 9, the AS matrices involved
 in $\C_\mathbf{d}$ are pre-computed, for instance based on a Monte Carlo
 estimation. For this, we prefer the active subspace to be learned (and not
 random, which usually requires several random AS to work well, in practice
 and theoretically, see \cite{cartis2023bound}). There, we follow the
 approach by \cite {Wycoff2021}, because the required number of
 hyperparameters to learn the linear embedding remains limited compared to,
 say, \cite {Garnett2013,letham2020re}.

 In Step 10, the low fidelity data is obtained by predicting with the additive
 model. The remaining step is dedicated to learning the linear embedding and
 the corresponding GP hyperparameters: AS dimension, lengthscales on the AS,
 and multi-fidelity correlation $\rho$. The eigen-vectors $\U$ of
 $\C_\mathbf{d} = \U \boldsymbol {\Lambda} \U^\top$ provide the new
 coordinate system, i.e., with rotated coordinates $\mathbf{X}_{E,r}
 = \XE \U_{1:r}$ (assuming that $\Xset$ is centered). It remains to select
 the number $r$ of eigen-vectors from $\U$. Since $r$ is a discrete
 parameter, one simple workaround is to optimize the lengthscale parameters
 for various values of $r$, before selecting the value with the best overall
 log-likelihood. One can consider that lengthscales for the inactive
 dimensions are set to infinity, such that this remains the same model
 defined on the full rotated space. Note that this can result in a noisy
 model, where the noise absorbs the contributions of the inactive
 dimensions.

\section{Empirical Evaluation} 

We conduct a comparative analysis between the multi-fidelity approach and
baseline methods on synthetic functions and datasets. Rather than taking very
large input dimensions $d$ and data sizes $n$, we focus on the low data
regime, considering $n$ up to $500$ and $d$ up to $32$.

\subsection{Setup}

As a baseline, we use a standard GP model (hereafter denoted by Ref), with a
tensor product kernel. The \texttt{R} \cite{rcore}
package \texttt {hetGP} \cite{Binois2021} is used for learning the
hyperparameters, where the initialization of the hyperparameters is
complemented by an initialization with the \texttt{R}
package \texttt{RobustGaSP} \cite {Gu2022} for a robust hyperparameter
estimation \cite{Gu2018}. We entertain an additional variant of standard GPs,
with an isotropic kernel (Iso). Additionally, we assess the individual
component models of the multi-fidelity approach: a first-order additive
model (Add) and a linearly embedded GP (AS). A multi-fidelity model with a
standard GP for the finer level is also considered (MF), in addition to the
version with an active subspace (ASMF, which is our main proposal). We
further add naive variants (n-) of the multi-fidelity models, involving a
direct summation of the additive model with the one on the residuals. The
implementation of the proposed models is in the Supplementary Material to
allow reproduction of the results. All use Matérn 5/2 kernels in these
experiments.

For test functions, we start with draws from GPs with $d = 8, 15$, avoiding
model mismatch. That is, we consider draws from standard GPs, GPs on active
subspaces (ASGP) and first-order additive GPs (AddGP). A weighted sum of
AddGP and ASGP with varying weights is added as well. The subsequent set of
tests is based on real
datasets \texttt{Concrete} ($d=8$) \cite {Newman1998}, \texttt{pumadyn} ($d=8,32$)
\cite{corke1996robotics} and
\texttt{CASP} ($d=9$) \cite{nr}.

Lastly, we use classical toy problems:
Sobol ($d=8$) \cite {marrel2009calculations},
penicillin ($d=7$) \cite{liang2021scalable}, Levy ($d=10,
20$) \cite{laguna2005experimental} and Cola ($d=17$) \cite {mathar1994class}.
We also embed lower dimensional test functions, Hartmann3 ($d_e = 3$) with a
random AS matrix with $d = 15$ and Branin ($d_e=2$) \cite{dixon1978global}
with a random hashing matrix with $d = 10$. We complement these by adding an
additive GP realization to the linearly embedded Hartmann3 function. 

From a thousand randomly sampled locations (unless the data set size is
smaller) where these benchmarks are evaluated, a training set of size $n_0$
is extracted. All test sets are centered, and rescaled to unit variance. As
for metrics, we rely on the root mean square error (RMSE) and a proper
score \citep{gneiting2007strictly} (i.e., log-predictive density). They are
computed on the remaining data after training.

\subsection{Results}

\begin{figure*}[htpb]
\centering
\includegraphics[width=0.3333\textwidth, trim= 30 45 30 40, clip]{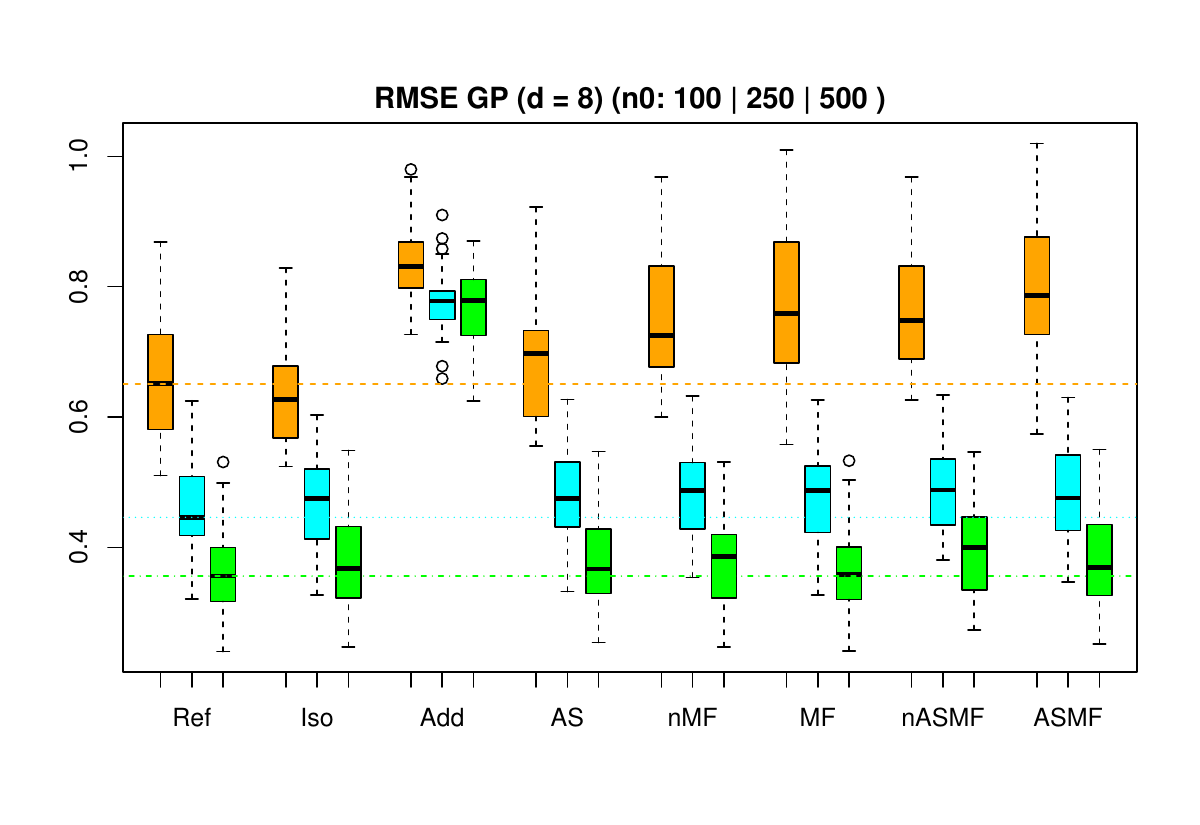}%
\includegraphics[width=0.3333\textwidth, trim= 20 45 30 40, clip]{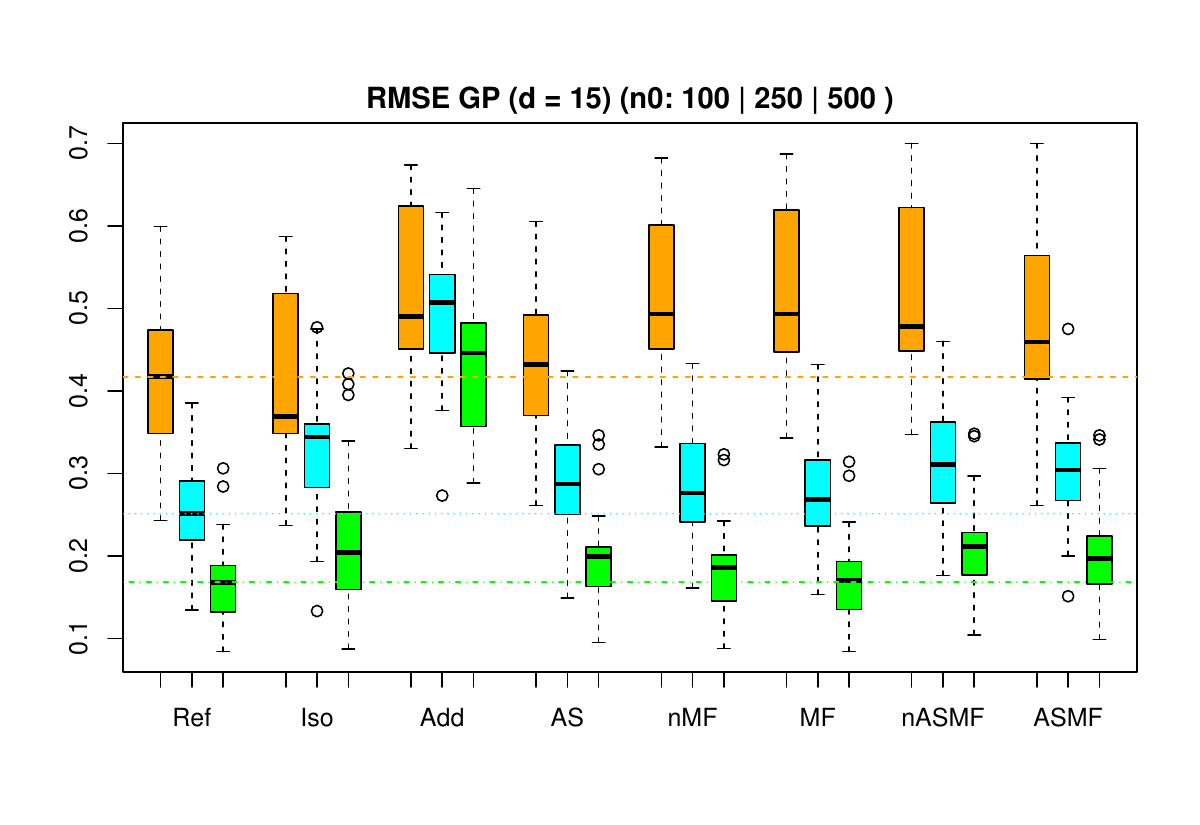}%
\includegraphics[width=0.3333\textwidth, trim= 20 45 30 40, clip]{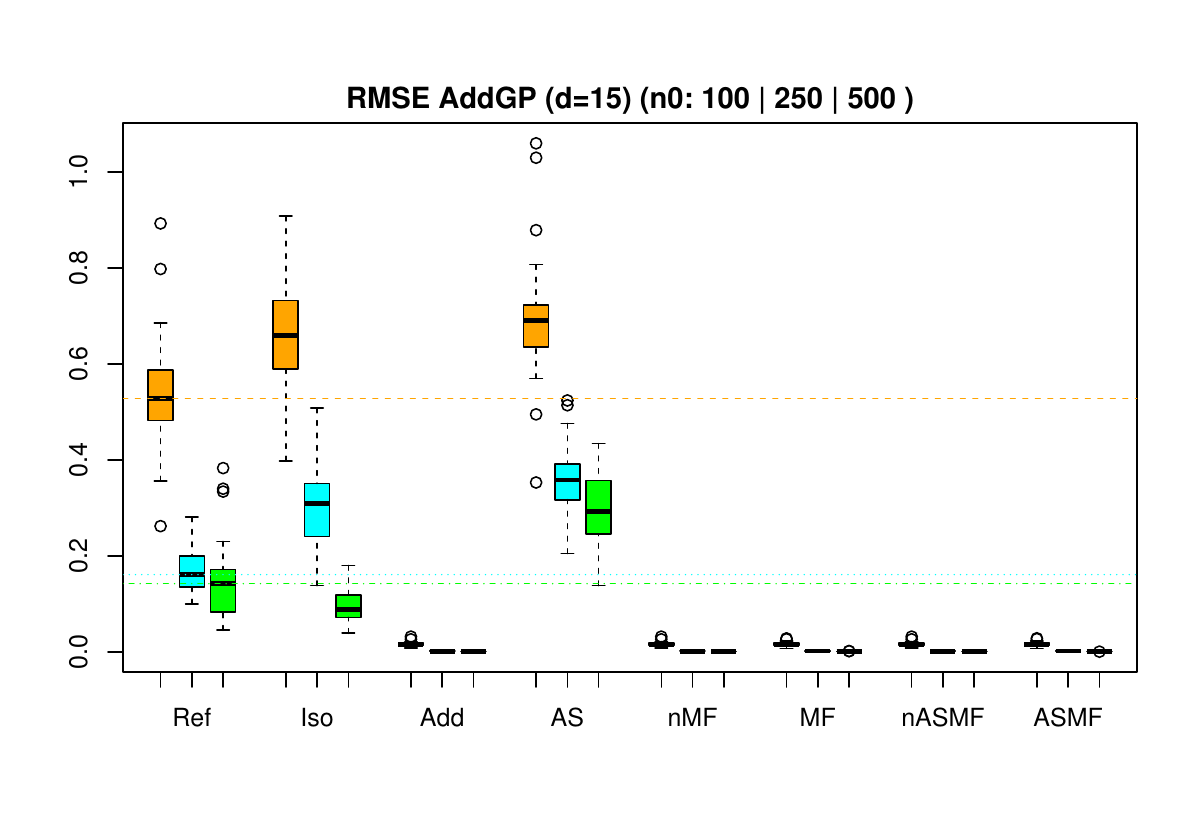}\\
\includegraphics[width=0.3333\textwidth, trim= 30 45 30 40, clip]{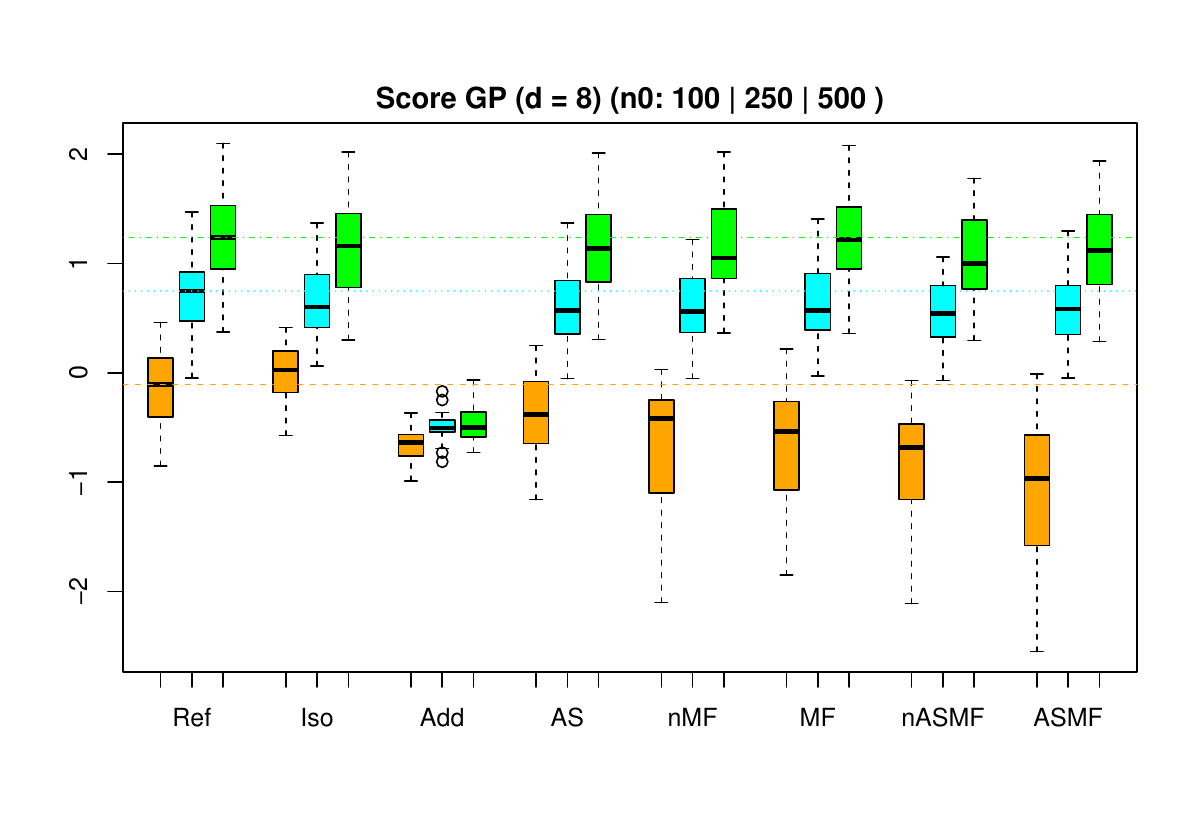}%
\includegraphics[width=0.3333\textwidth, trim= 20 45 30 40, clip]{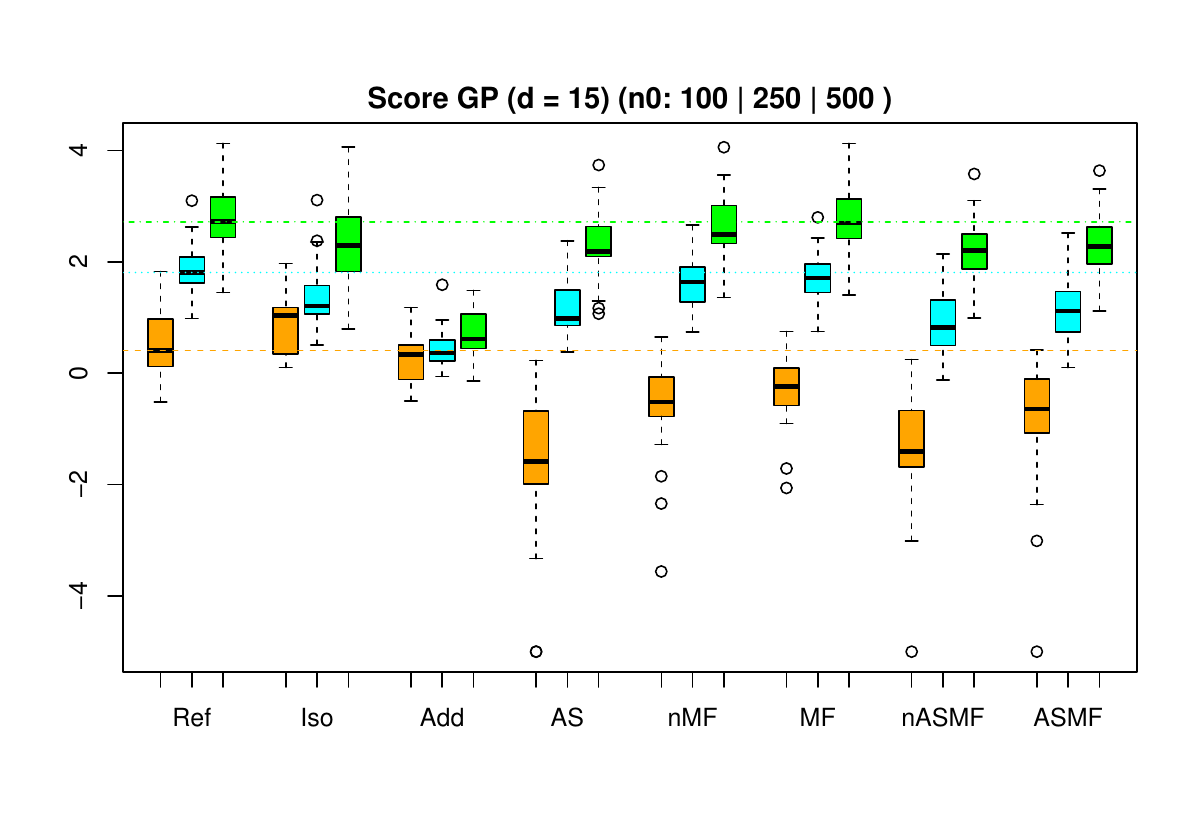}%
\includegraphics[width=0.3333\textwidth, trim= 20 45 30 40, clip]{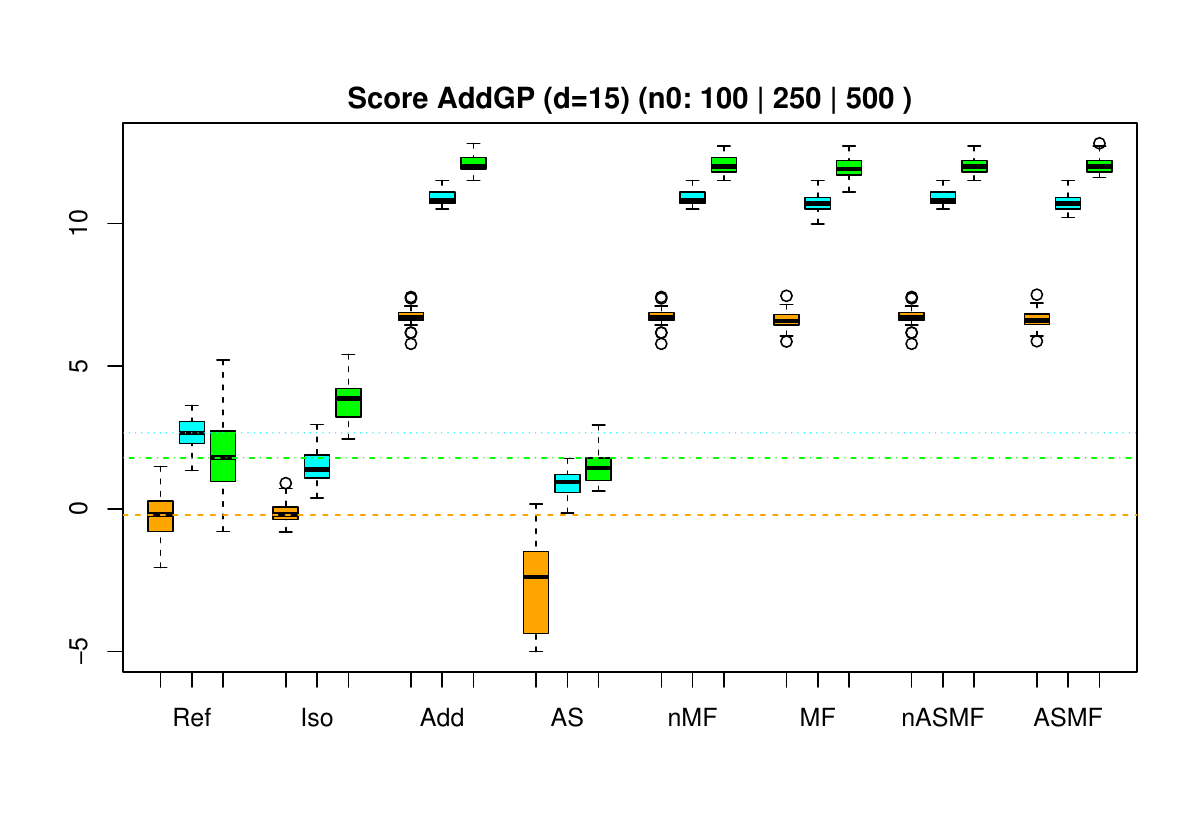}\\
\includegraphics[width=0.3333\textwidth, trim= 30 45 30 40, clip]{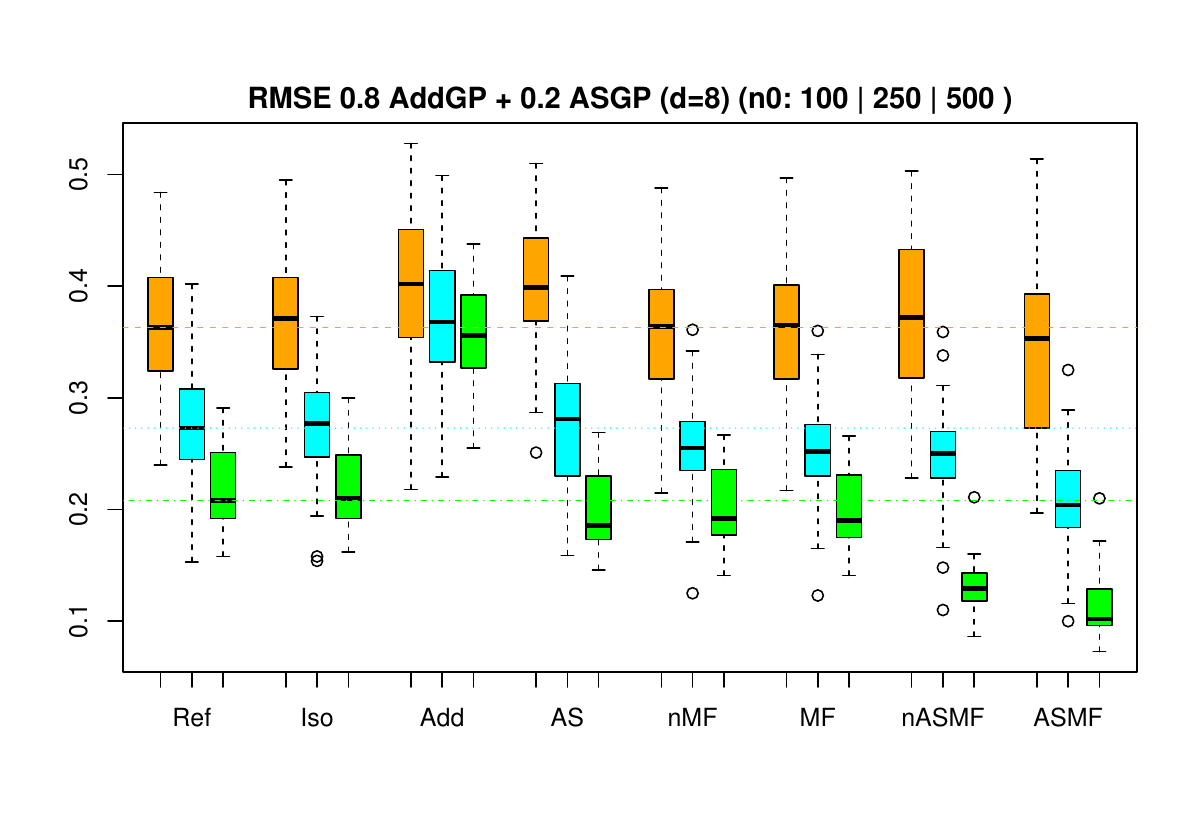}%
\includegraphics[width=0.3333\textwidth, trim= 20 45 30 40, clip]{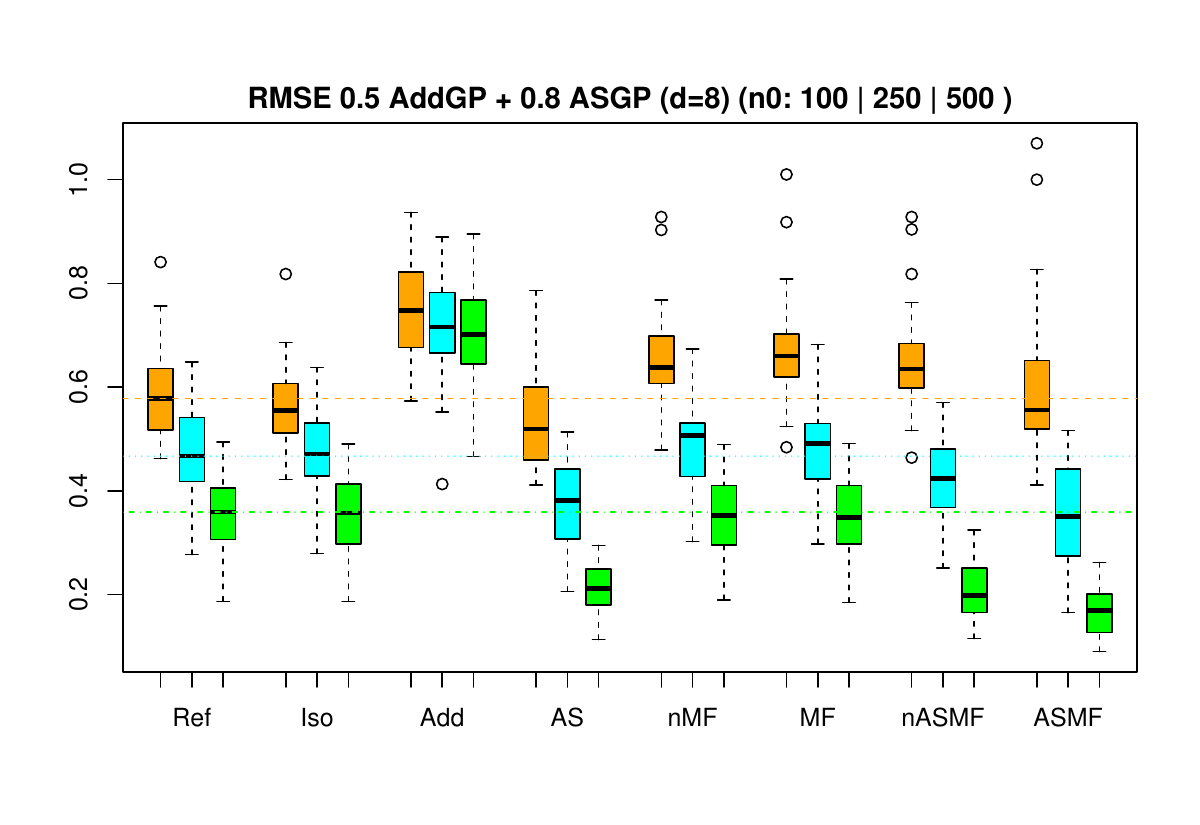}%
\includegraphics[width=0.3333\textwidth, trim= 20 45 30 40, clip]{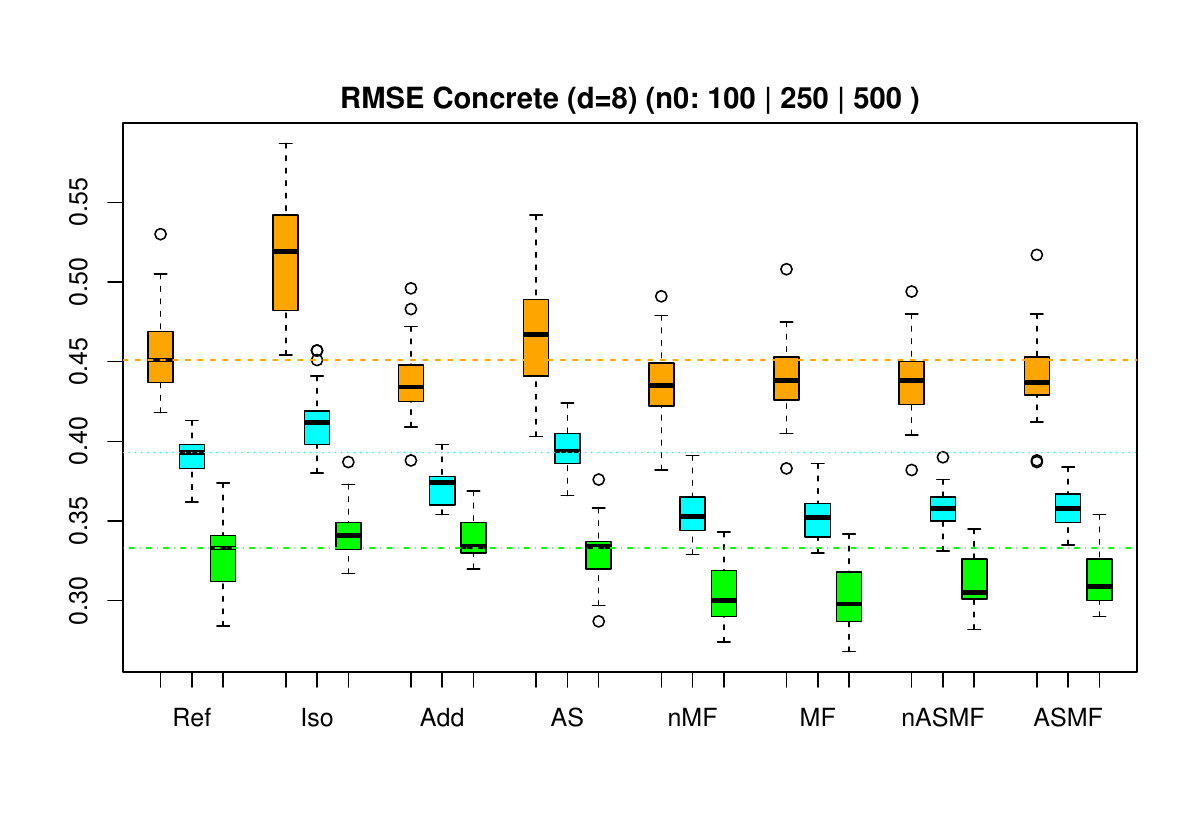}\\
\includegraphics[width=0.3333\textwidth, trim= 30 45 30 40, clip]{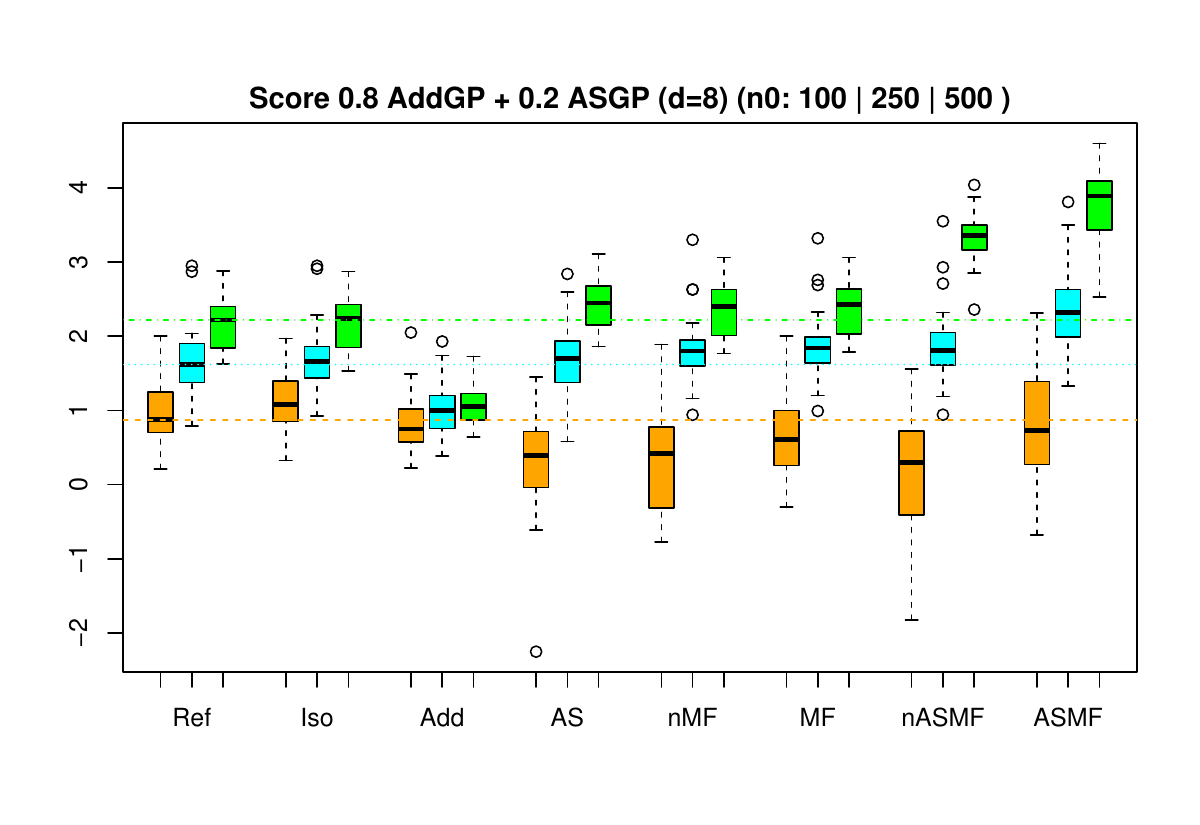}%
\includegraphics[width=0.3333\textwidth, trim= 20 45 30 40, clip]{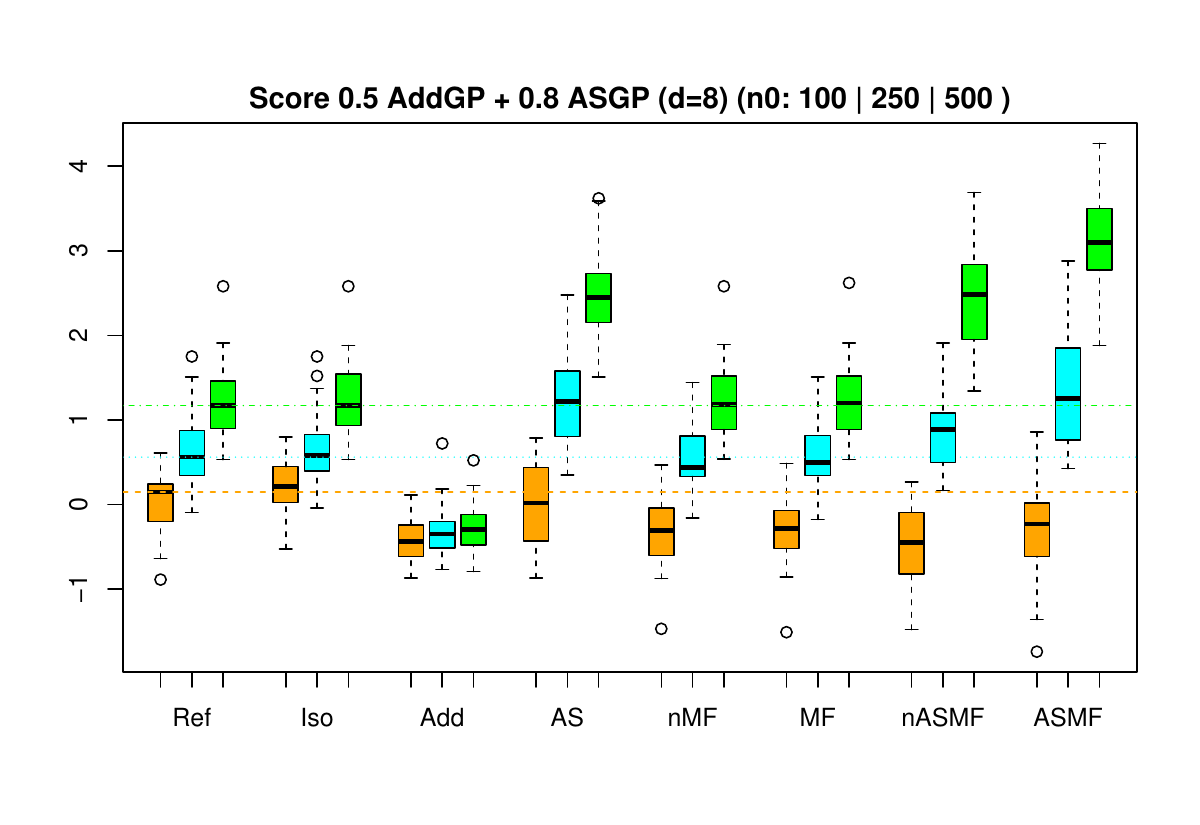}%
\includegraphics[width=0.3333\textwidth, trim= 20 45 30 40, clip]{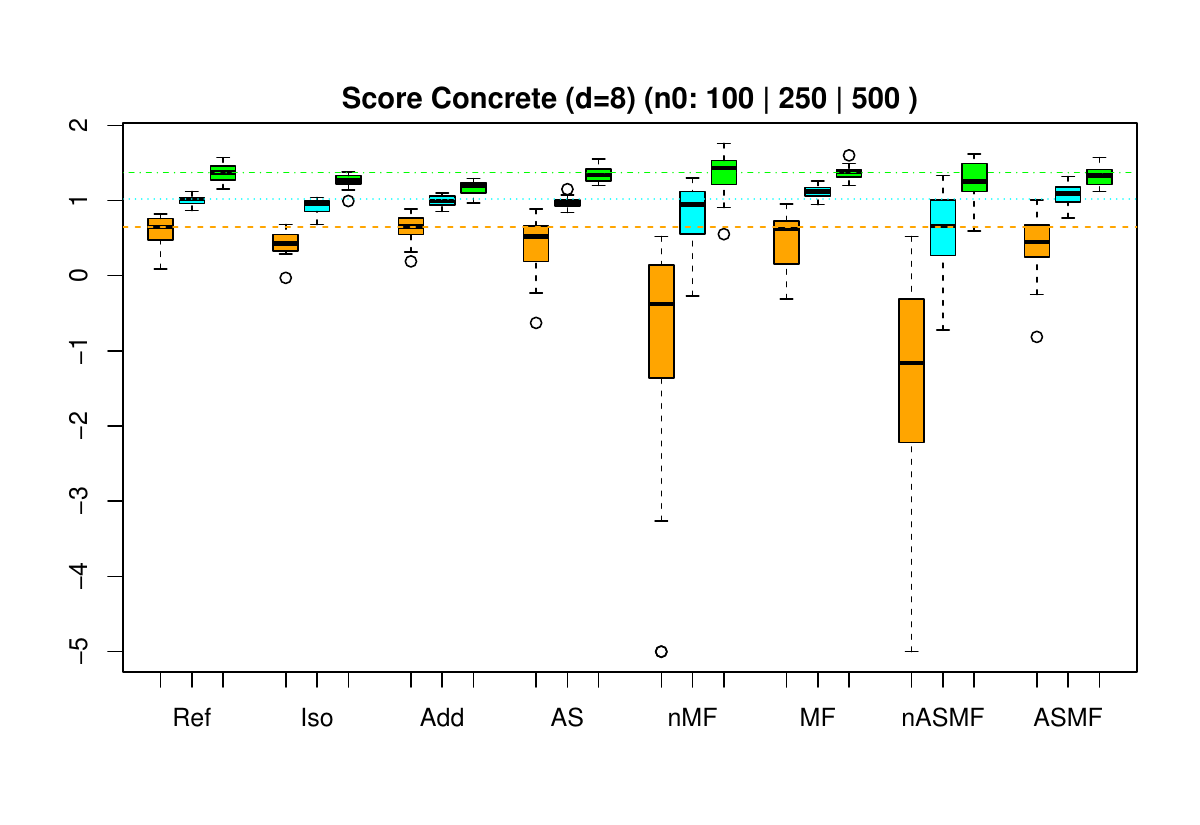}\\
\includegraphics[width=0.3333\textwidth, trim= 30 45 30 40, clip]{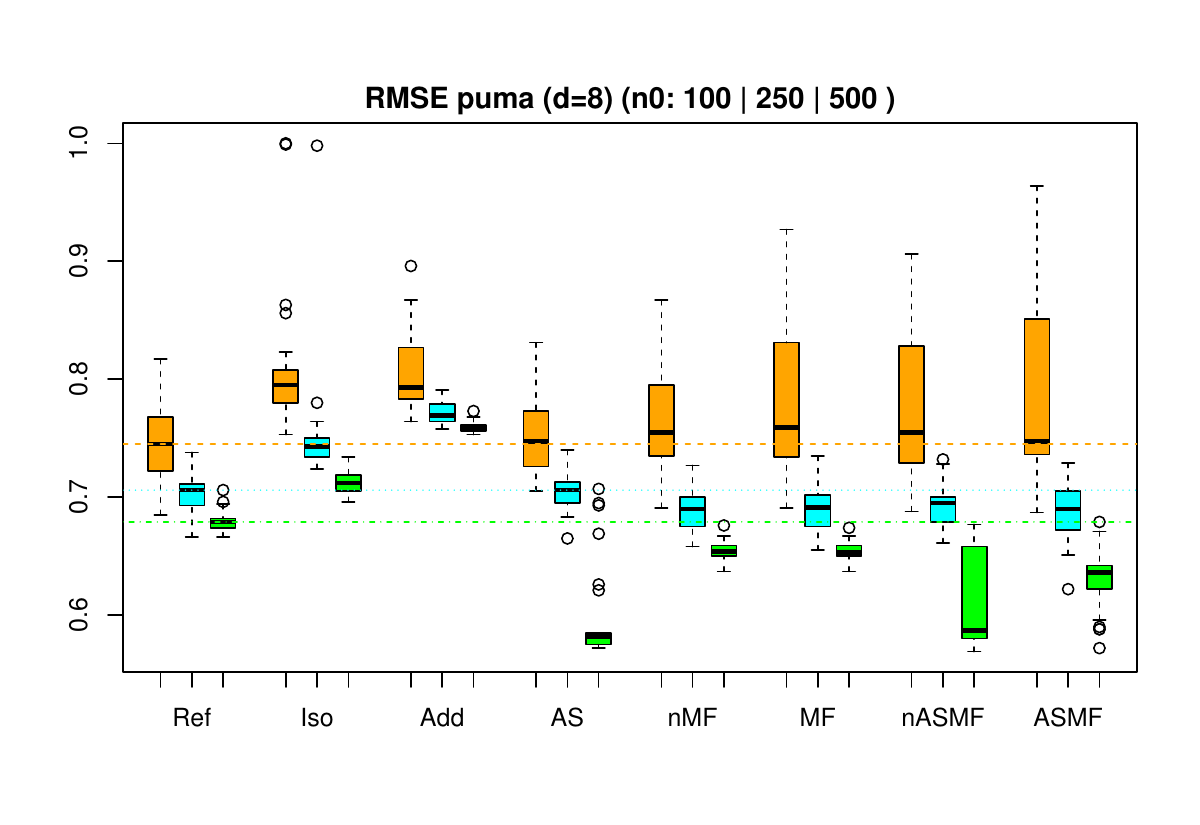}%
\includegraphics[width=0.3333\textwidth, trim= 20 45 30 40, clip]{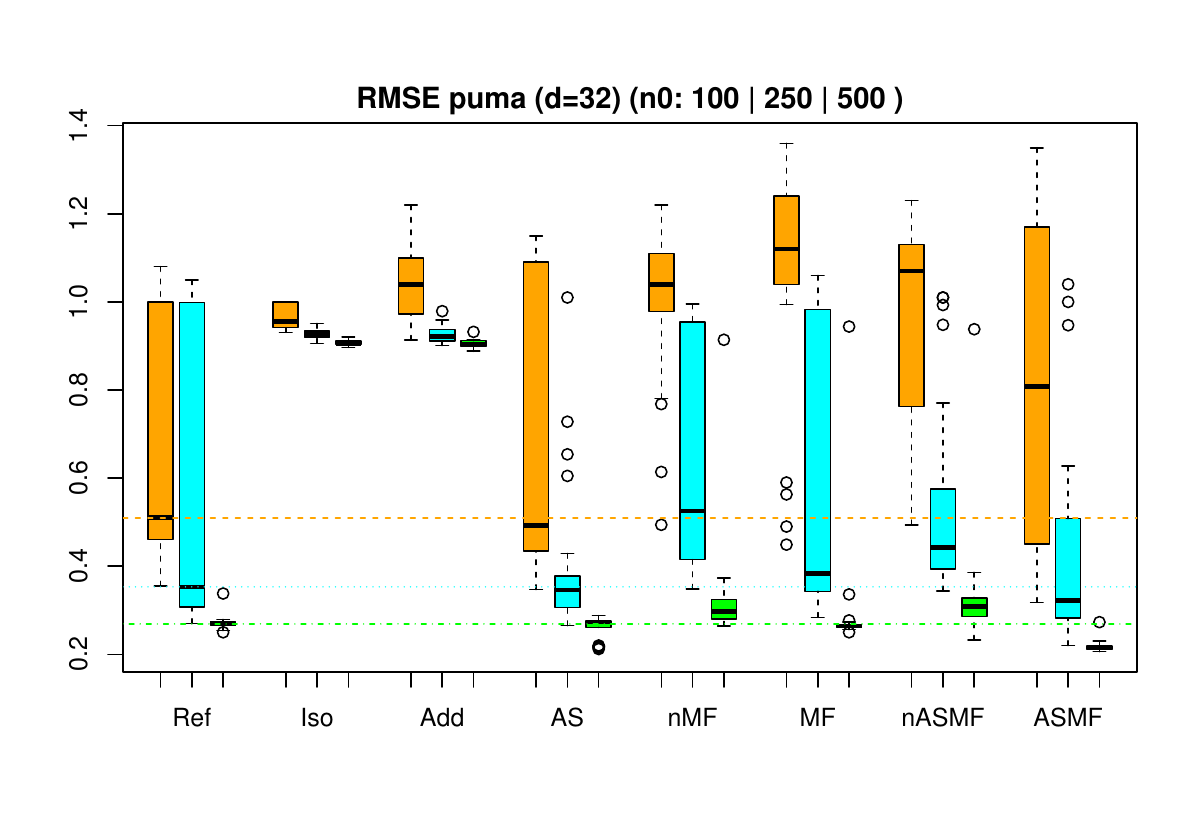}%
\includegraphics[width=0.3333\textwidth, trim= 20 45 30 40, clip]{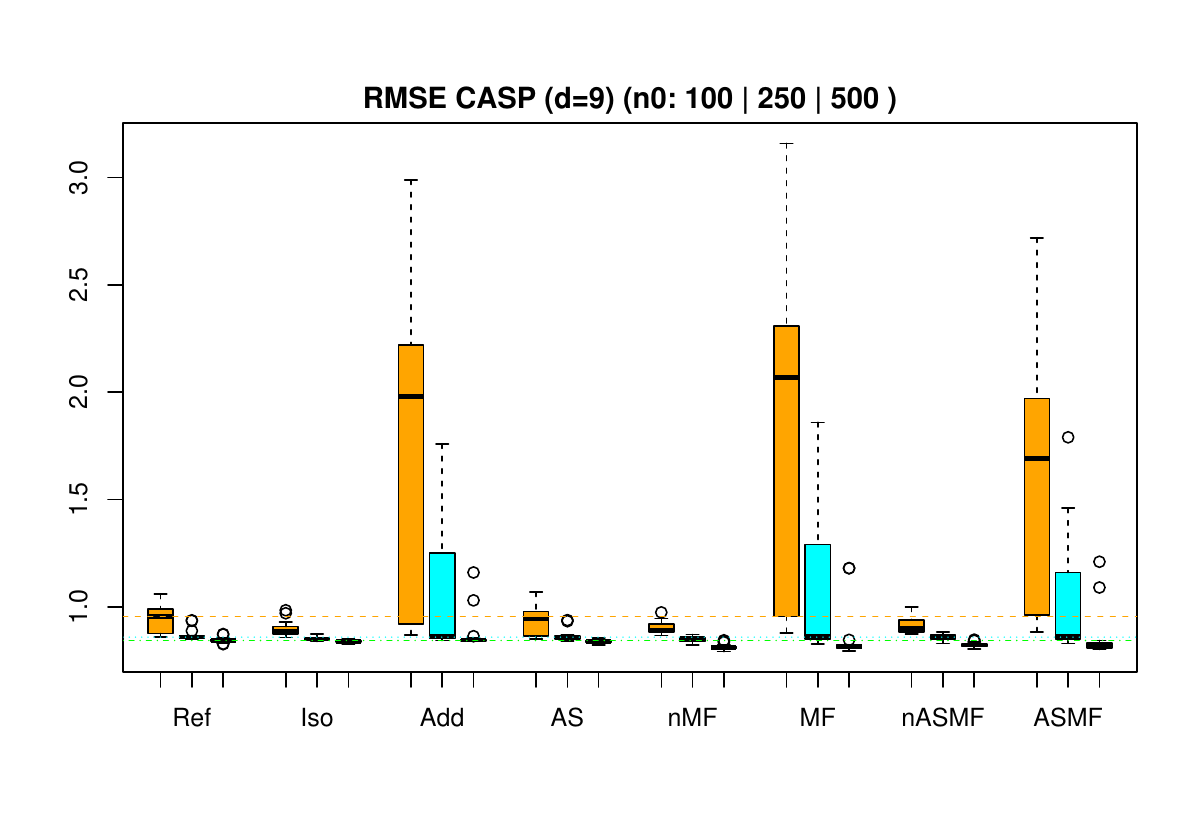}\\
\includegraphics[width=0.3333\textwidth, trim= 30 45 30 40, clip]{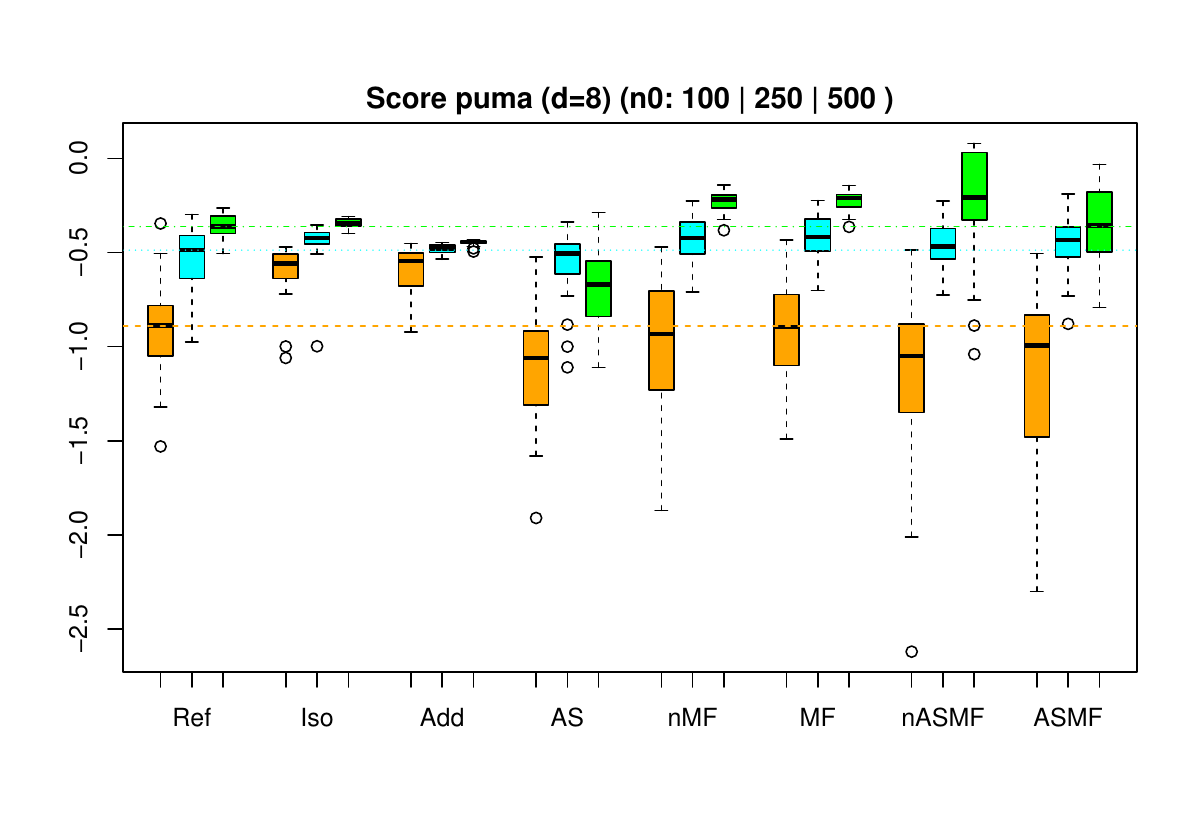}%
\includegraphics[width=0.3333\textwidth, trim= 20 45 30 40, clip]{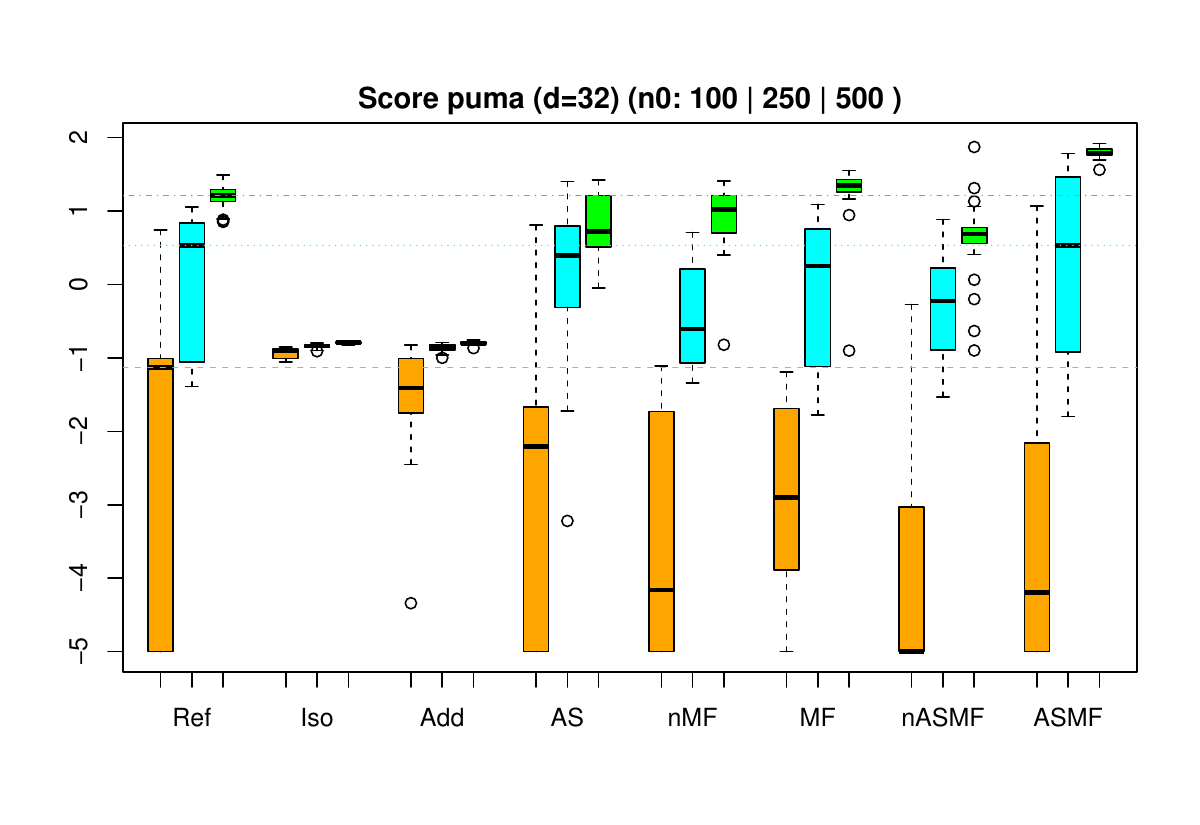}%
\includegraphics[width=0.3333\textwidth, trim= 20 45 30 40, clip]{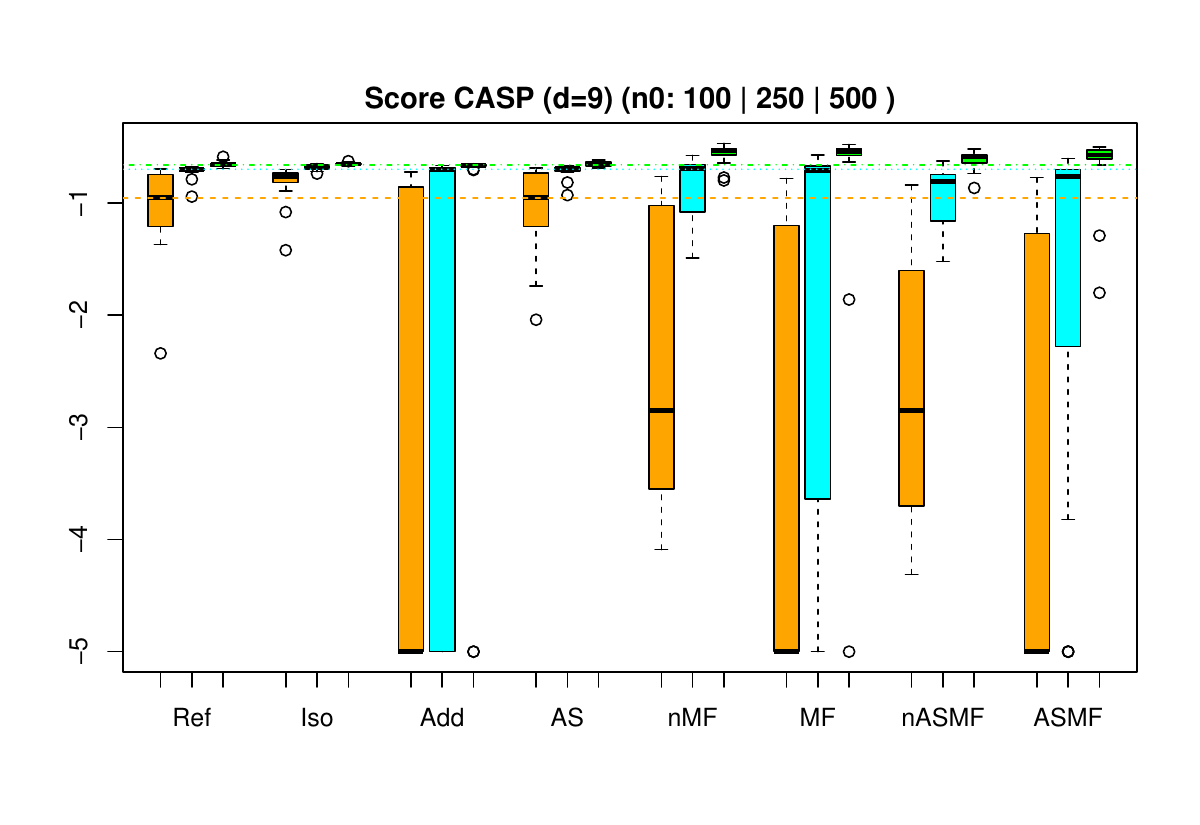}%
\caption{First part RMSE and score results. The different $n_0$ values appear
 with different colors, the corresponding lines indicate the baseline result
 from standard GP models.}
\label{fig:rmseres}
\end{figure*}

\begin{figure*}[htpb]
\centering
\includegraphics[width=0.3333\textwidth, trim= 30 45 30 40, clip]{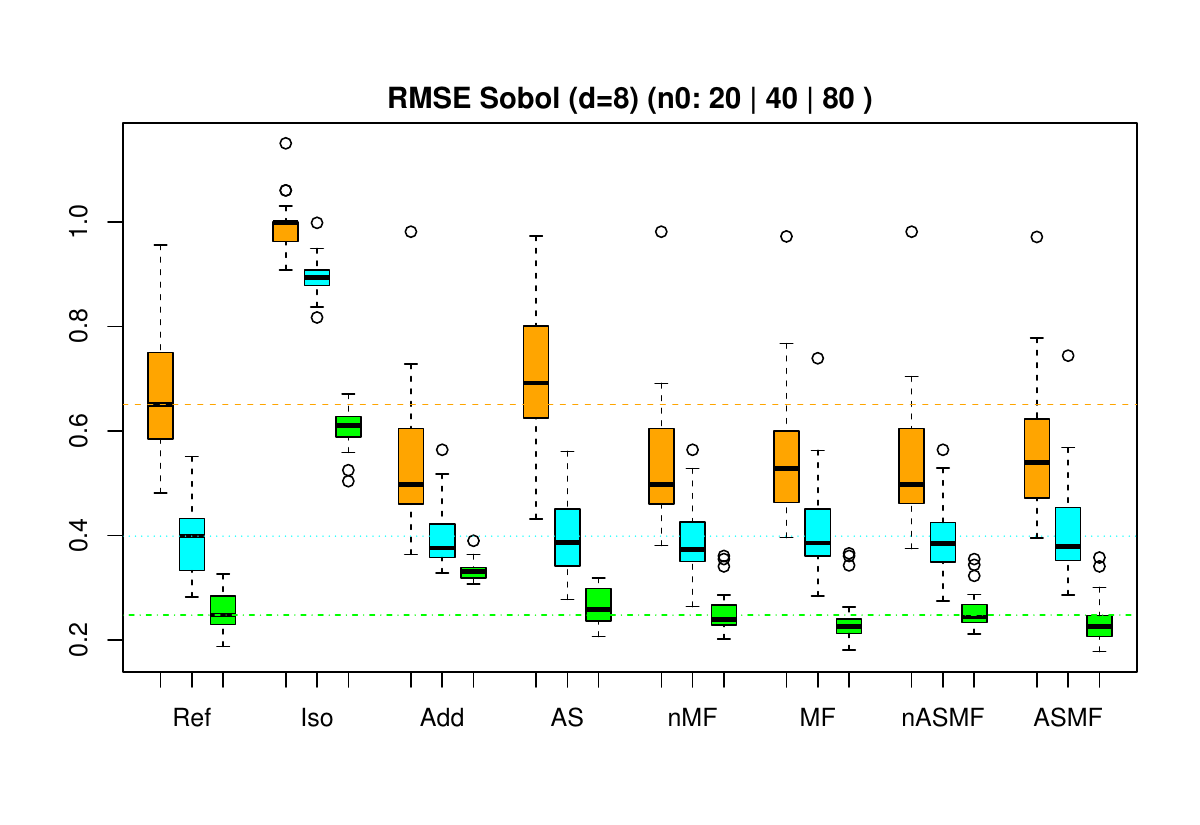}%
\includegraphics[width=0.3333\textwidth, trim= 20 45 30 40, clip]{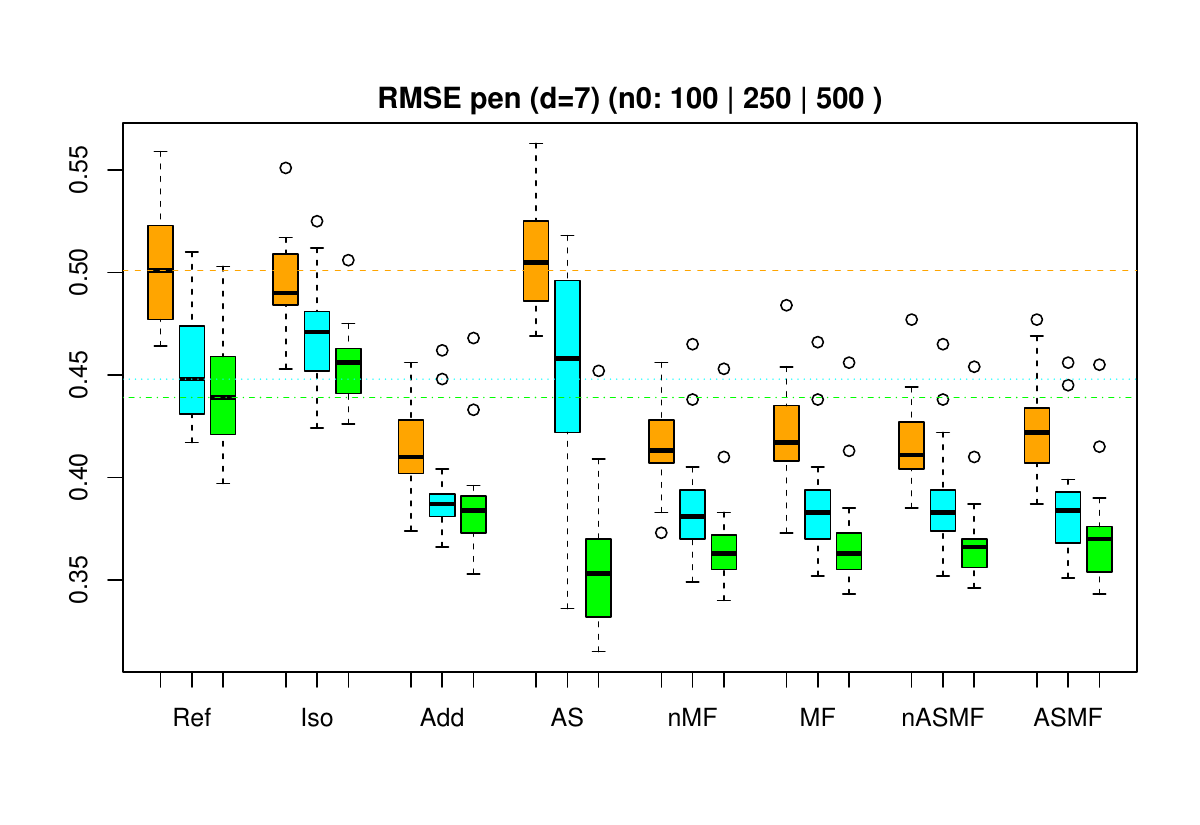}%
\includegraphics[width=0.3333\textwidth, trim= 20 45 30 40, clip]{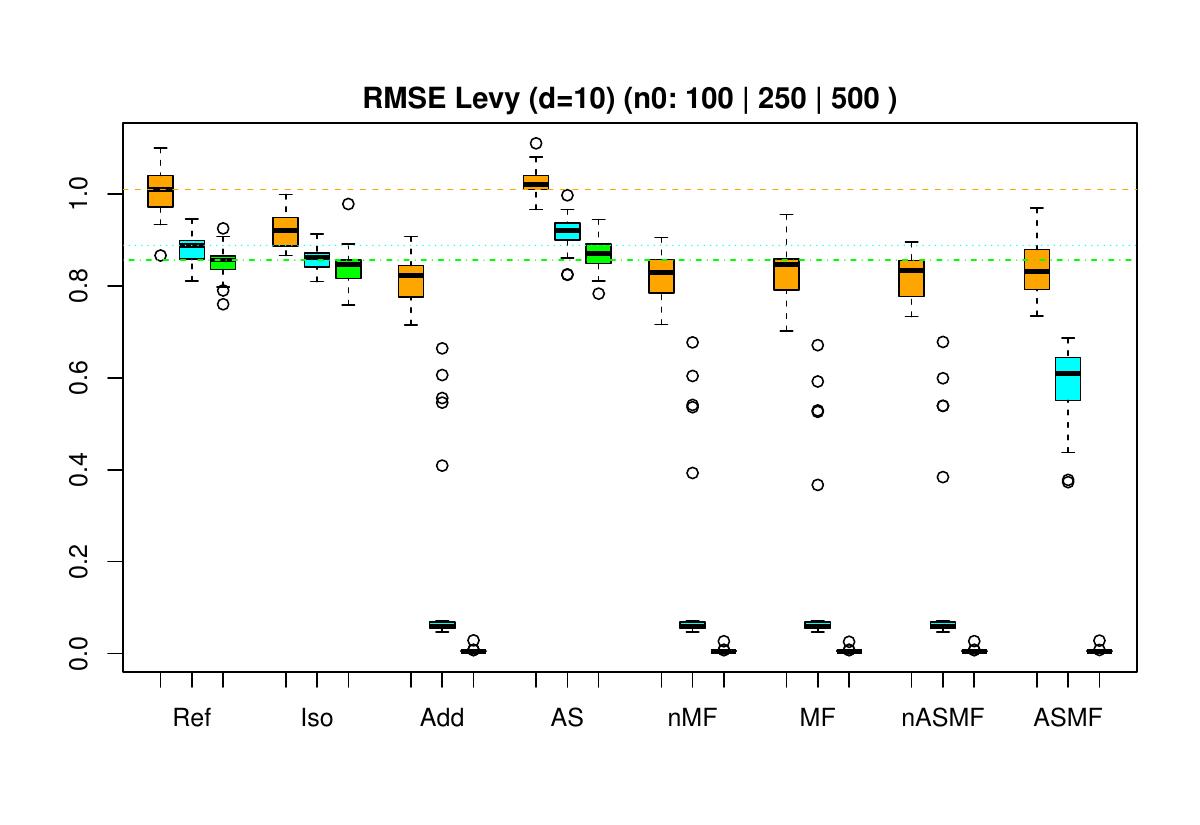}\\
\includegraphics[width=0.3333\textwidth, trim= 30 45 30 40, clip]{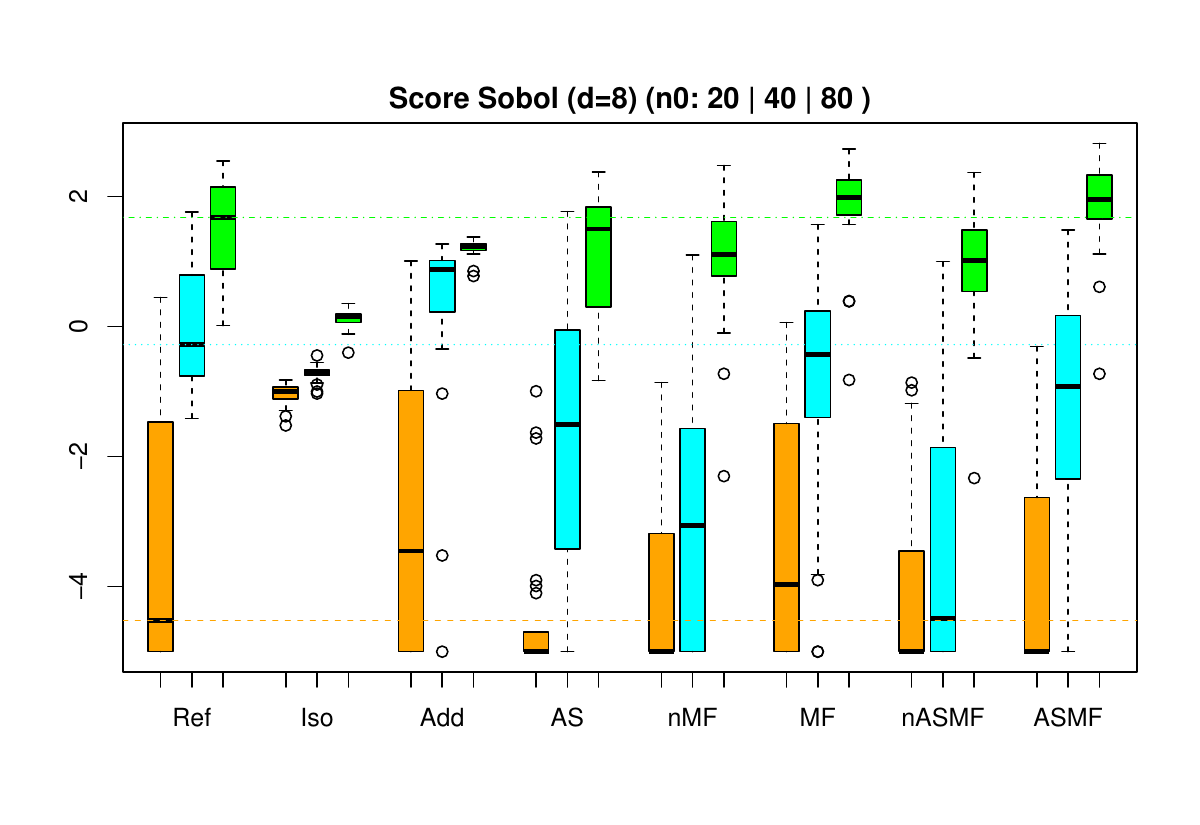}%
\includegraphics[width=0.3333\textwidth, trim= 20 45 30 40, clip]{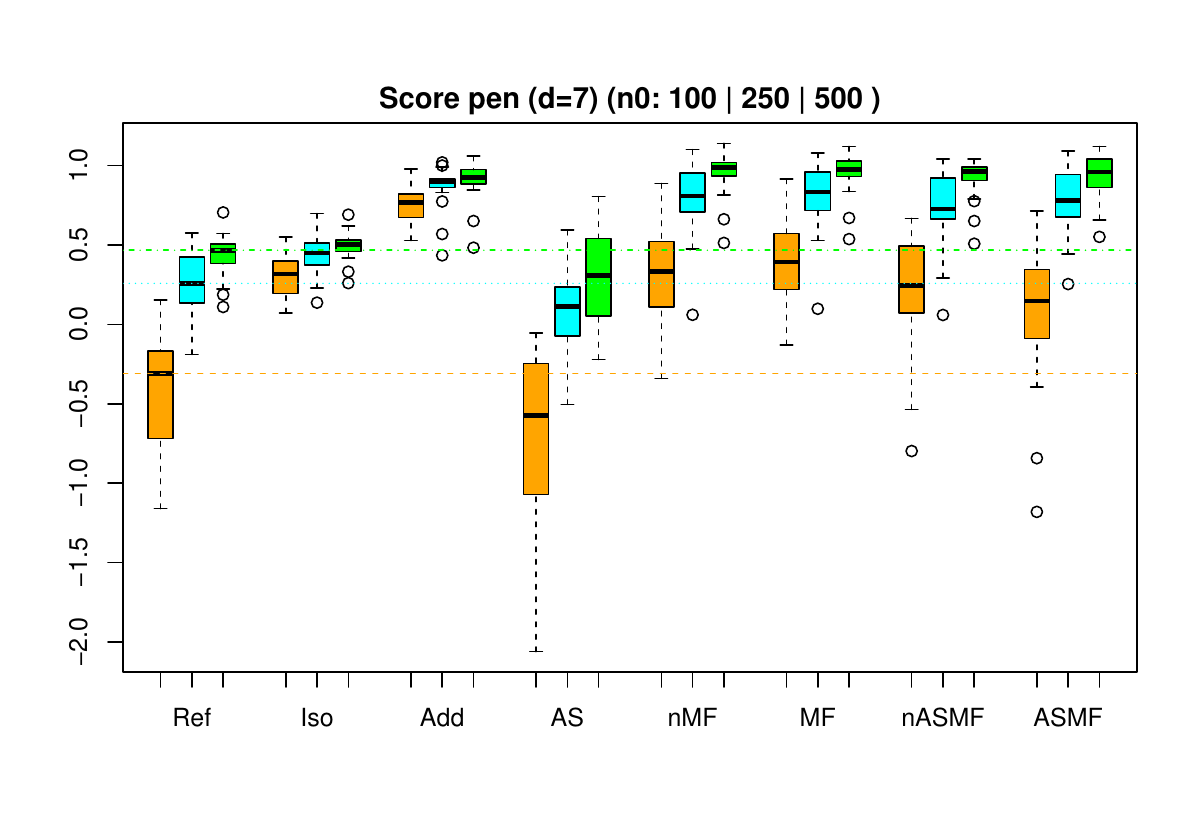}%
\includegraphics[width=0.3333\textwidth, trim= 20 45 30 40, clip]{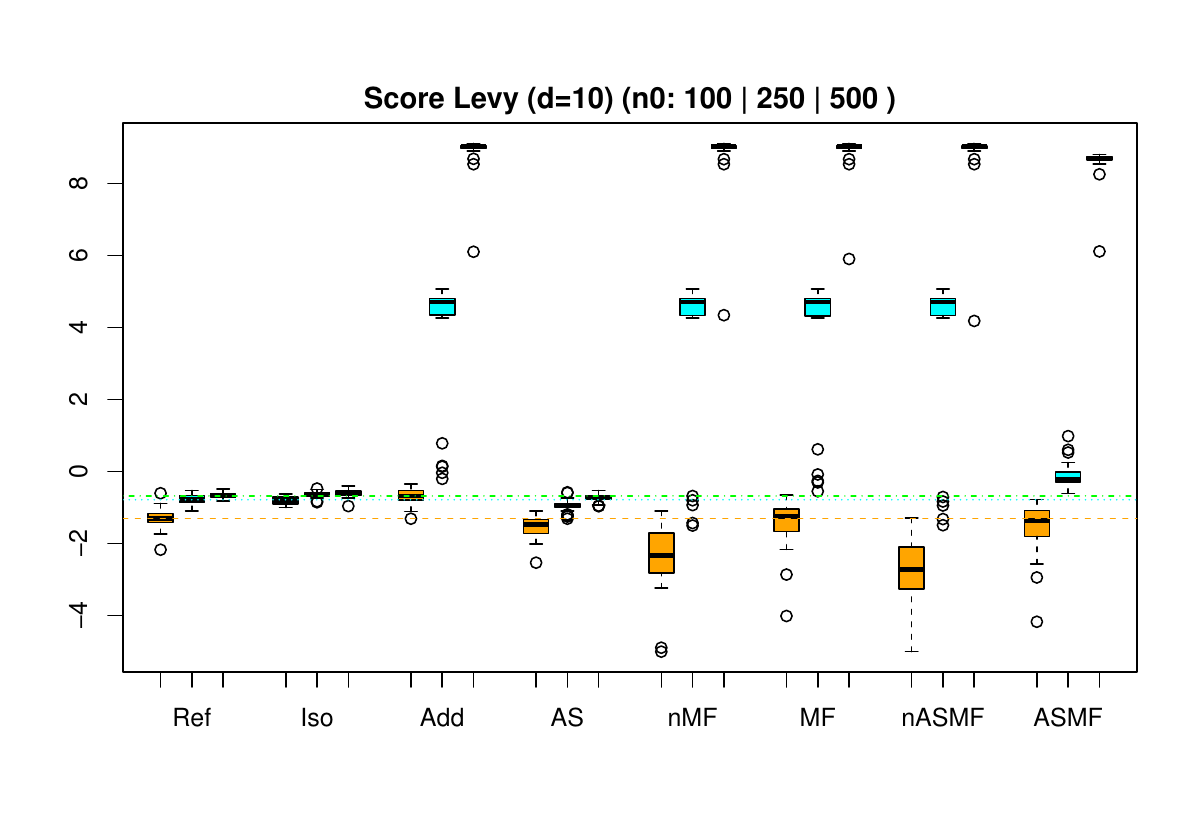}\\
\includegraphics[width=0.3333\textwidth, trim= 30 45 30 40, clip]{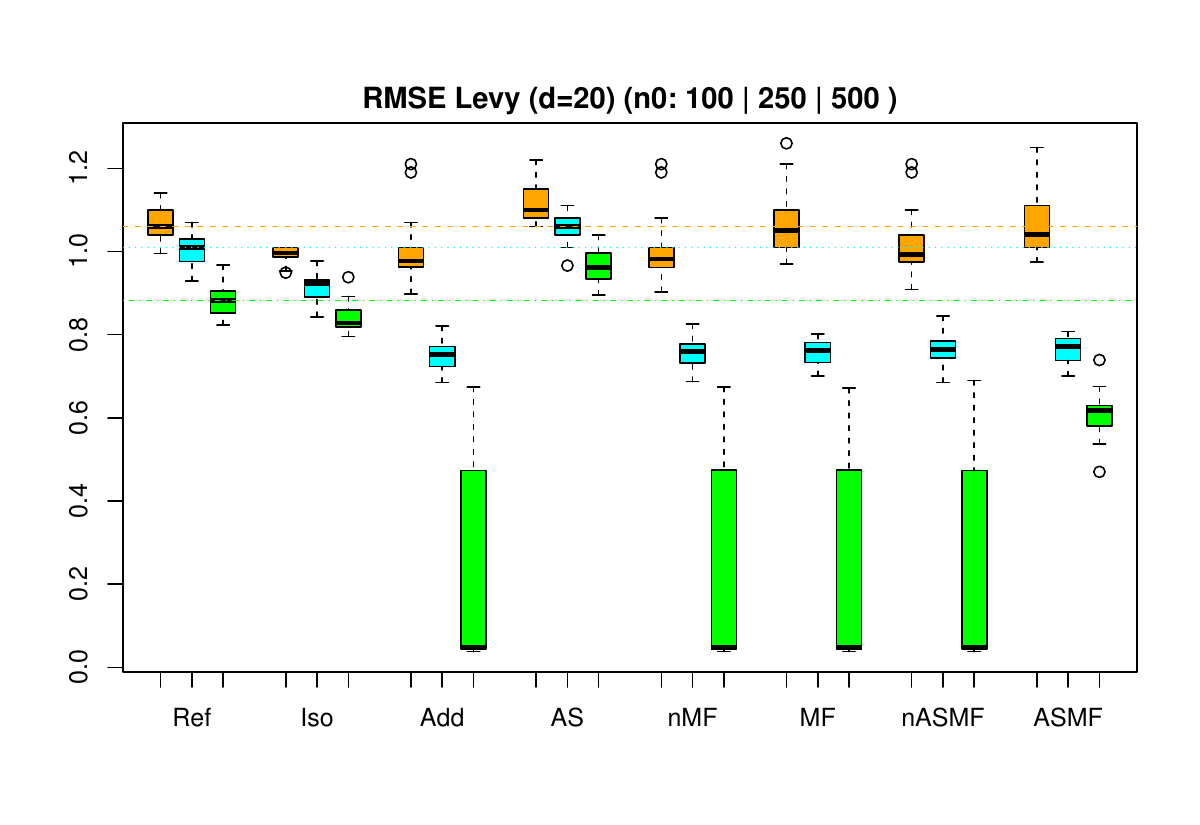}%
\includegraphics[width=0.3333\textwidth, trim= 20 45 30 40, clip]{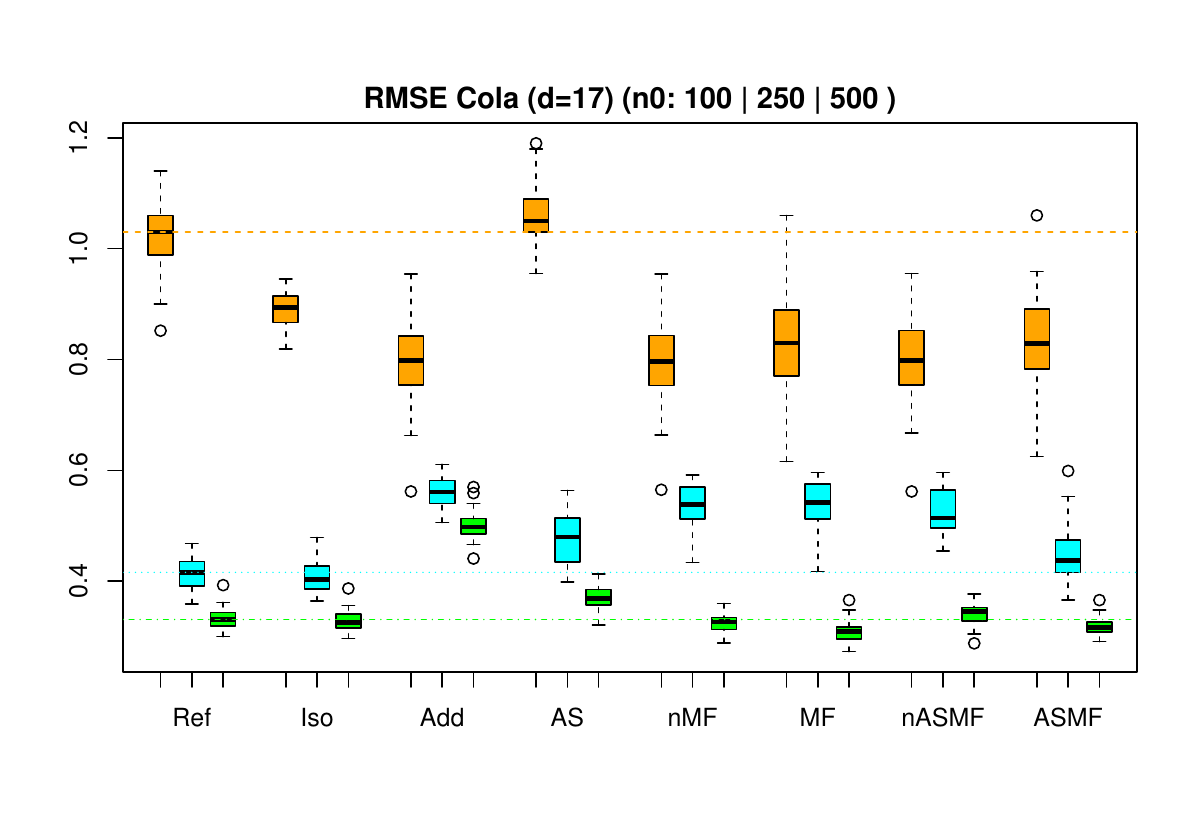}%
\includegraphics[width=0.3333\textwidth, trim= 20 45 30 40, clip]{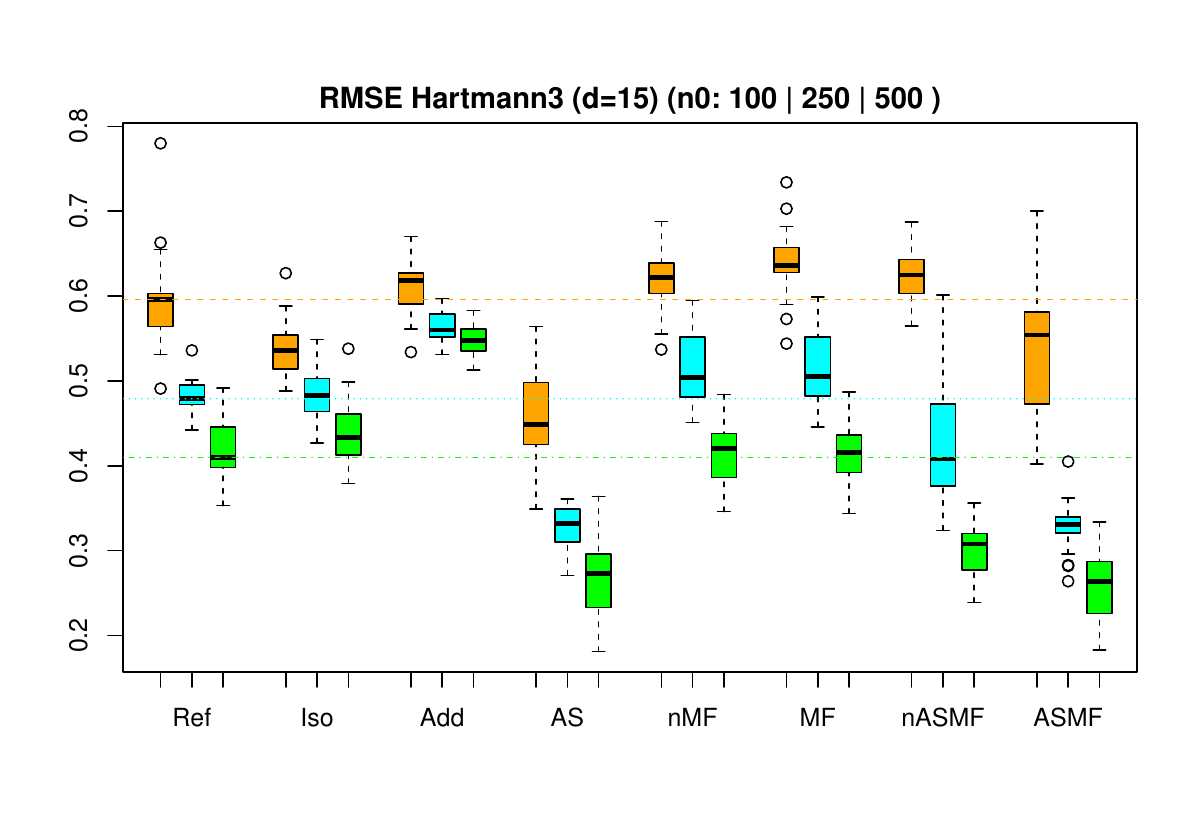}\\
\includegraphics[width=0.3333\textwidth, trim= 30 45 30 40, clip]{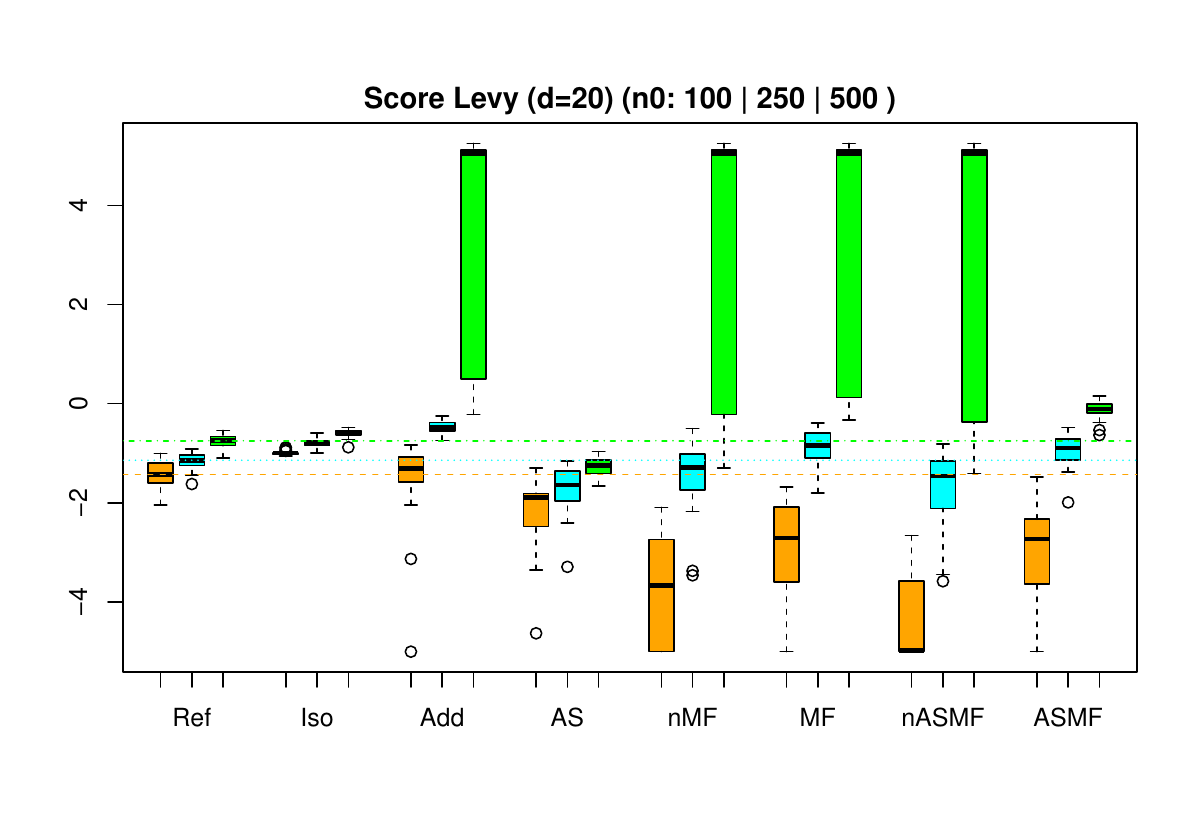}%
\includegraphics[width=0.3333\textwidth, trim= 20 45 30 40, clip]{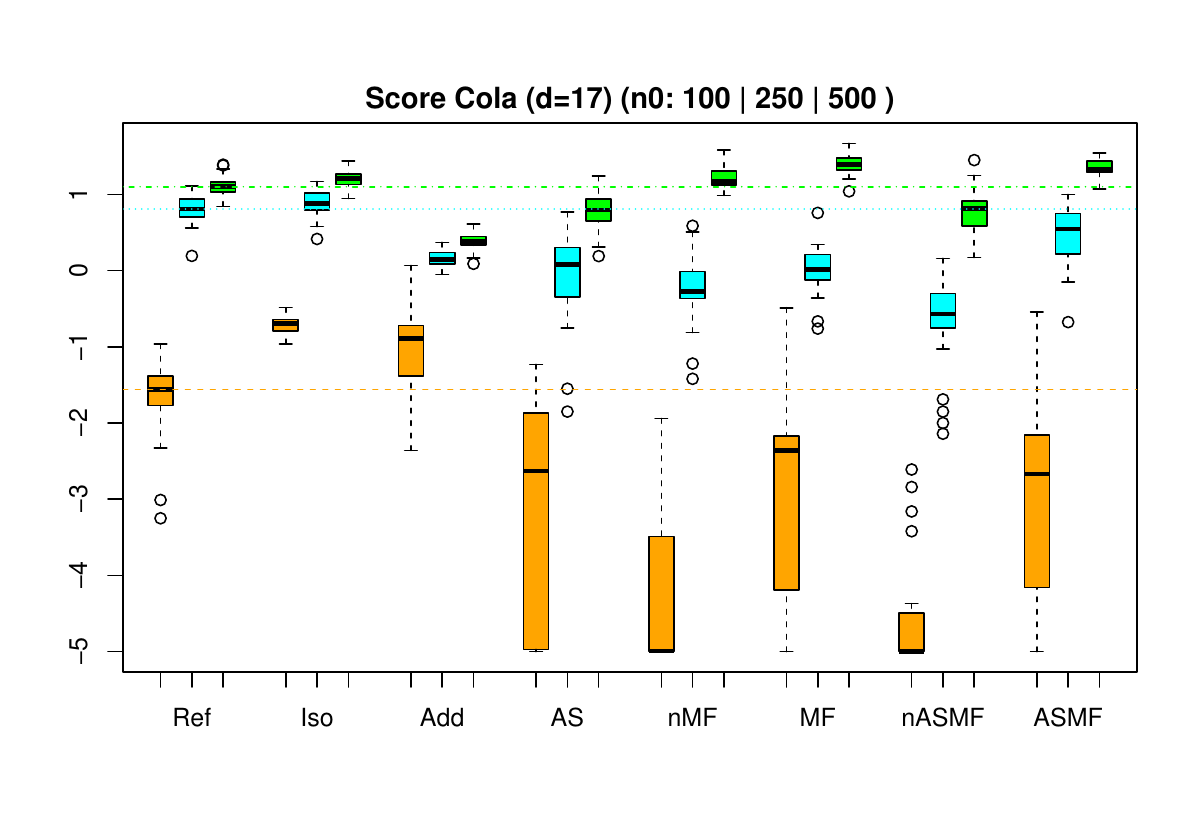}%
\includegraphics[width=0.3333\textwidth, trim= 20 45 30 40, clip]{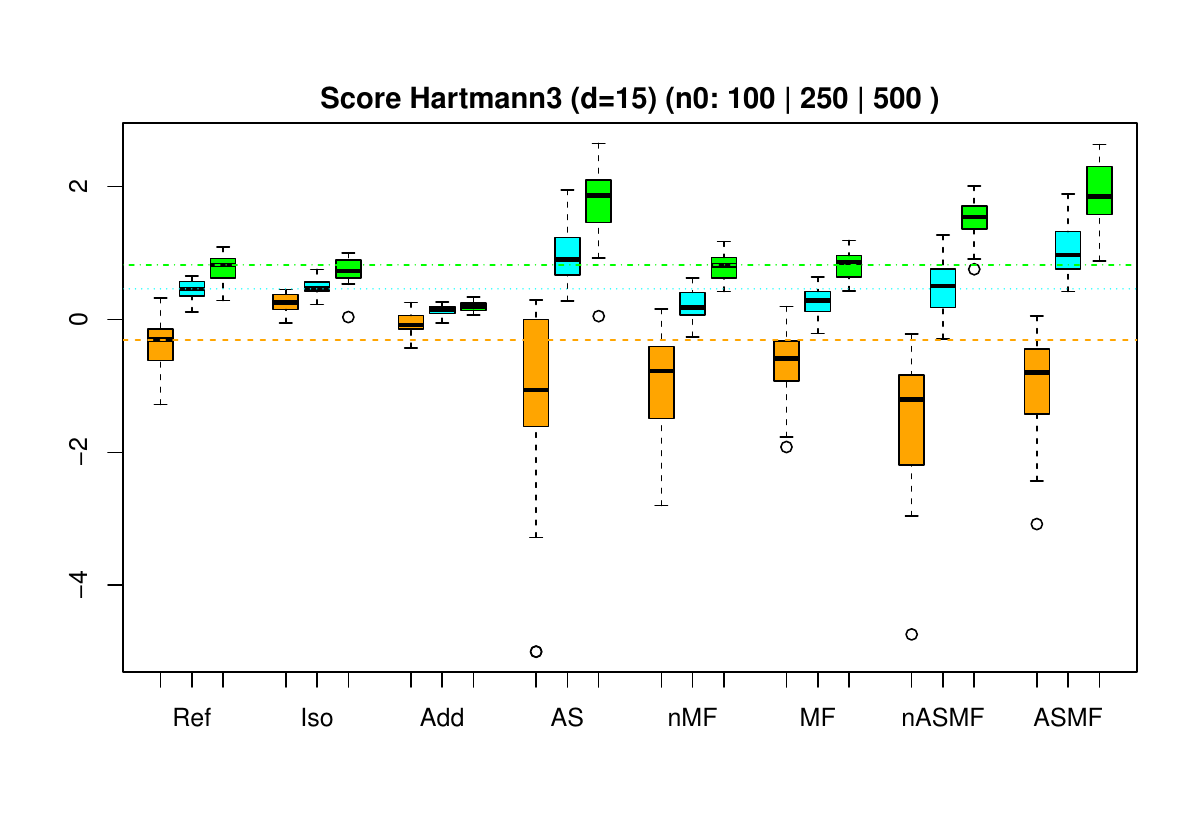}\\
\includegraphics[width=0.3333\textwidth, trim= 30 45 30 40, clip]{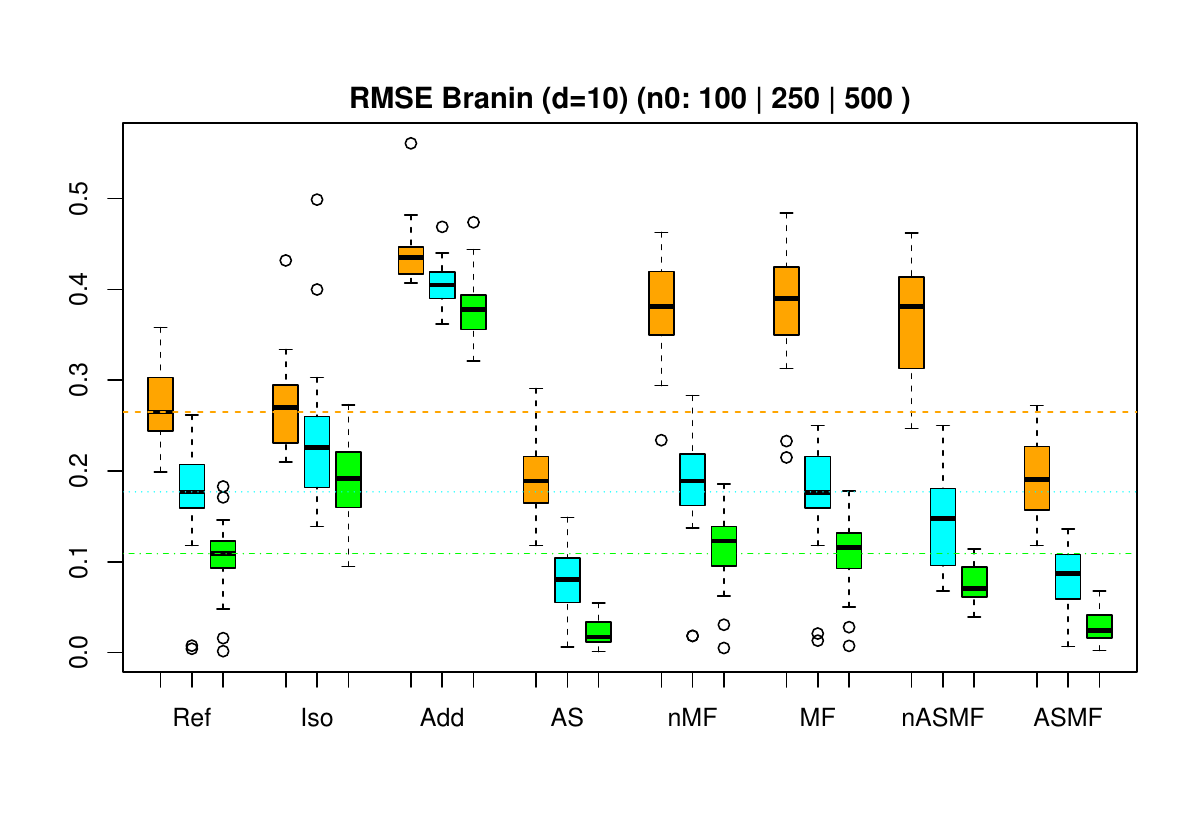}%
\includegraphics[width=0.3333\textwidth, trim= 20 45 30 40, clip]{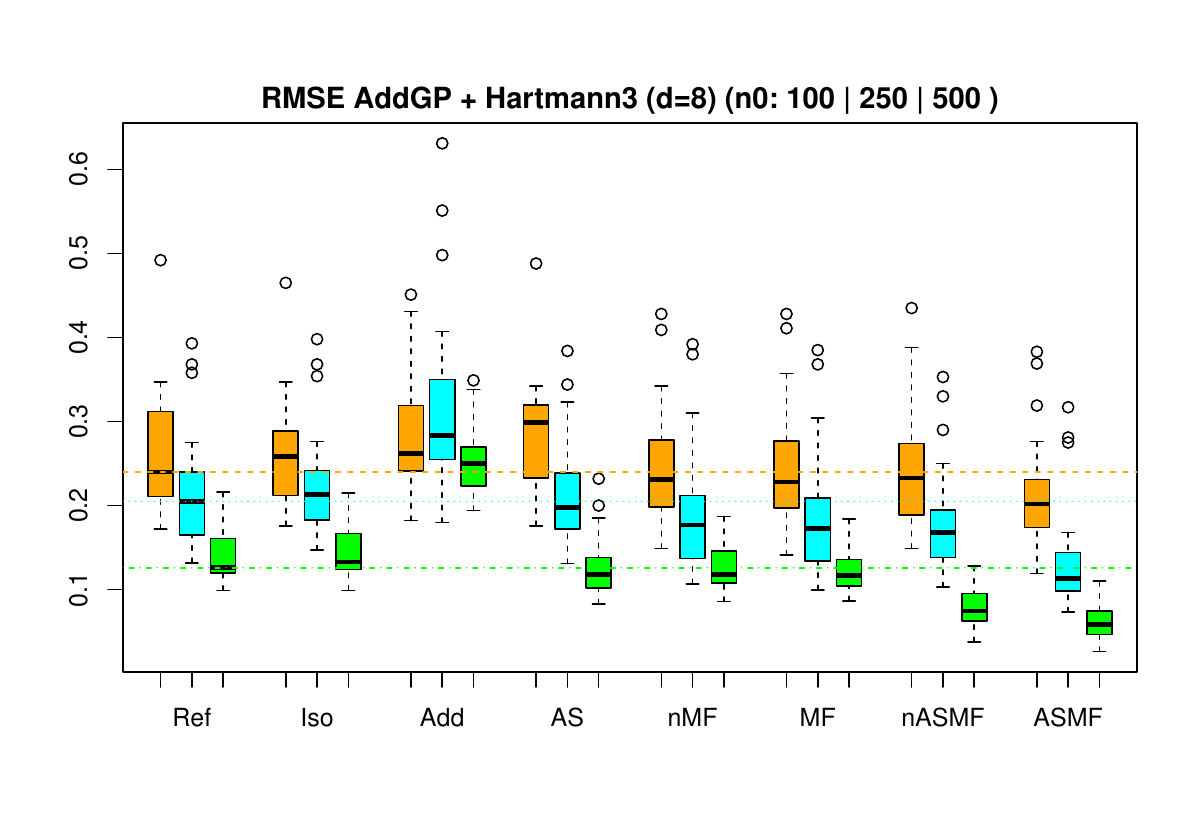}%
\includegraphics[width=0.3333\textwidth, trim= 20 45 30 40, clip]{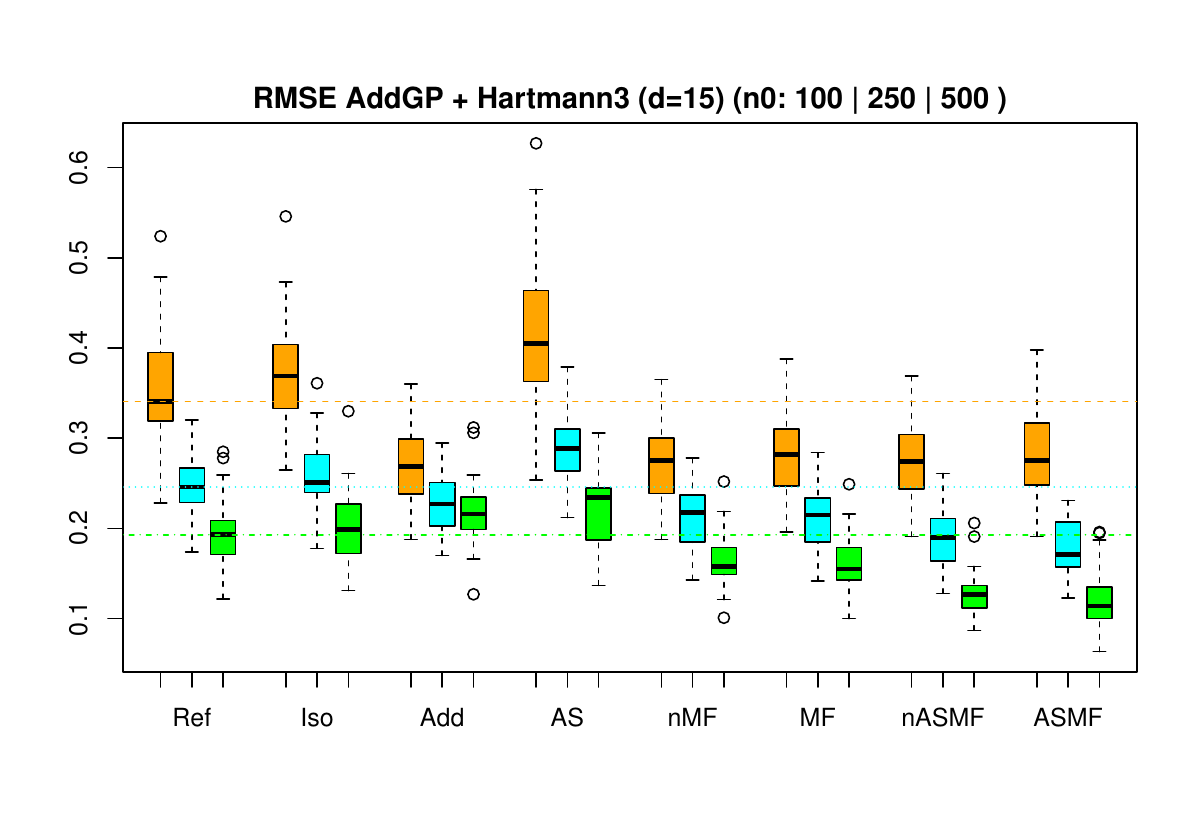}\\
\includegraphics[width=0.3333\textwidth, trim= 30 45 30 40, clip]{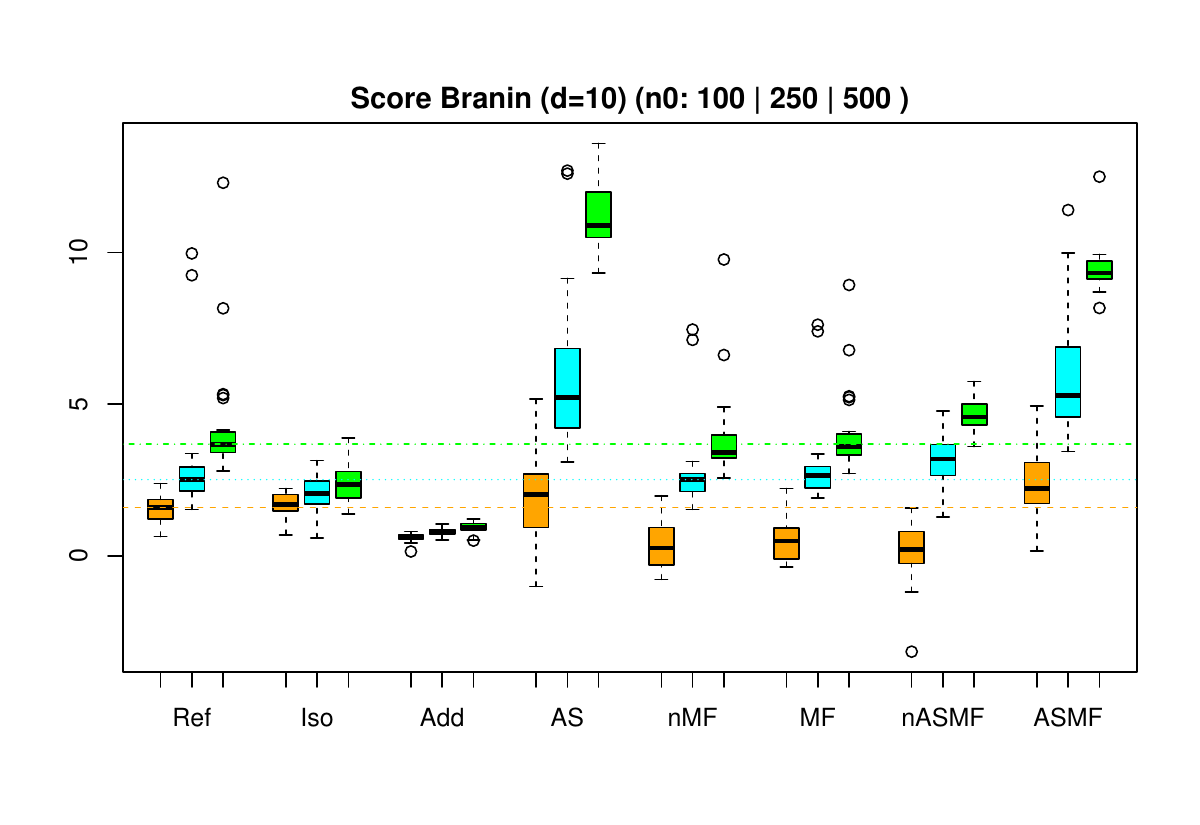}%
\includegraphics[width=0.3333\textwidth, trim= 20 45 30 40, clip]{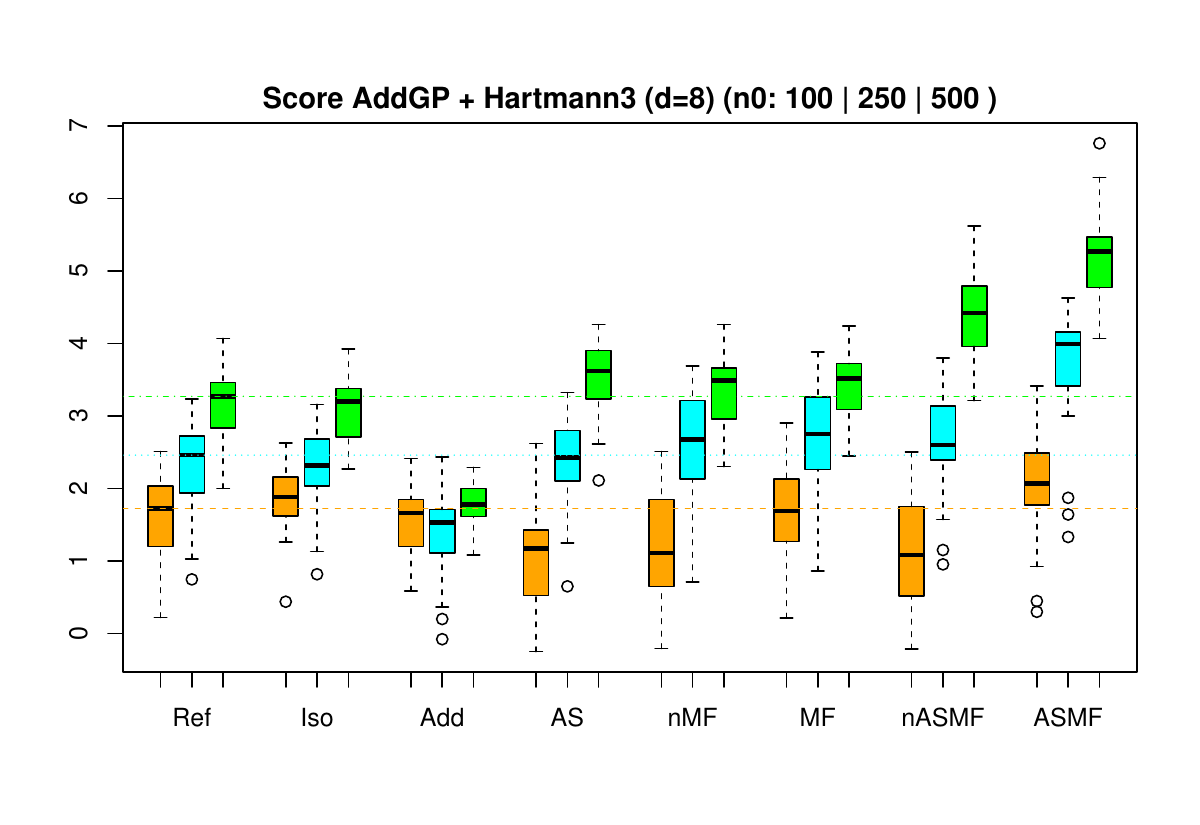}%
\includegraphics[width=0.3333\textwidth, trim= 20 45 30 40, clip]{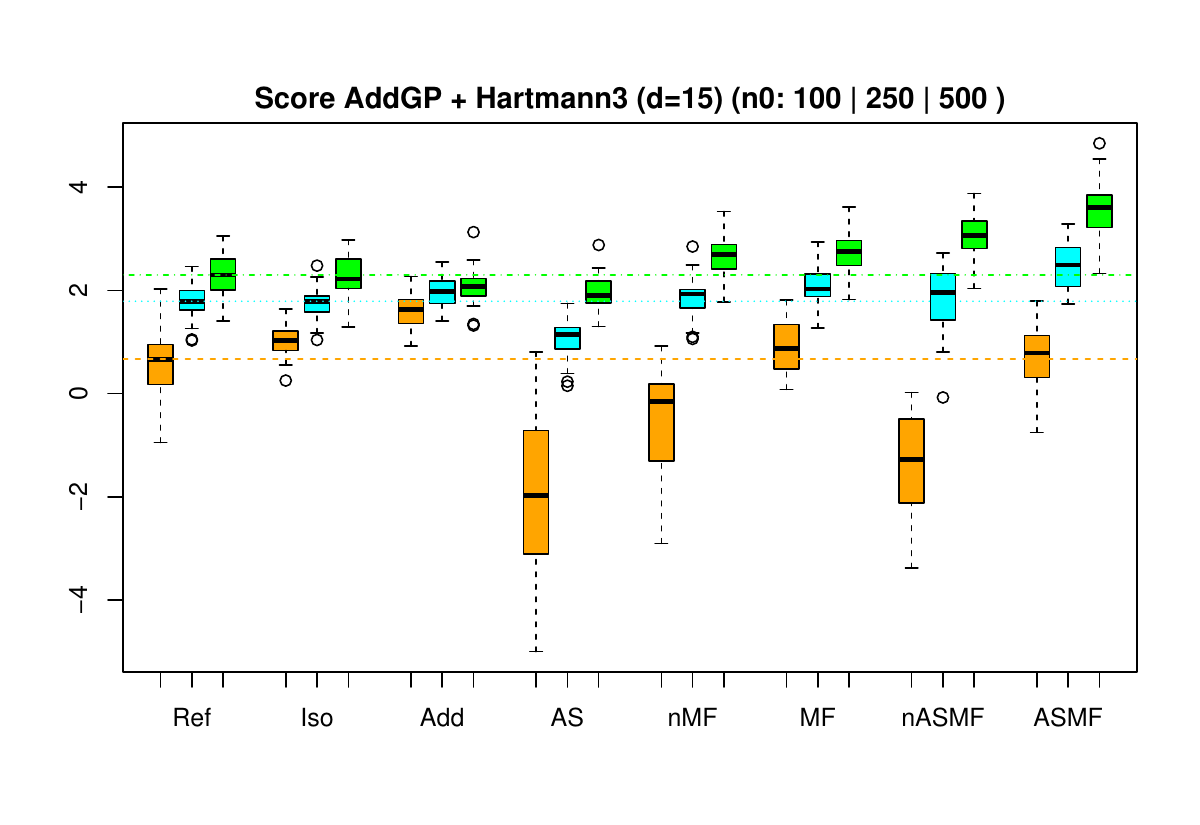}
\caption{Second part RMSE and score results. The different $n_0$ values appear
 with different colors, the corresponding lines indicate the baseline result
 from standard GP models.}
\label{fig:rmseres2}
\end{figure*}

The results over 25 replicates are presented for the RMSE (lower is better)
and score (higher is better) in Figures \ref{fig:rmseres}
and \ref{fig:rmseres2}. We threshold scores at -5 for better visualization.
Before delving into specific details, the multi-fidelity plus AS is in
general at least as good as regular GPs, especially on the largest training
data size. The inference is more difficult on the smallest budget, with a
large variance showing for puma ($d=32$) or CASP ($d=9$). Notably, it can
improve significantly the results when additive or low intrinsic
dimensionality is present. The few exceptions, e.g., on standard GP samples,
occur mostly for the lowest budgets and with a small difference. It is
noteworthy that regular GPs are in general not worse than most alternatives,
and in particular when first-order additive or AS-alone structures are not
present. This underscores the importance of meticulous hyperparameter tuning
for high-dimensional GPs, where larger values can be taken to offset greater
distances, without resulting in conditioning issues of the covariance
matrix.

The naive and AR multi-fidelity models seem to perform similarly on the RMSE,
but the AR model generally yields better scores, both with full and AS
kernels. Perhaps unsurprisingly, the best results of the MF models are
obtained when the additive model performs well: e.g., for Sobol, penicillin,
addGP. Then the AS version of the MF GPs tends to outperform the full GP one,
where the best results are obtained when AS structure is present as well:
penicillin, puma, addGP + Hartmann, or Levy. This improvement can be
attributed to either truncating the number of active dimensions, or simply to
the rotation applied on the input space when keeping all variables. This can
be seen in Appendix in Figure \ref{fig:nbactive}, representing the number of
variables kept by AS and ASMF GP models (also note that the selected
dimension is usually larger than the actual embedding one, which is
recommended in practice \cite{letham2020re} and theoretically
\cite{cartis2022dimensionality}). While additive and AS models excel when such
structure is present, they tend to perform poorly when this assumption does
not hold, especially in terms of score. This effect is more pronounced in
additive models, as AS GP ones can compensate using $r=d$ dimensions,
ultimately not reducing the dimension.

Then, for problems with simultaneously additive and AS structures, the
dedicated GP models perform best. Isotropic GPs may make a reasonable initial
choice with low data, as they are much easier to infer, but they quickly
become less effective than anisotropic ones. Timings are provided in Appendix
\ref{ap:adres}, Figure \ref{fig:timings}, where the repeated optimization of
 the likelihood to estimate the best low dimension $r$ in the AS models
 shows. This could be alleviated by not considering all possible values or by
 warm starting the search. Considering multi-fidelity does not add much
 computational effort. In terms of budget, more data is beneficial to all
 models, as reflected by the larger boxplots for lower budgets. The effect of
 the budget is the most striking on the AS models, suggesting that a minimal
 amount of data is essential for robust inference. Conversely, the additive
 model exhibits the least change with increasing budget and can perform well
 even with the lowest budgets.

Finally, we include an indicative comparison with the higher-order additive
model (OAK) from \cite{Lu2022a}, as depicted in Figure \ref {fig:supprmseres}
in Appendix \ref{ap:adres}. However, the results are harder to interpret: the
predictive variance is not returned to compute scores, additional scalings
are performed, and a measure on the input space is needed. OAK can perform
better than the alternatives on some test cases, but worse on the cases with
active subspaces. This highlights the difficulty of learning high-order
interactions in additive models.

\section{Conclusion and Perspectives}

We propose a simple solution to properly combine the predominant structural
assumptions for high-dimensional modeling: additivity and low intrinsic
dimensionality. The resulting multi-fidelity model is simple to construct and
robust to incorrect assumptions. The promising results obtained open
perspectives in several main directions. First, the inference of GP
hyperparameters for high-dimensional problems may be improved, potentially
starting with GPs that have pre-selected lengthscales, as suggested
in \cite {appriou2023combination}. Then, a comparison could be conducted with
direct inference of the active subspaces matrix within GPs, see,
e.g., \cite {letham2020re} or \cite{Garnett2013}. Given that larger datasets
may be necessary for more precise inference of such features, a combination
with sparse GP models, e.g., as in \cite{moss2023inducing}, would be
considered. This could further include input and output warpings, as already
advocated by \cite{lin2020transformation,Lu2022a}. A second direction to
explore is the use of less linear multi-fidelity models as summarized
by \cite {brevault2020overview}, or even multi-source models, see, e.g.,
\cite{poloczek2017multi}. This would be beneficial when combining models whose
 structural assumptions have no natural ordering, like BOCK \cite
 {oh2018bock}, and higher order additive models \cite{Lu2022a}. Non-linear
 dimension reduction is another appealing candidate, see, e.g.,
\cite{guhaniyogi2016compressed}. Lastly, GP modeling shines in sequential
procedures, where existing works only focus on individual aspects, say
additivity \cite{schwabe1995designing}, multi-fidelity \cite{le2015cokriging}
or active subspace estimation \cite{Wycoff2021}. Future research could delve
into the alignment of these goals compared to Bayesian optimization, exploring
how these aspects synergize in sequential decision-making processes.

\appendix

\section{One Shot Active Subspace Gaussian Process Learning}
\label{ap:asgp}

In practice, $\C = \int_\Xset
\nabla(f(\x))^\top
\nabla(f(\x)) \lambda(d\x)$ is often estimated by Monte Carlo, either directly
on $f$ when the gradient is available, or on its surrogate. In
\cite{Wycoff2021}, this AS matrix is expressed in closed form for a GP, which
can be further used to reduce the dimension. Let $k$ be a twice-differentiable
kernel, with derivatives $\kappa$, $\Wij :=
    \int_\X \boldsymbol{\kappa}_i(X)^\top \boldsymbol{\kappa}_j(X) d\lambda$,
    and $E_{i,j} :=
    \int_\X \frac{\partial^2 k (X, X)}{\partial x_i \partial x_j}
    d\lambda$. Then, $C^{(n)}_{i,j} = E_{i,j} - tr\left(\K_n^{-1} \Wij \right) +
    \y^\top \K_n^{-1} \Wij \K_n^{-1} \y$.
As a result, estimating a GP with reduced dimension is possible via a two-step
process: first, fit a high-dimensional GP, then use the corresponding AS
matrix to learn a low dimensional GP on the projected data, with kernel $k(\x,
\x') = \tilde{k}(\A^\top \x, \A^\top \x')$, assuming a centered $\Xset$, where
$\A = \U_r$, the first $r$ eigen-vectors of $\C^{(n)}$.

Nevertheless, the AS matrix $\C^{(n)}$ in fact only depends on the $d$
lengthscale hyperparameters, and the data. What we propose here is to use the
same parameterization of the AS matrix of the GP, but learn the parameters via
the likelihood of a low dimensional GP, i.e., learn all
hyperparameters: $l_1, \dots, l_r, \theta_1, \dots
\theta_d$  at once, where the $l_i$ (resp.\ $\theta_i$) are the low (high)
dimensional GP lengthscales.

We rely on the work made previously with the derivative of an AS kernel in \cite{Wycoff2021}, and give the additional required expressions, that is $\frac{\partial \C}{\partial \theta_i}$:
$$\frac{\partial C^{(n)}_{i,j}}{\partial \theta_i} = \frac{ \partial E_{i,j} - tr\left( \K_n^{-1} \Wij \right) +
    \y^\top \K_n^{-1} \Wij \K_n^{-1} \y}{\partial \theta_i}.$$

Hence we need $\frac{\partial E_{i,j}}{\partial \theta_i}$, $\frac{\partial \K_n^{-1}}{\partial \theta_i}$ and $\frac{\partial \Wij}{\partial \theta_i}$. These are combined to get:

\small
\begin{align*}
  \begin{autobreak}
\frac{\partial C^{(n)}_{i,j}}{\partial \theta_i} = \frac{\partial E_{i,j}}{\partial \theta_i} - 
tr\left( \frac{\partial \K_n^{-1}}{\partial \theta_i} \Wij + \K_n^{-1} \frac{\partial \Wij}{\partial \theta_i} \right) +
\y^\top \frac{\partial \K_n^{-1}}{\partial \theta_i} \Wij \K_n^{-1} \y + \y^\top \K_n^{-1} \frac{\partial \Wij}{\partial \theta_i} \K_n^{-1} \y + \y^\top \K_n^{-1} \Wij \frac{\partial \K_n^{-1}}{\partial \theta_i} \y.
  \end{autobreak}
\end{align*}

\normalsize

and up to the likelihood level:

$$\frac{\partial \log L}{\partial \theta_i} = \frac{\partial \left(const. - \frac{n}{2} \log \y^\top \K^{-1} \y - \frac{1}{2} \log |\K| \right)}{\partial \theta_i} = \frac{n}{2 \hat{\sigma}^2} \y^\top \K^{-1} \frac{\partial \K}{\partial \theta_i} \K^{-1} \y -\frac{1}{2} Tr \left( \K^{-1} \frac{\partial \K}{\partial \theta_i} \right)$$

where, using the chain rule $\frac{\partial K_{i,j}}{\partial \theta_i} = \frac{\partial K_{i,j}}{\partial \U} \frac{\partial \U}{\partial \theta_i}$. More precisely, using \cite{Petersen2008}, involving the eigen-vectors $\U_l$ and corresponding eigen-value $\lambda_l$ of $\C$, and pseudo-inverses $\dagger$: 
$$\frac{\partial \U_l}{\partial \theta_i} = (\lambda_l \Id - \C)^{\dagger} \frac{\partial \C}{\partial \theta_i} \U_l.$$

As an example, for a Gaussian kernel in the reduced dimension too, $\mathbf{h} = (\x_i - \x_j)$, such that:
$\frac{\partial K_{i,j}}{\partial \W} = 2 Diag(\boldsymbol{l}) \W \mathbf{h} \mathbf{h}^\top K_{i,j}.$

For this Gaussian kernel case, parameterized by $k(\x, \x') = \prod \limits_{i = 1}^d \exp \left(- \left(\frac{x_i - x'_i}{\theta_i} \right)^2 \right)$ (denoting $a = x_i$, $b = x'_i$ and $t=\theta_i$):

$\frac{\partial E_{i,j}}{\partial \theta} =  \frac{\partial}{\partial \theta} \frac{\delta_{i=j}}{\theta^{2}} = -2 \delta_{i=j} \theta^{-3}$
\scriptsize
\begin{align*}
  \begin{autobreak}
  \MoveEqLeft
  \frac{\partial w_{i,i}(a, b, t)}{\partial t} = 
  -\dfrac{\left(2t\left(\left(b+a-2\right)\left(2t^2-\left(b+a-2\right)^2\right)\mathrm{e}^\frac{b+a}{t^2}-\left(b+a\right)\left(2t^2-\left(b+a\right)^2\right)\mathrm{e}^\frac{1}{t^2}\right)\mathrm{e}^\frac{-b^2-a^2 - 2}{2t^2}\right)}{16t^6}
    -\dfrac{\left(\sqrt{{\pi}}\left(\operatorname{erf}\left(\frac{b+a}{2t}\right)
  -\operatorname{erf}\left(\frac{b+a-2}{2t}\right)\right)\left(2t^2-2bt+2at-b^2+2ab-a^2\right)\left(2t^2+2bt-2at-b^2+2ab-a^2\right)\mathrm{e}^\frac{-b^2+2ab-a^2-1}{4t^2} \right)}{16t^6}
  \end{autobreak}
\end{align*}

\begin{align*}
  \begin{autobreak}
  \MoveEqLeft
  \frac{\partial w_{i,j}(a, b, t)}{\partial t} = 
\dfrac{-2\left(a-b\right)t\mathrm{e}^{-\frac{\left(a-b\right)^2}{4t^2}}\left(\left(b+a\right)\mathrm{e}^{-\frac{\left(b+a\right)^2}{4t^2}}-\left(b+a-2\right)\mathrm{e}^{-\frac{\left(b+a-2\right)^2}{4t^2}}\right) }{8t^4}
+\dfrac{4t\left(-\left(b^2+a^2\right)\mathrm{e}^{-\frac{b^2+a^2}{2t^2}}-\left(-b^2+2\left(b+a-1\right)-a^2\right)\mathrm{e}^\frac{-b^2+2\left(b+a-1\right)-a^2}{2t^2}\right)}{8t^4}
+ \dfrac{
2\sqrt{{\pi}}\left(a-b\right)\left(\operatorname{erf}\left(\frac{b+a-2}{2t}\right)-\operatorname{erf}\left(\frac{b+a}{2t}\right)\right)t^2\mathrm{e}^{-\frac{\left(a-b\right)^2}{4t^2}}-\sqrt{{\pi}}\left(a-b\right)^3\left(\operatorname{erf}\left(\frac{b+a-2}{2t}\right)-\operatorname{erf}\left(\frac{b+a}{2t}\right)\right)\mathrm{e}^{-\frac{\left(a-b\right)^2}{4t^2}}}{8t^4}
  \end{autobreak}
\end{align*}

\begin{align*}
  \begin{autobreak}
  \MoveEqLeft
  \frac{\partial I_{l,l}(a, b, t)}{\partial t} = 
\dfrac{2\sqrt{{\pi}}t^2\mathrm{e}^{-\frac{\left(b-a\right)^2}{4t^2}}\left(\operatorname{erf}\left(\frac{b+a}{2\left|t\right|}\right)-\operatorname{erf}\left(\frac{b+a-2}{2\left|t\right|}\right)\right)+\sqrt{{\pi}}\left(b-a\right)^2\mathrm{e}^{-\frac{\left(b-a\right)^2}{4t^2}}\left(\operatorname{erf}\left(\frac{b+a}{2\left|t\right|}\right)-\operatorname{erf}\left(\frac{b+a-2}{2\left|t\right|}\right)\right)}{4t^2}
+\dfrac{2\mathrm{e}^{-\frac{\left(b-a\right)^2}{4t^2}}\left(\left(b+a-2\right)\mathrm{e}^{-\frac{\left(b+a-2\right)^2}{4t^2}}-\left(b+a\right)\mathrm{e}^{-\frac{\left(b+a\right)^2}{4t^2}}\right)}{4t}
  \end{autobreak}
\end{align*}

\normalsize

\section{Complements on the Auto-regressive Multi-fidelity Model}
\label{ap:mufi}

For the implementation, we start by giving some log-likelihood derivatives.
Then we discuss the link to the recursive formulation of the AR multi-fidelity
model.

\subsection{Log-likelihood Derivatives}

Derivatives of the log-likelihood are given, e.g., in \cite{Forrester2008}.
The only change is the log-likelihood derivative that requires adaptation in
the additive case.

For the coarse level model with an additive kernel where $\K = \sum \limits_{i = 1}^d \alpha_i \K^{(i)}(\theta_i) + g \Id$:
$$L = -n/2 \log(2 \pi) -1/2 \y_C^\top \K^{-1} \y_C - 1/2 \log |\K|$$ 
$$\frac{\partial L}{\partial \theta_i} = 1/2 \y_C^\top \K^{-1} \frac{\partial \K}{\partial \theta_i} \K^{-1} \y_C -1/2 Tr \left(\K^{-1} \frac{\partial \K}{\partial \theta_i} \right)$$
$$\frac{\partial L}{\partial \alpha_i} = 1/2 \y_C^\top \K^{-1} \frac{\partial \K}{\partial \alpha_i} \K^{-1} \y_C -1/2 Tr \left(\K^{-1} \frac{\partial \K}{\partial \alpha_i} \right) = 1/2 \y_C^\top \K^{-1} \K^{(i)} \K^{-1} \y_C -1/2 Tr \left(\K^{-1} \K^{(i)} \right)$$

Then for the AR multi-fidelity kernel:

$$\frac{\partial \tilde{L}}{\partial \rho} = - n/2 \frac{\partial \log(\hat{\sigma}_d^2)}{\partial \rho} = %-1/2 \frac{\partial (\y_E - \rho \y_C)^\top \K^{-1} (\y_E - \rho \y_C)}{\partial \rho} = 
\frac{\y_C \K^{-1} (\y_E - \rho \y_C)}{\hat{\sigma}_d^2}$$

\subsection{Recursive Formulation}

A recursive formulation of the multi-fidelity AR model was proposed
by \citep{le2014recursive}, which is equivalent to the one
by \cite{kennedy2000predicting} but only in the deterministic case. This
formulation writes:
\begin{align}
\begin{split}
m_{n,E}(\x) &= m_{n, C}(\x) + k_E(\x, \XE) \KE^{-1} (\y^{(E)} - \rho \y^{(C)}),\\
s_{n,E}^{2}(\x) &= \rho^2 s_{n,C}(\x) + k_E(\x, \x) - k_E(\x, \XE) \KE^{-1} k_E(\XE, \x)
\label{eq:recarmufi}
\end{split}
\end{align}
which reduces the computational complexity of fitting the finer level and
gives the predictive quantities at all fidelity levels.
We give here a simple proof of the equivalence
in this case, then show that it does not apply in the noisy one.

\subsubsection{Deterministic Case}

In the deterministic case, when both designs are equal, $\XC = \XE$ (of size $n$):

$\tilde{\K} = \begin{bmatrix}
k_C(\XC, \XC) & \rho  k_C(\XC, \XE) \\
 \rho k_C(\XE, \XC) & \rho^2  k_C(\XE, \XE) + k_E(\XE, \XE)
\end{bmatrix} 
:= \begin{bmatrix} \KC & \rho \KC\\
\rho \KC & \rho^2 \KC + \KE 
\end{bmatrix}$ 

Similarly, $\tilde{\veck}(\x) := [\rho \kC, \rho^2 \kC + \kE]^\top$ for shorter notation (dropping the dependence on $\x$).

Then the block-matrix inverse formula \cite{Petersen2008} gives, following the notations there:
$\C_1 = \KC - \rho^2 \KC (\KE + \rho^2 \KC )^{-1} \KC$ and $\C_2 = \rho^2 \KC + \KE - \rho^2 \KC = \KE$ where this second equality is used for expressing $\tilde{\K}^{-1}$:
$\tilde{\K}^{-1} = \begin{bmatrix}
\rho^2 \KE^{-1} + \KC^{-1} & -\rho \KE^{-1}\\
-\rho \KE^{-1} & \KE^{-1}  
\end{bmatrix}.$

Consequently, for the predictive equations:

$\tilde{\K}^{-1} \tilde{\veck} = [\rho^3 \KE^{-1} \kC + \rho \KC^{-1} \kC - \rho^3 \KE^{-1} \kC - \rho \KE^{-1} \kE,  -\rho^2 \KE^{-1} \kC + \rho^2 \KE^{-1} \kC + \KE^{-1} \kE] = [\rho \KC^{-1} \kC - \rho \KE^{-1} \kE, \KE^{-1} \kE]$

such that 

$m_{n,E}(\x) = \tilde{\veck}^\top \tilde{\K}^{-1} \tilde{\y} = \rho \kC \KC^{-1} \y_C - \rho \kE \KE^{-1} \y_C + \kE \KE^{-1} \y_E = m_{n,C}(\x) + \kE \KE^{-1} (\y_E - \rho \y_C)$

and 

$s_{n,E}^2(\x) = \rho^2 k_C(\x, \x) + k_E(\x, \x) - \tilde{\veck}^\top \tilde{\K}^{-1} \tilde{\veck} = \rho^2 k_C(\x, \x) + k_E(\x, \x) -$\\  $\left(\rho^2 \kC \KC^{-1} \kC - \rho^2 \kC \KE^{-1} \kE + \rho^2 \kC \KE^{-1} \kE + \kE \KE^{-1} \kE \right) = %\rho^2 \kC \KC^{-1} \kC + \kE \KE^{-1} \kE = 
\rho^2 s^2_{n,C}(\x) + s^2_{n,d}(\x)$

For the non equal DoE, $\XC$ can be split into $[\XE, \XR]$ where $\XR$ are the designs where only the cheap level is evaluated. In this case, we still have that $\C_2 = \rho^2 \KC(\XE, \XE) + \KE - \rho^2 \KC(\XE, \XE) = \KE$. The last term is obtained by realizing that $\A_{12}$ is the first $n_E$ lines of $\KC(\XC, \XC) = \A_{11}$. Then $\A_{12} \A_{11}^{-1} = [\Id, \mathbf{0}]$ and $\A_{11}^{-1} \A_{21} = [\Id, \mathbf{0}]^\top$ such that:

$\tilde{\K}^{-1} = \begin{bmatrix}
\begin{bmatrix} \rho^2 \KE^{-1} & \mathbf{0} \\ \mathbf{0} & \mathbf{0} \end{bmatrix} + \KC(\XC, \XC)^{-1} & \begin{bmatrix} -\rho \KE^{-1}\\ \mathbf{0} \end{bmatrix}\\
\begin{bmatrix} -\rho \KE^{-1} & \mathbf{0} \end{bmatrix} & \KE^{-1}  
\end{bmatrix}$

Hence, with $\tilde{\veck} = [\rho k_C(\XE, \x), \rho k_C(\XR, \x), \rho^2
k_C(\XE, \x) + \kE]$, the same simplifications as above occur, giving the
equivalence between the two formulations. From the expression of
$\tilde{\K}^{-1}$, only $n\times n$ inverses and determinants need to be
computed.

\subsubsection{Noisy Coarse Function}

Now assume that the low fidelity function is noisy, as is usually the case for
an additive model. This time, the application of the block inverse matrix when
$\XE = \XC$ on:

$\check{\K} = \begin{bmatrix}
k_C(\XC, \XC) + g \Id & \rho k_C(\XC, \XE) \\
 \rho  k_C(\XE, \XC) & \rho^2  k_C(\XE, \XE) + k_E(\XE, \XE)
\end{bmatrix} 
:= \begin{bmatrix} \KC + \mathbf{D} & \rho \KC \\
\rho \KC & \rho^2 \KC + \KE 
\end{bmatrix}$

gives, (following again notations from \cite{Petersen2008}): $\C_1 = (\KC +
\mathbf{D}) - \rho^2 \KC (\KE + \rho^2 \KC )^{-1} \KC$ and $\C_2 = \rho^2 \KC
+ \KE - \rho^2 \KC (\KC + \mathbf{D})^{-1} \KC = \KE + \rho^2 (\KC^{-1} +
  g^{-1} \Id)^{-1} = \KE + g \rho^2 \KC (\KC + g \Id)^{-1}$ using the Woodbury
  identity. Even though it allows to reduce the computational complexity of
  the direct multi-fidelity approach, it does not lead to the expressions from
  the recursive formulation. In particular, the recursive variance expression
  does not equal zero at $\XE$ since the low fidelity variance is greater
  than zero.

\section{Additional Results}
\label{ap:adres}

\begin{figure*}[htpb]
\centering
\includegraphics[width=\textwidth, trim= 20 40 20 20, clip]{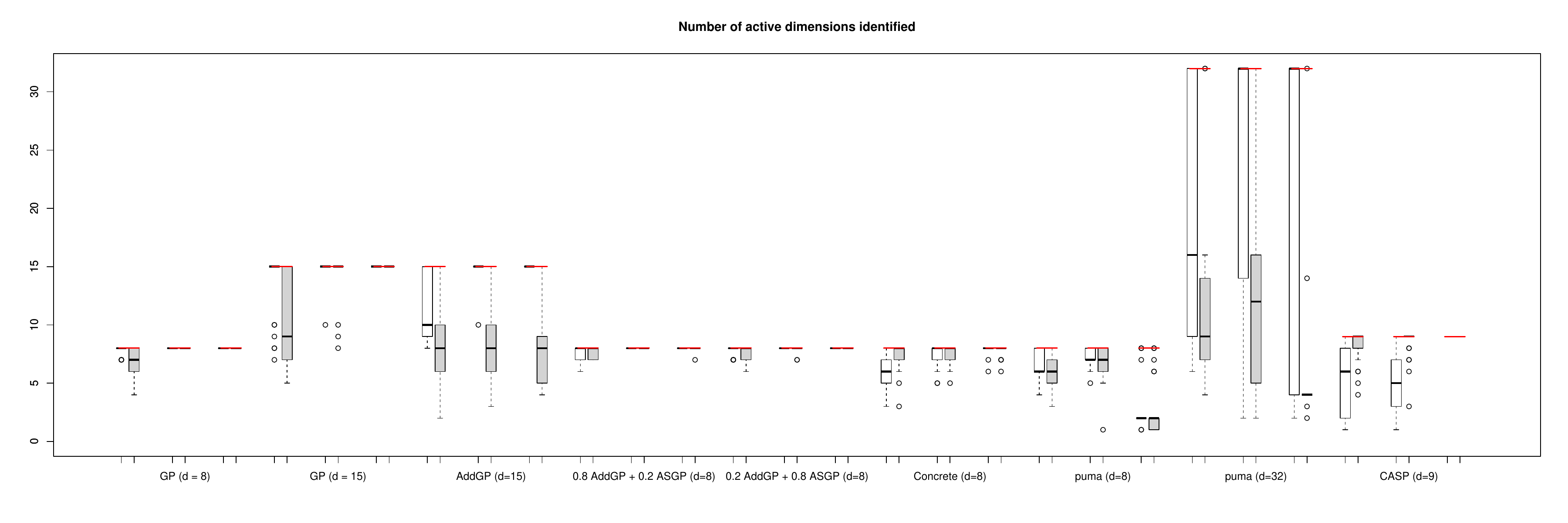}\\
\includegraphics[width=\textwidth, trim= 20 40 20 20, clip]{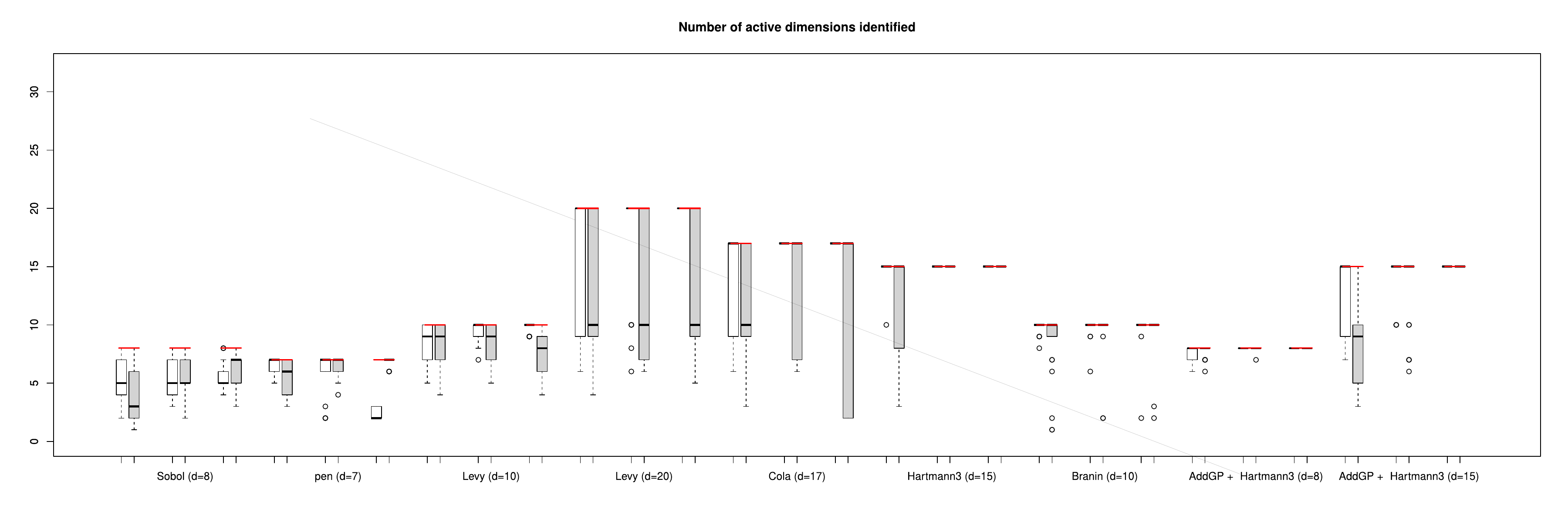}\\
\caption{Number of active dimensions kept for the AS GP (left boxplots) and MF
AS GP (right gray-filled boxplots). The red segments indicate the number of
variables of the problem ($d$).}
\label{fig:nbactive}
\end{figure*}

\begin{figure*}[htpb]
\centering
\includegraphics[width=0.3333\textwidth, trim= 30 45 30 40, clip]{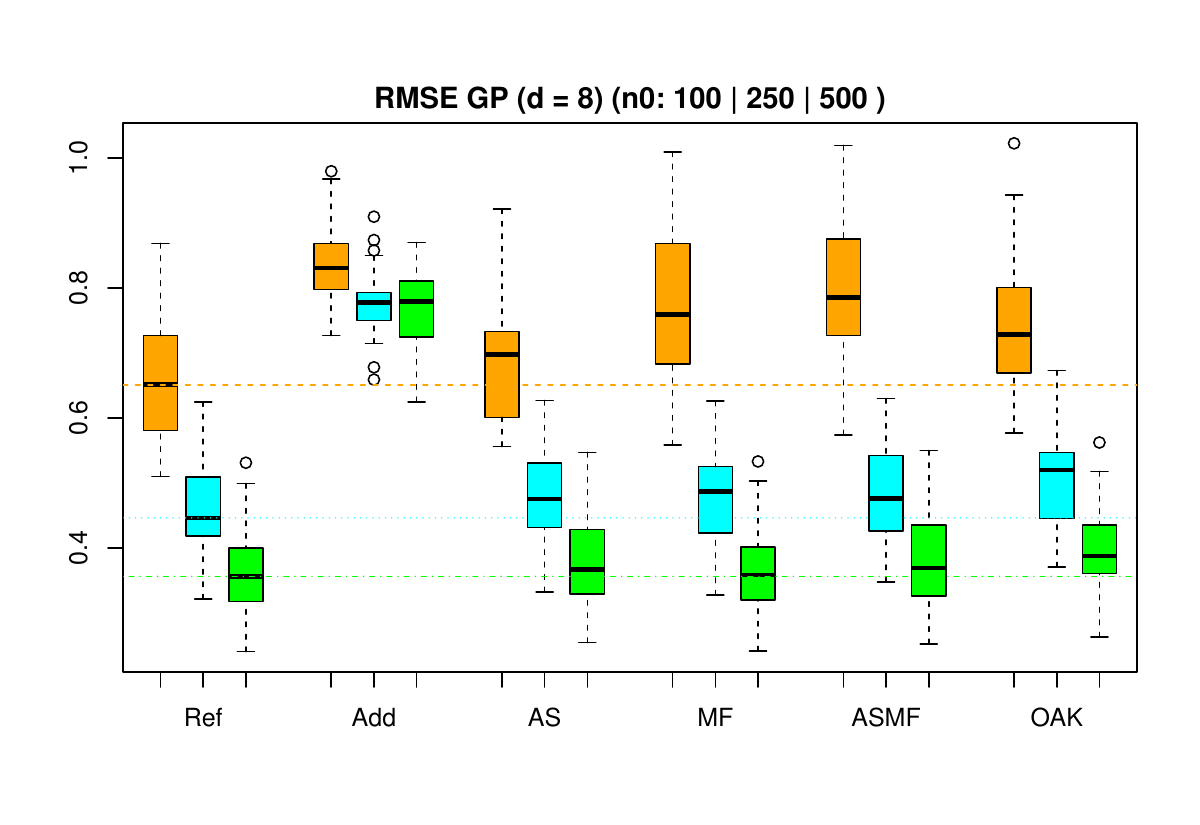}%
\includegraphics[width=0.3333\textwidth, trim= 20 45 30 40, clip]{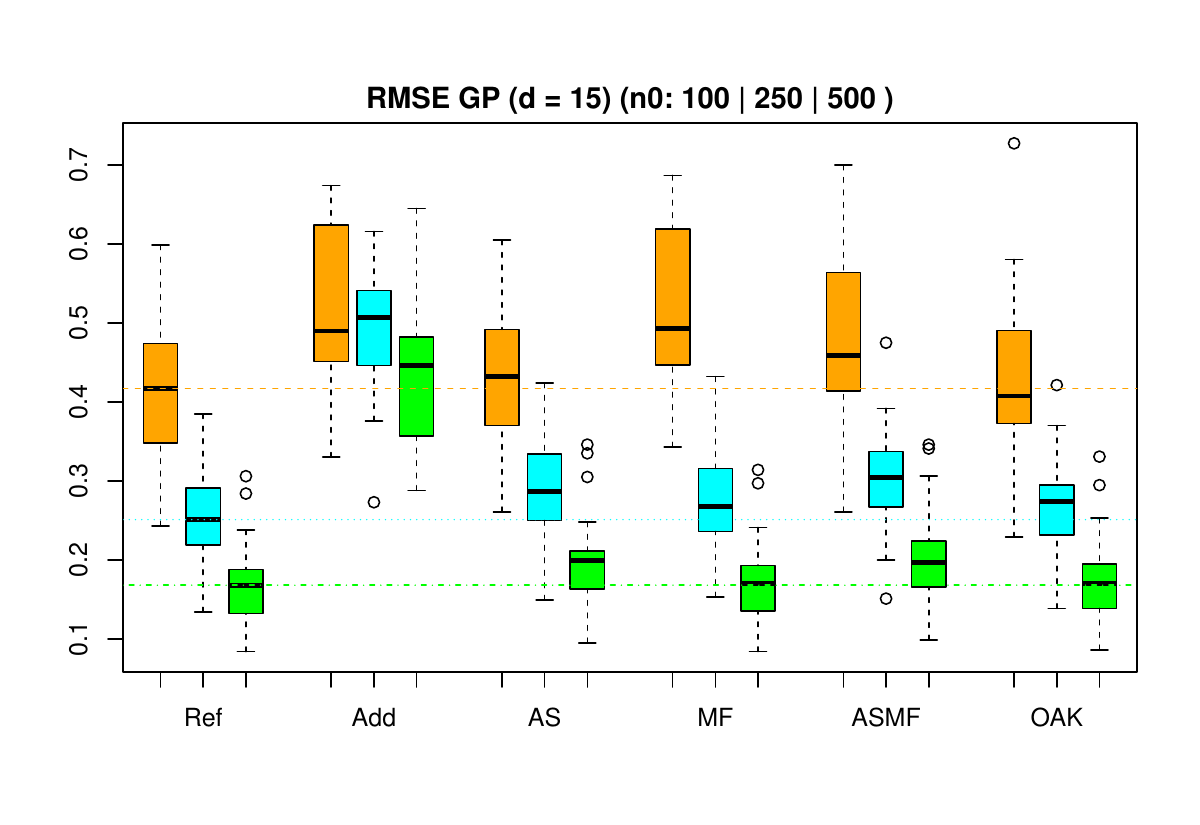}%
\includegraphics[width=0.3333\textwidth, trim= 20 45 30 40, clip]{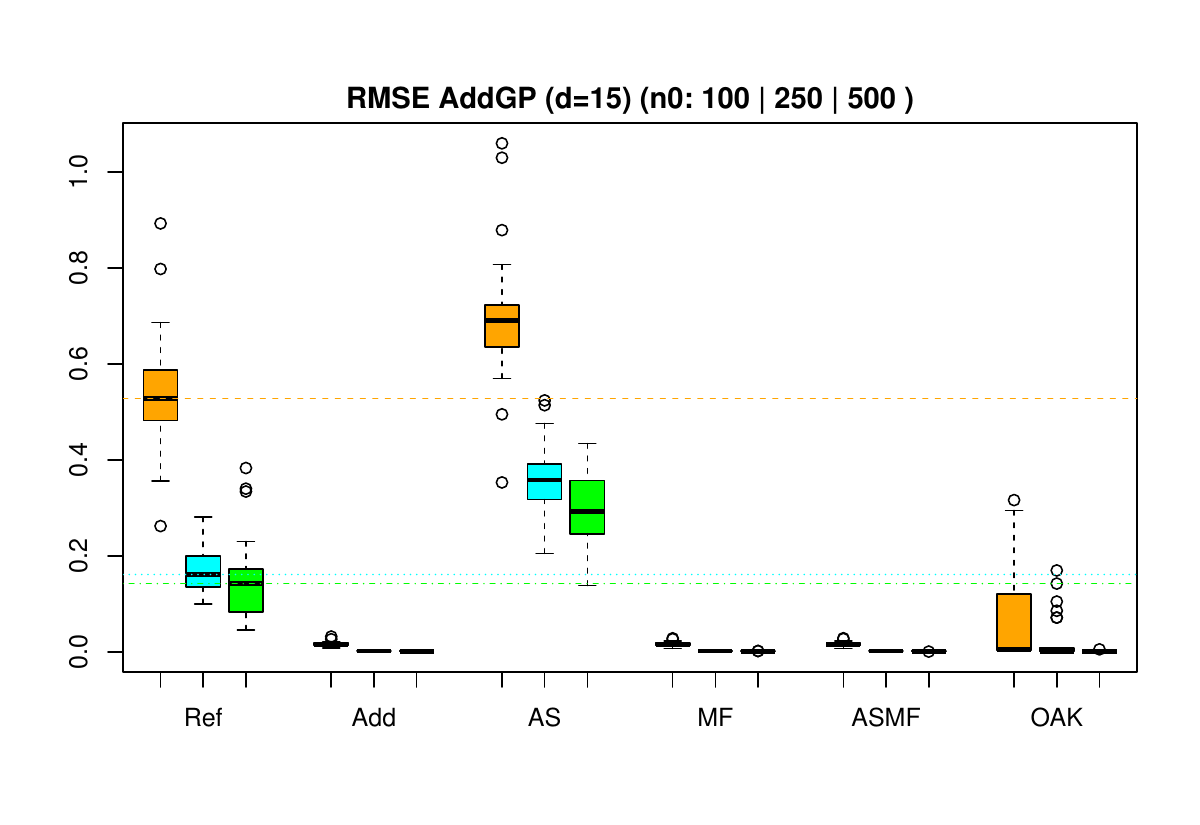}\\
\includegraphics[width=0.3333\textwidth, trim= 30 45 30 40, clip]{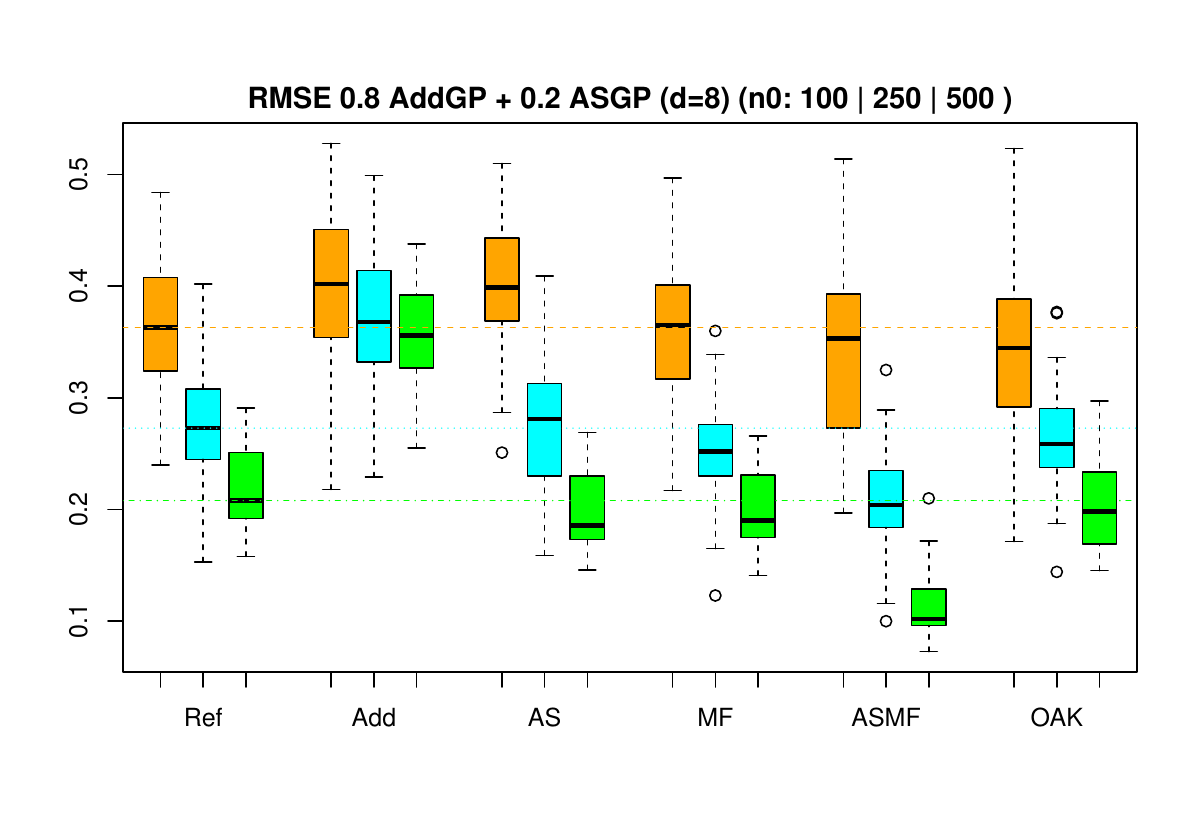}%
\includegraphics[width=0.3333\textwidth, trim= 20 45 30 40, clip]{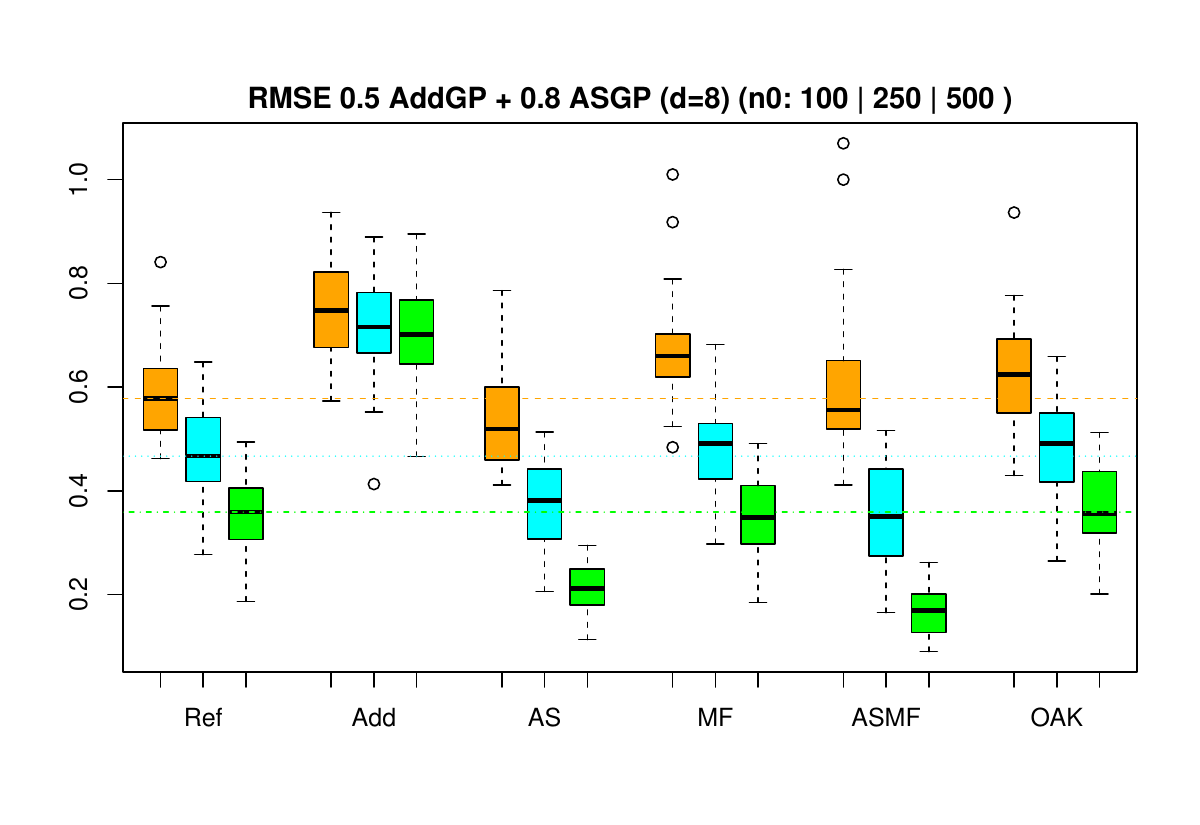}%
\includegraphics[width=0.3333\textwidth, trim= 20 45 30 40, clip]{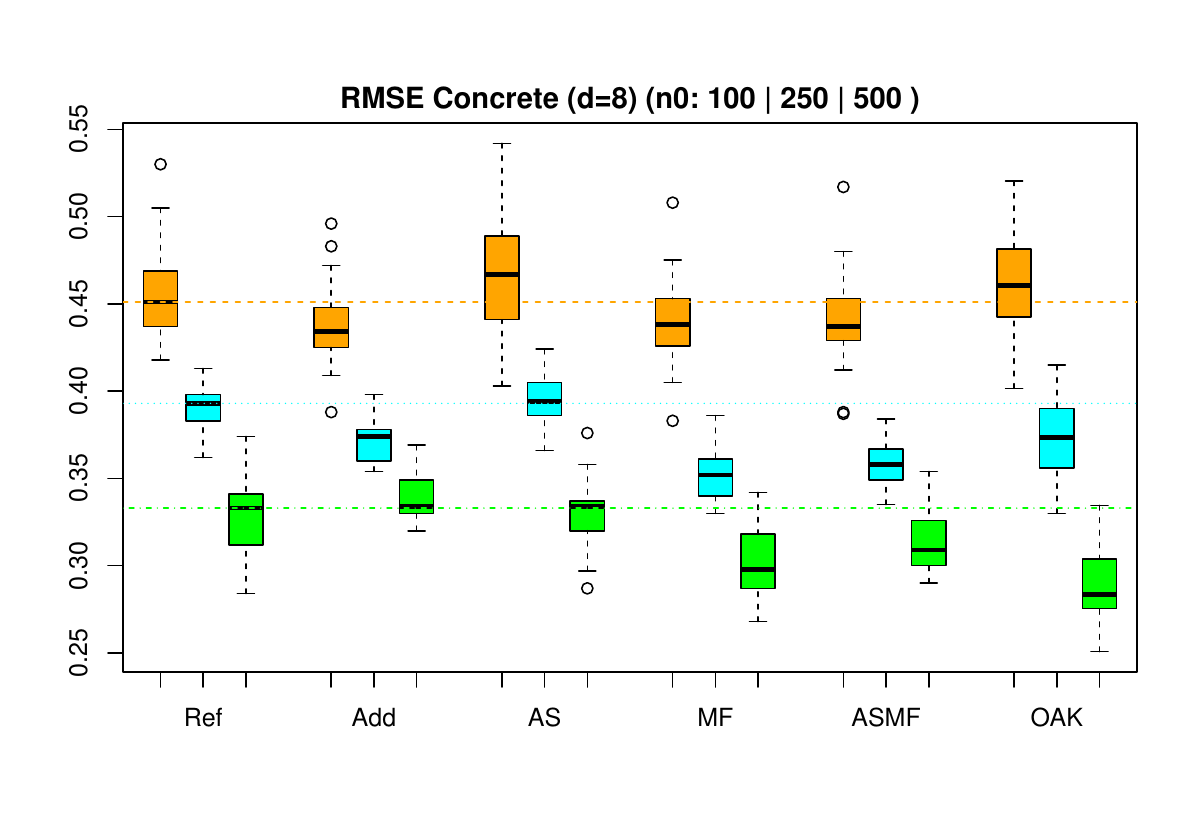}\\
\includegraphics[width=0.3333\textwidth, trim= 30 45 30 40, clip]{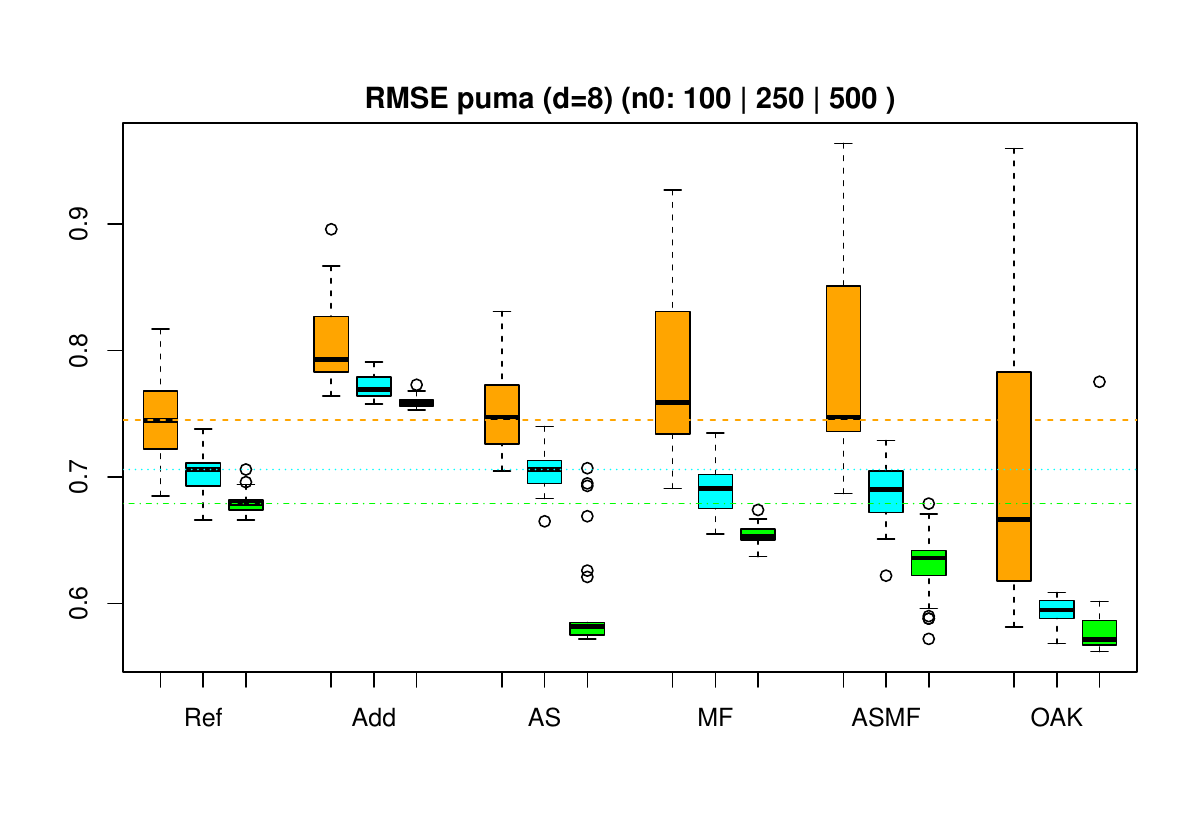}%
\includegraphics[width=0.3333\textwidth, trim= 20 45 30 40, clip]{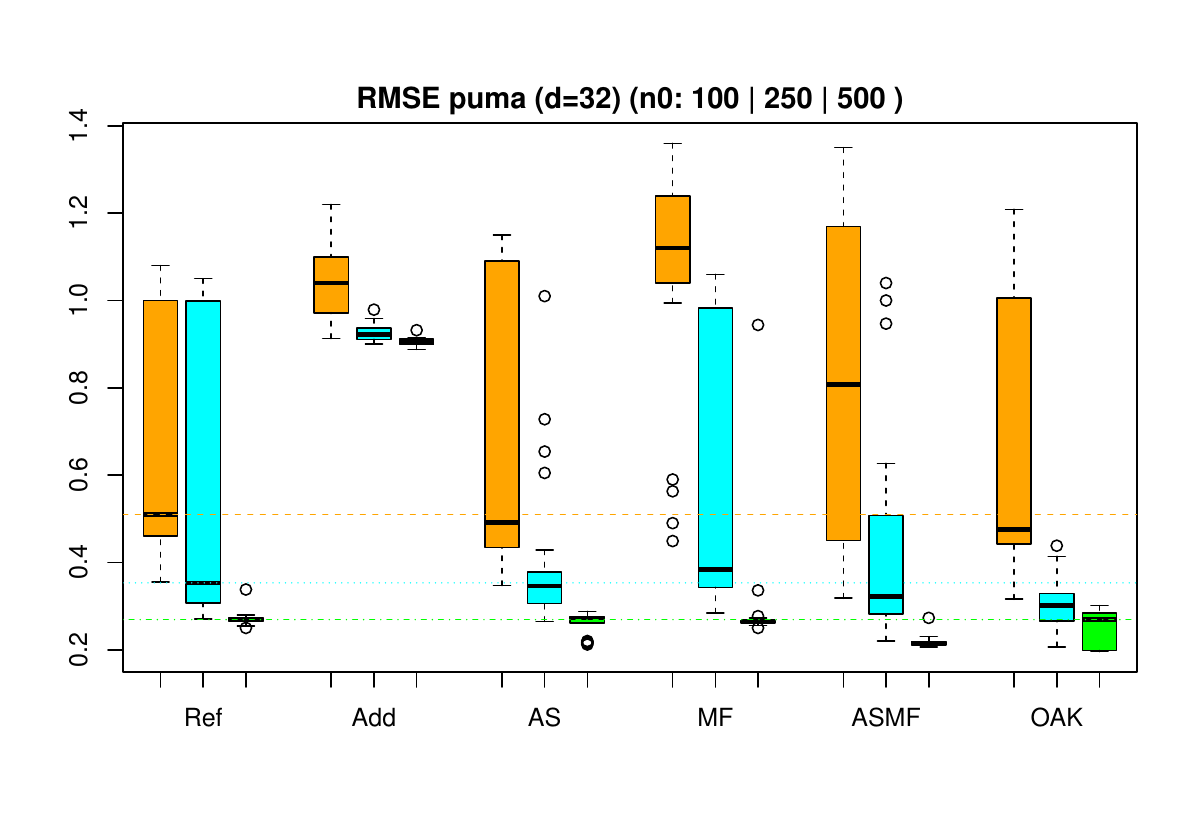}%
\includegraphics[width=0.3333\textwidth, trim= 20 45 30 40, clip]{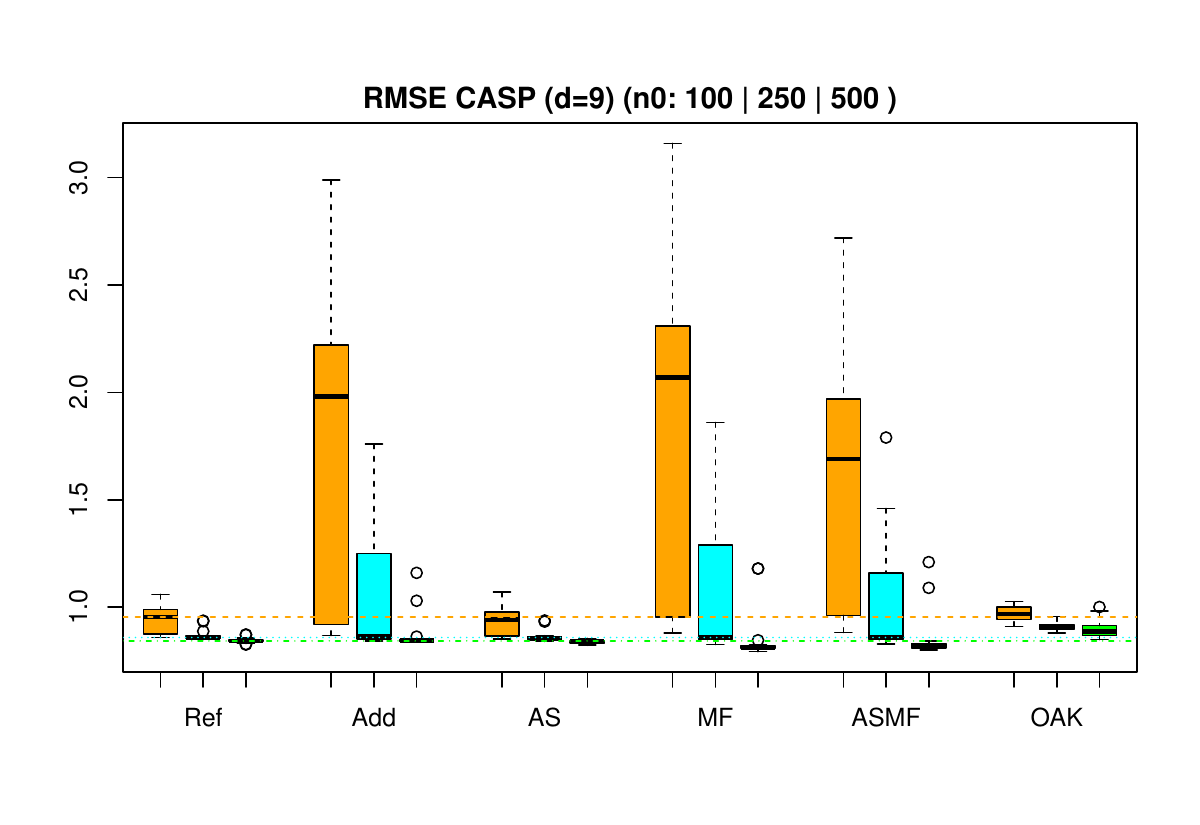}\\
\includegraphics[width=0.3333\textwidth, trim= 30 45 30 40, clip]{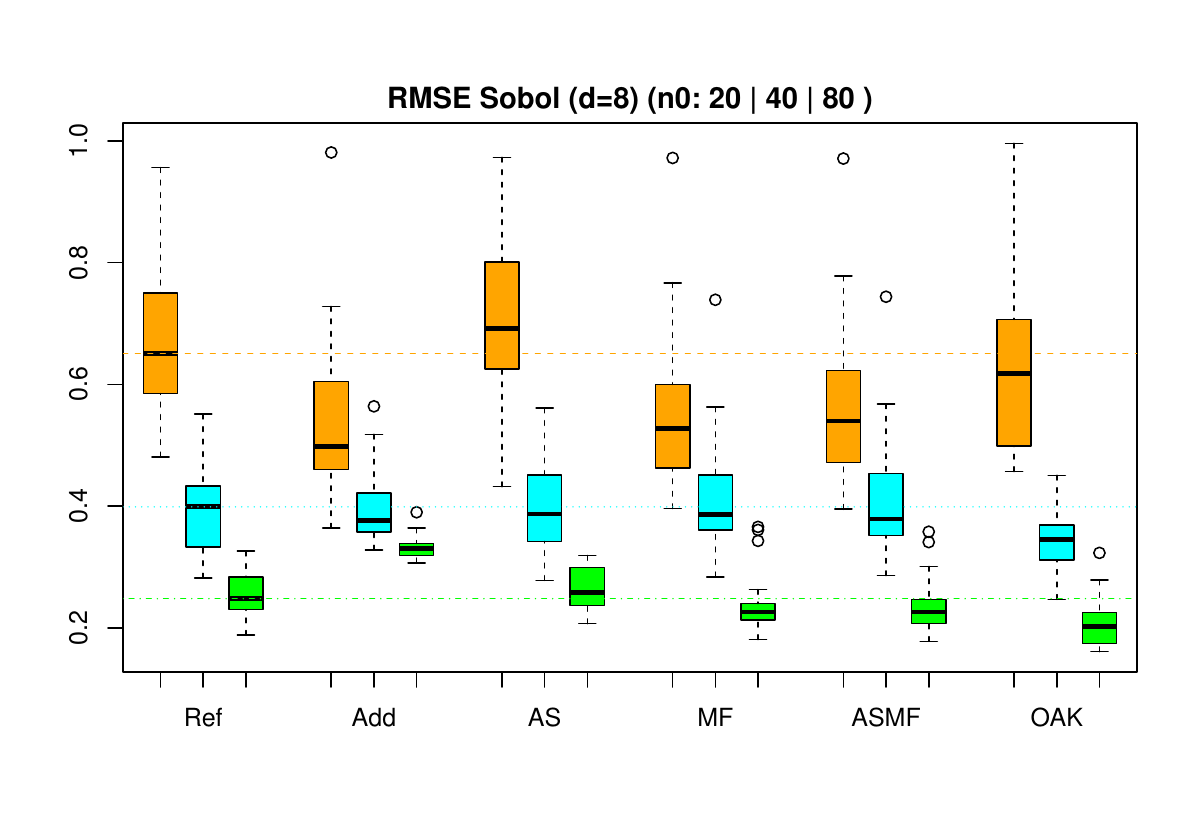}%
\includegraphics[width=0.3333\textwidth, trim= 20 45 30 40, clip]{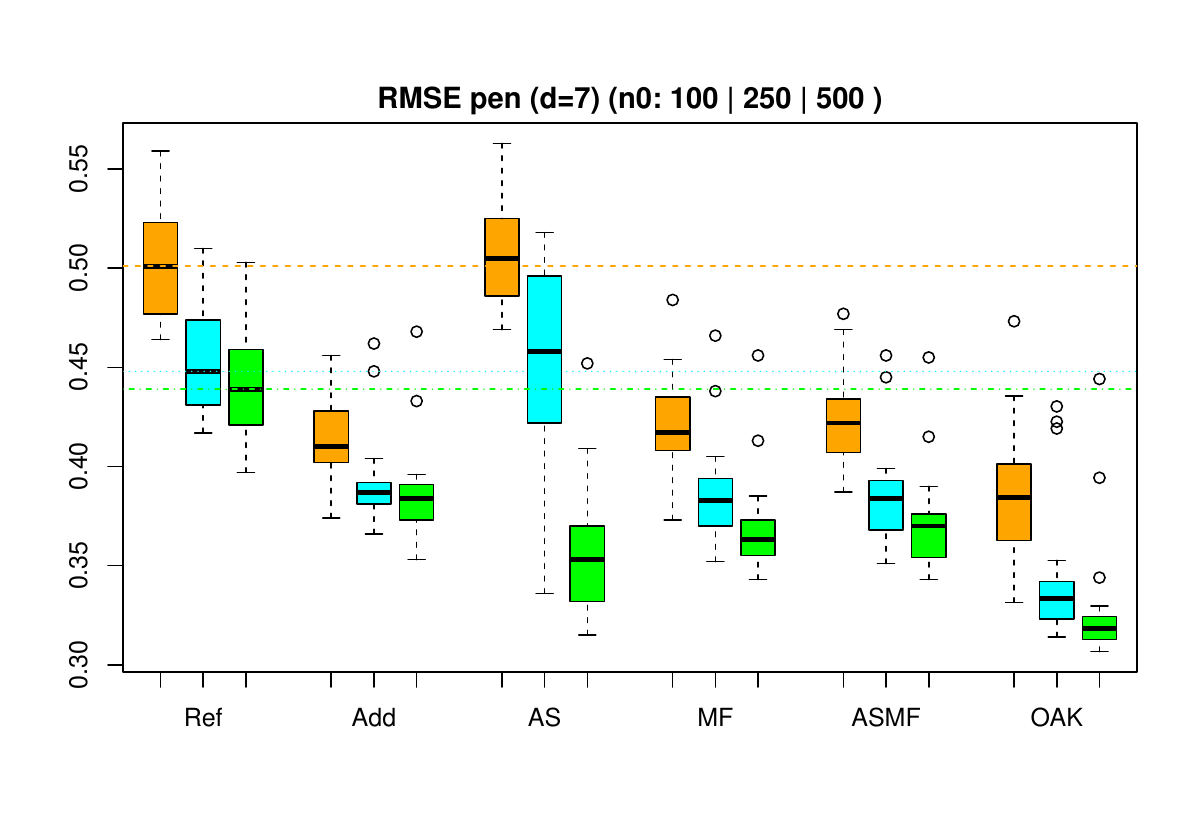}%
\includegraphics[width=0.3333\textwidth, trim= 20 45 30 40, clip]{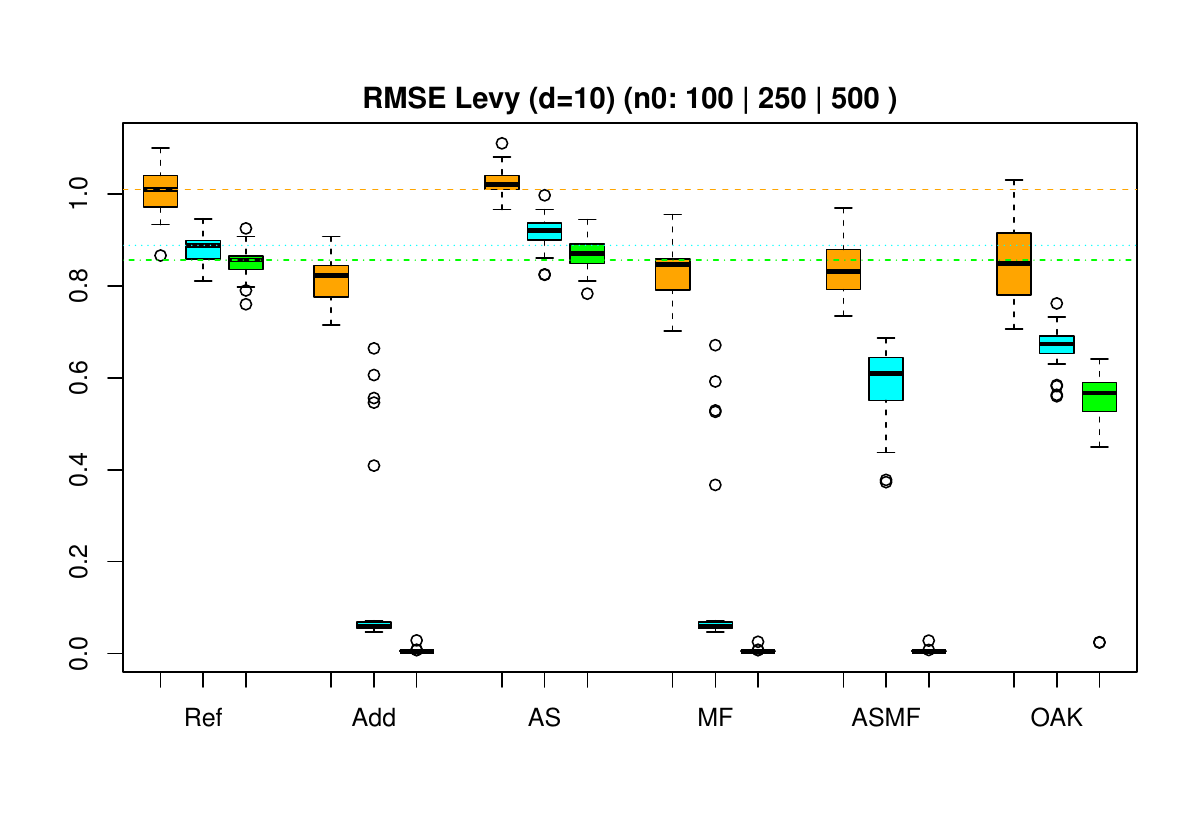}\\
\includegraphics[width=0.3333\textwidth, trim= 30 45 30 40, clip]{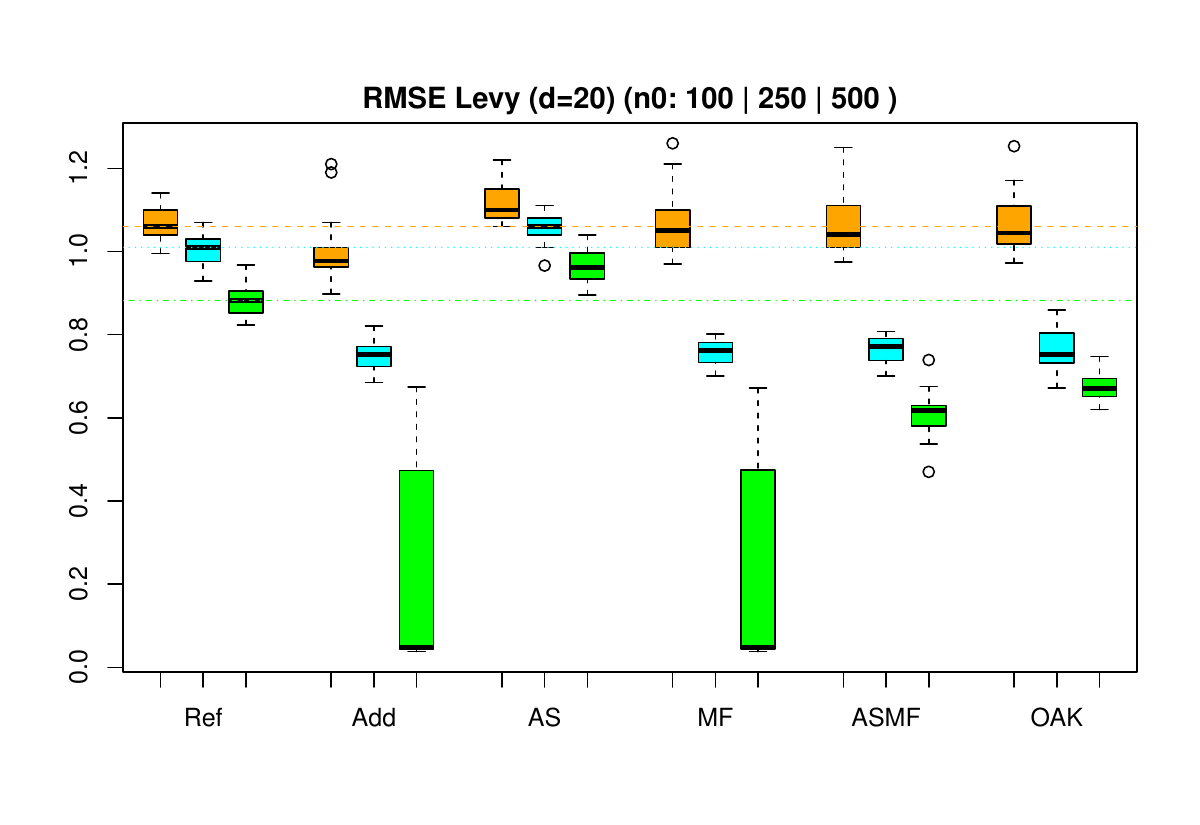}%
\includegraphics[width=0.3333\textwidth, trim= 20 45 30 40, clip]{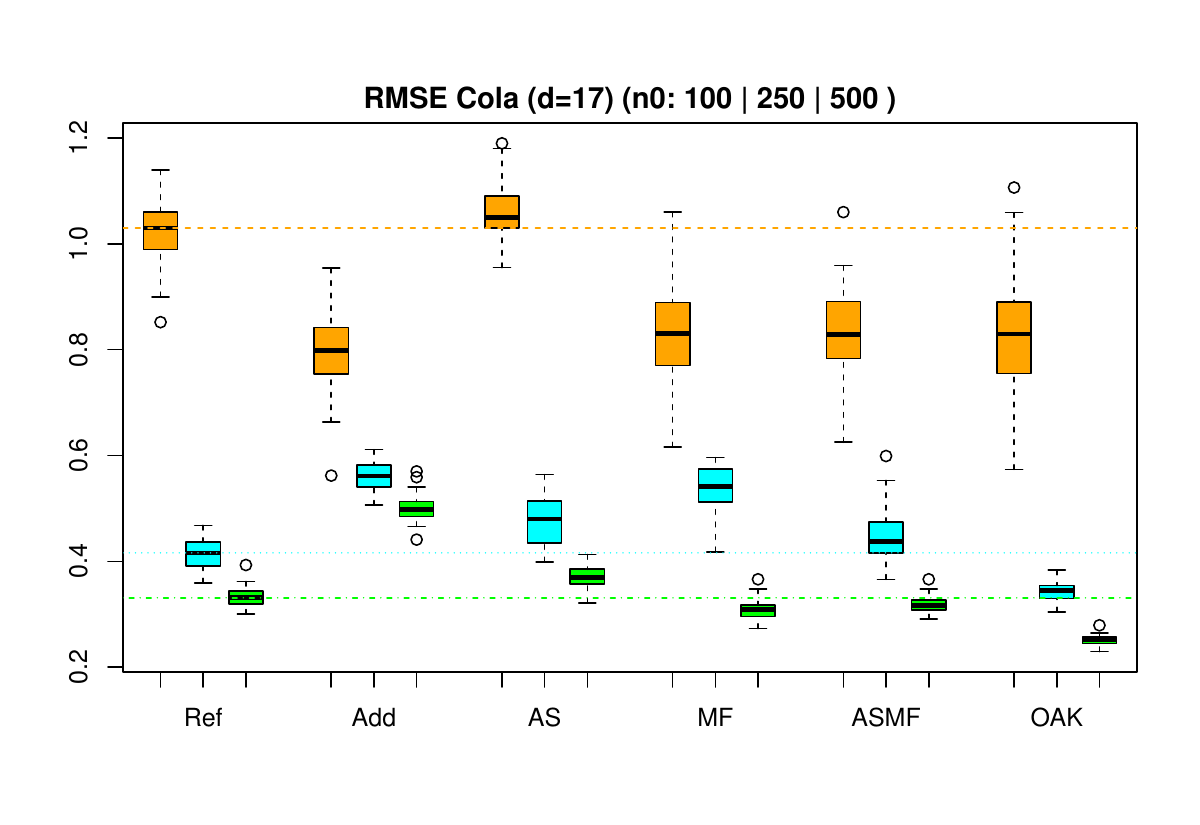}%
\includegraphics[width=0.3333\textwidth, trim= 20 45 30 40, clip]{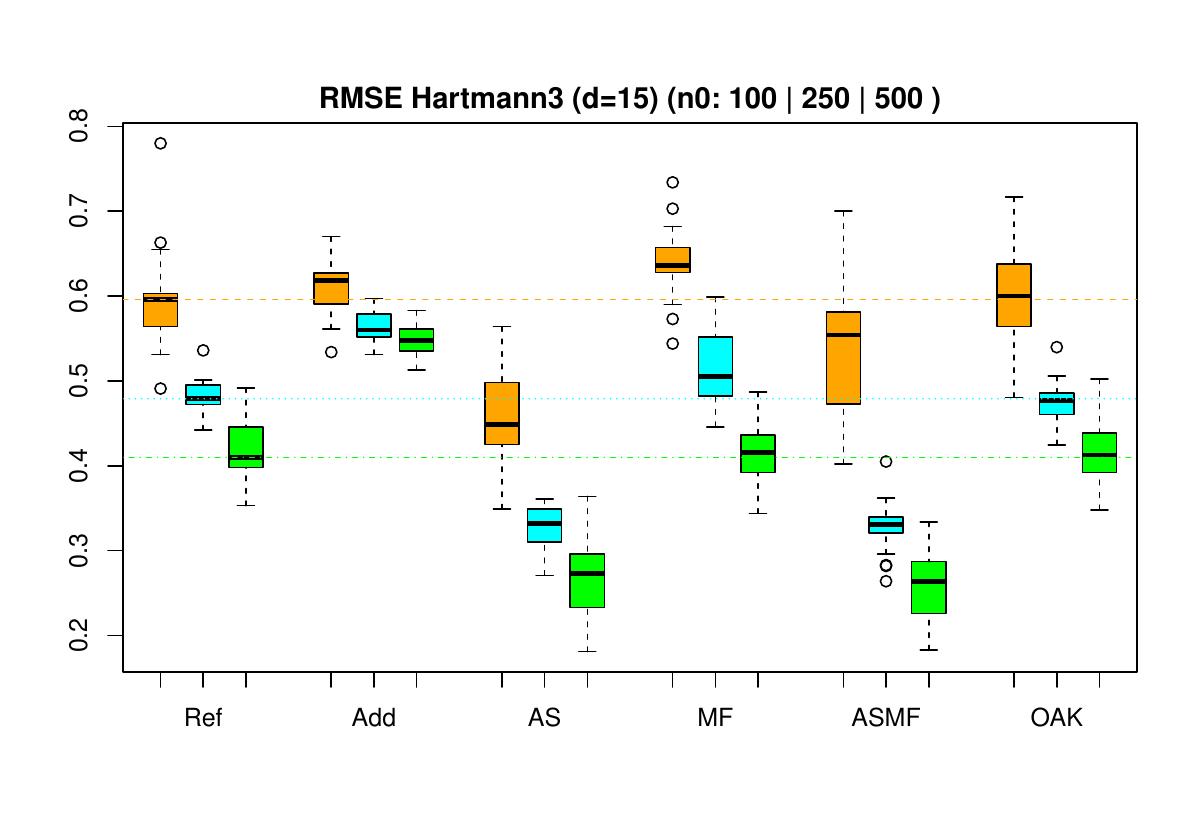}\\
\includegraphics[width=0.3333\textwidth, trim= 30 45 30 40, clip]{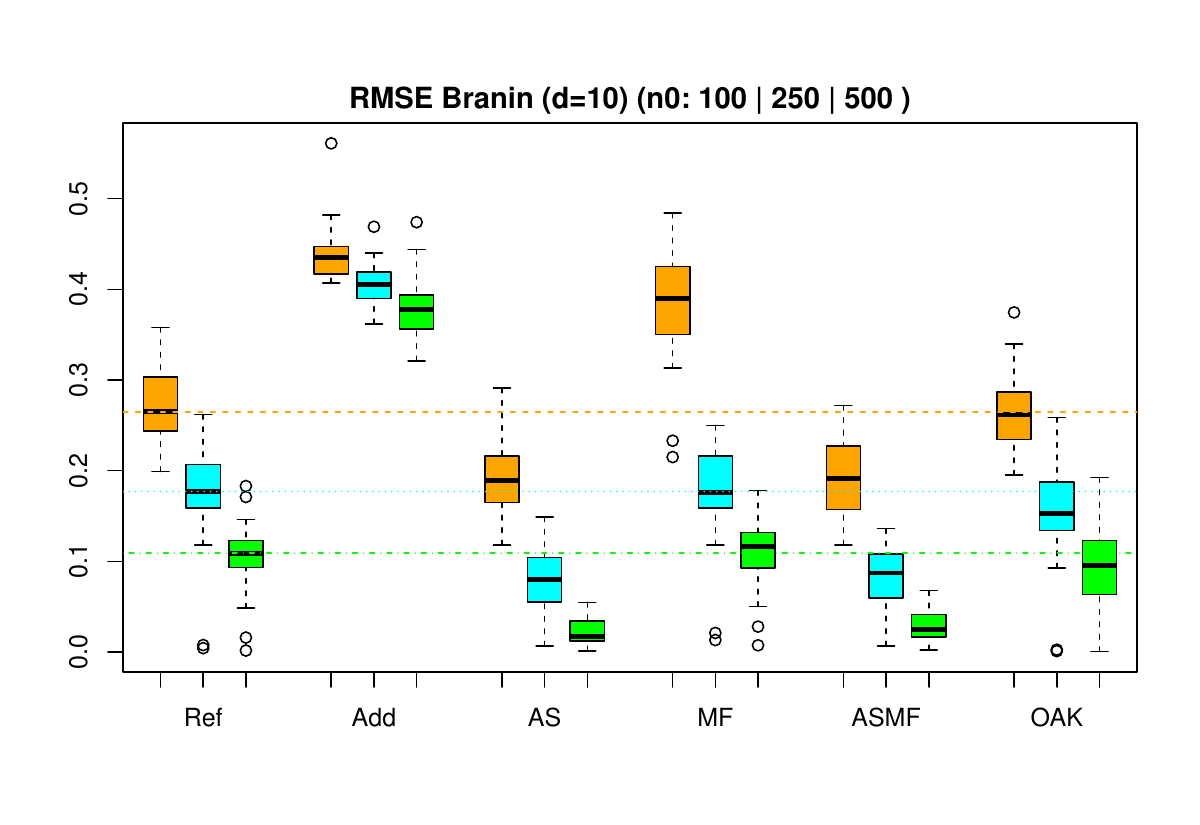}%
\includegraphics[width=0.3333\textwidth, trim= 20 45 30 40, clip]{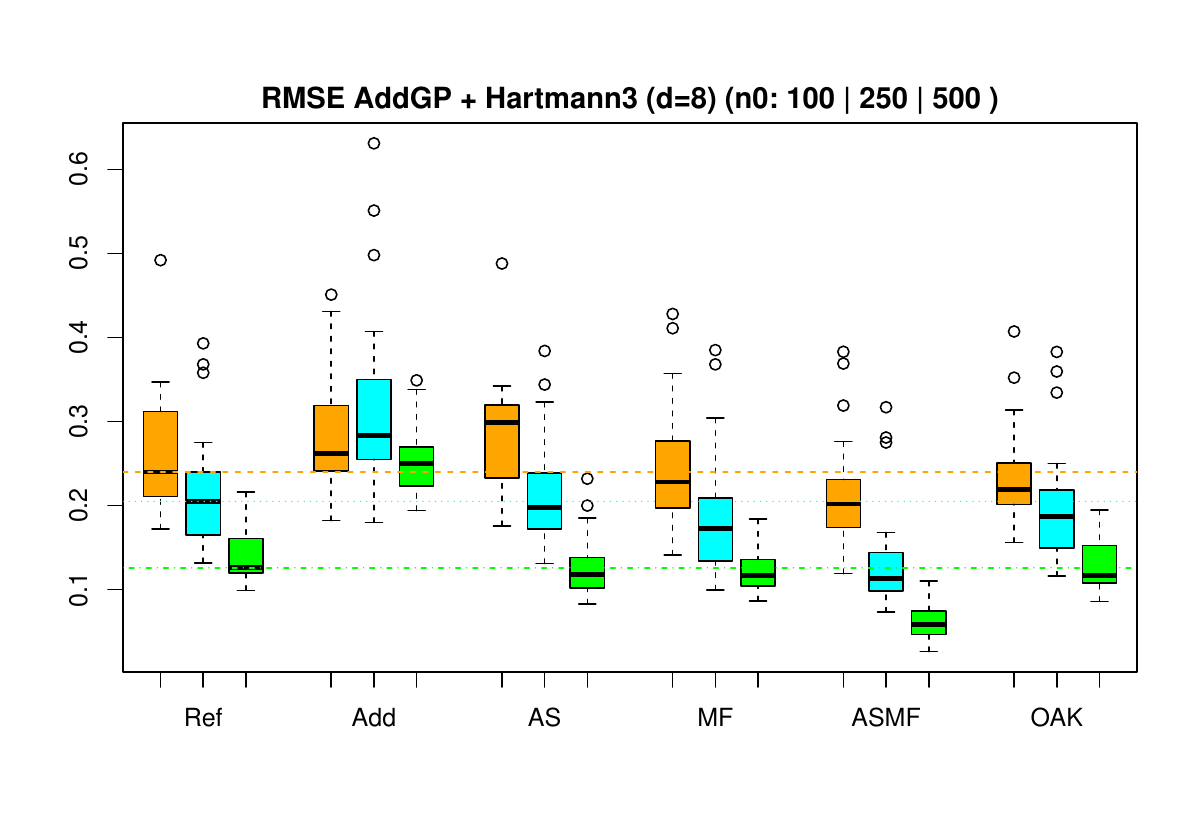}%
\includegraphics[width=0.3333\textwidth, trim= 20 45 30 40, clip]{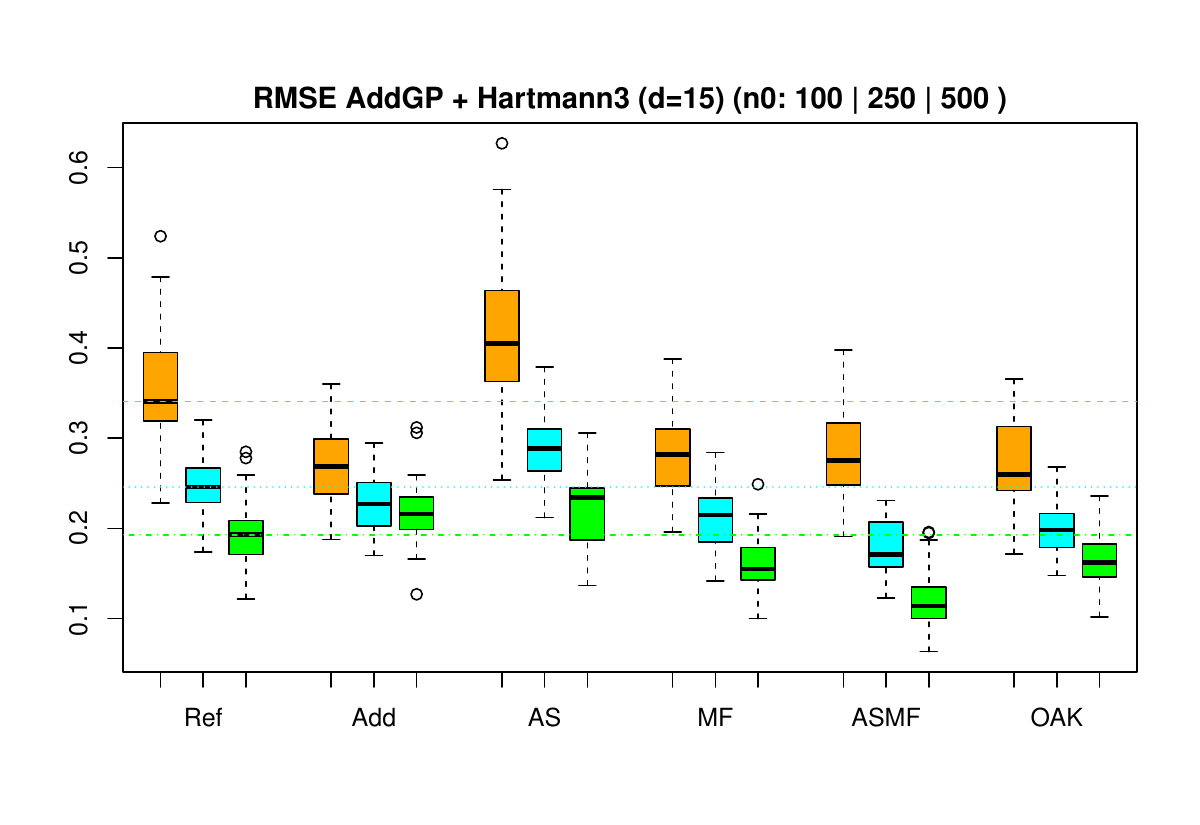}
\caption{Additional RMSE results including the OAK model. The different $n_0$
 values appear with different colors, the corresponding lines indicate the
 baseline result from standard GP models.}
\label{fig:supprmseres}
\end{figure*}

\begin{figure*}[htpb]
\centering
\includegraphics[width=0.3333\textwidth, trim= 30 45 30 40, clip]{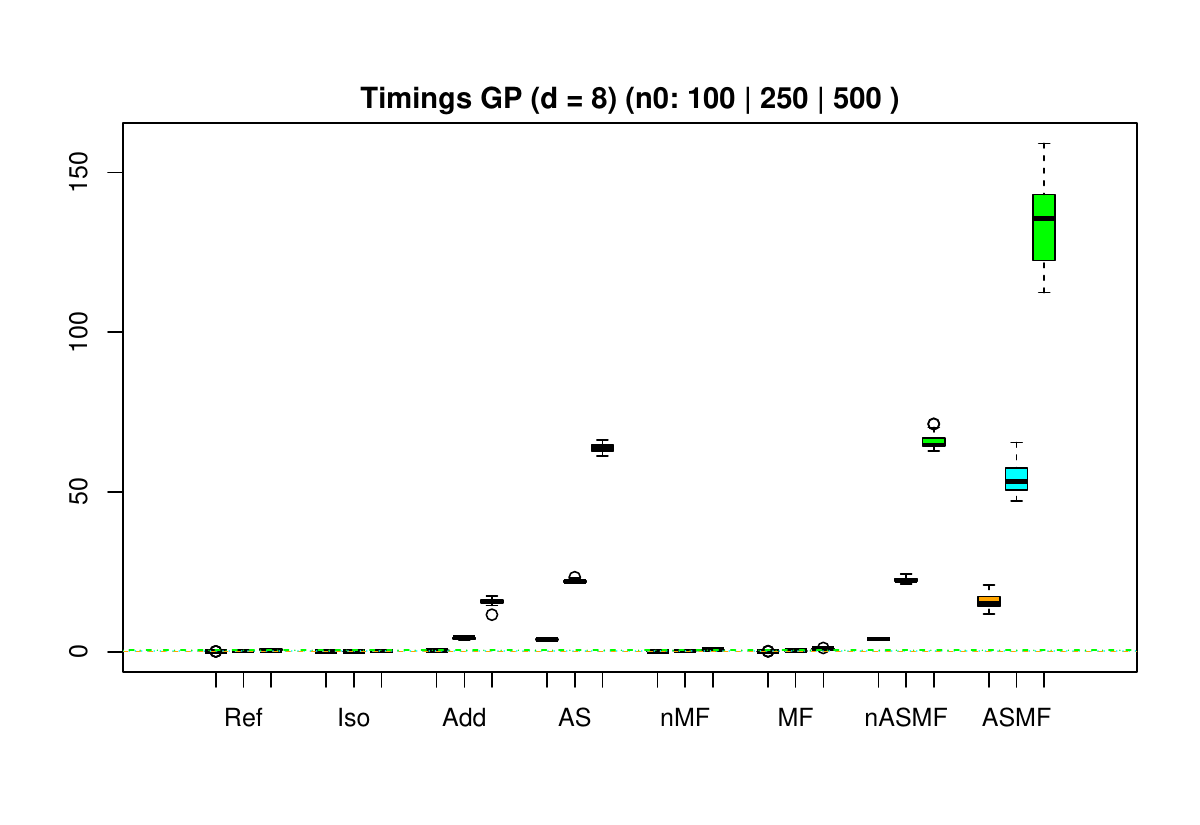}%
\includegraphics[width=0.3333\textwidth, trim= 20 45 30 40, clip]{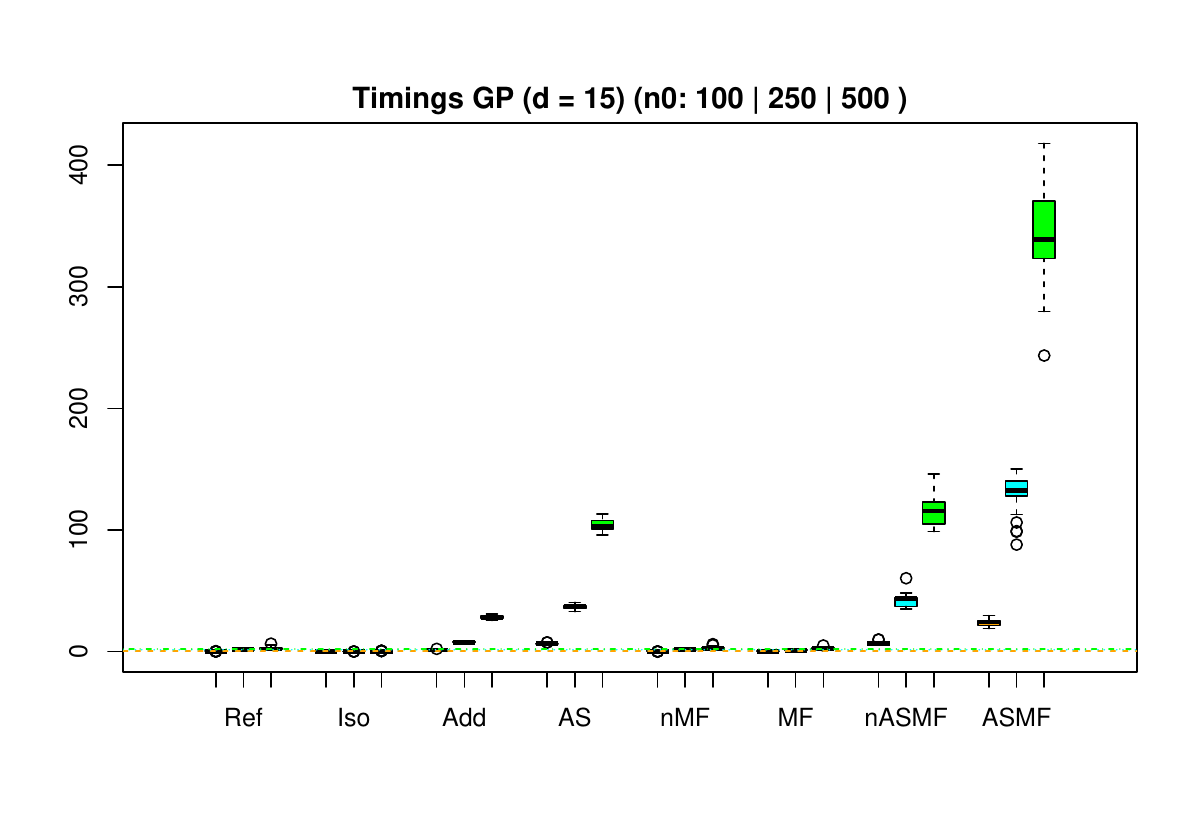}%
\includegraphics[width=0.3333\textwidth, trim= 20 45 30 40, clip]{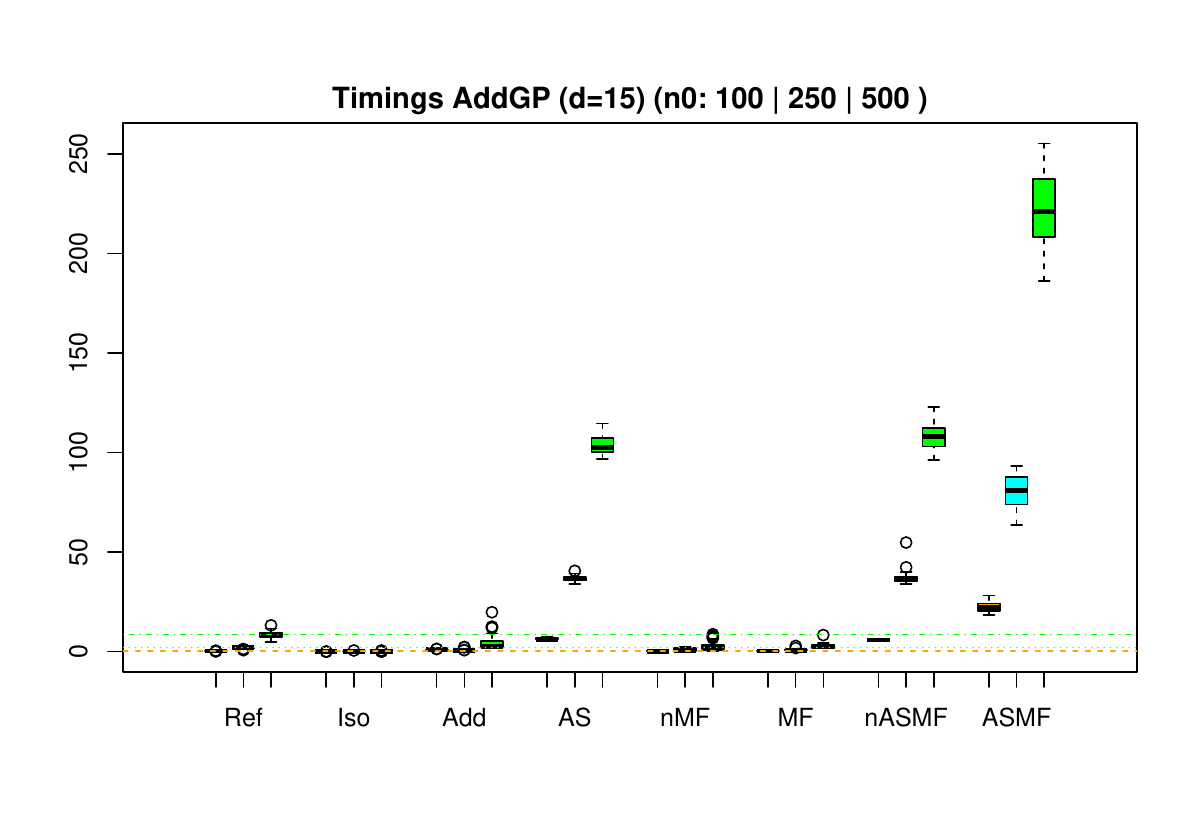}\\
\includegraphics[width=0.3333\textwidth, trim= 30 45 30 40, clip]{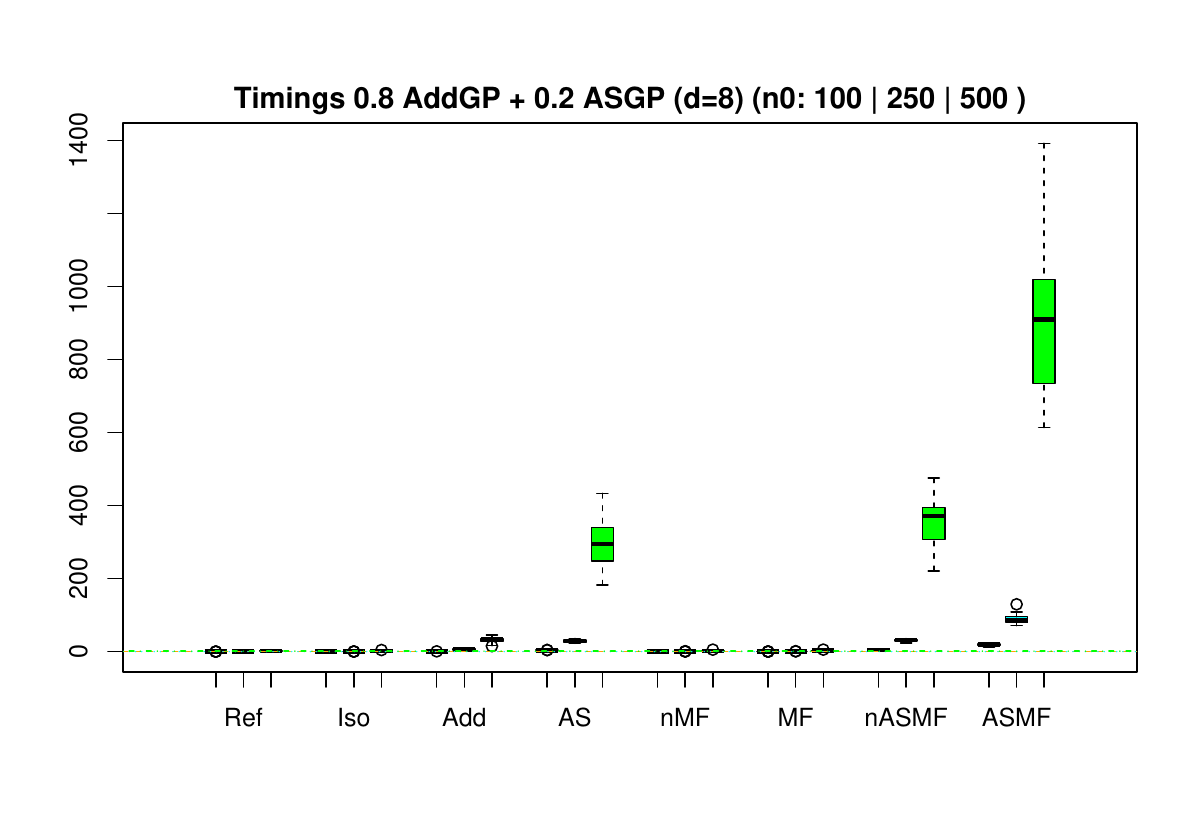}%
\includegraphics[width=0.3333\textwidth, trim= 20 45 30 40, clip]{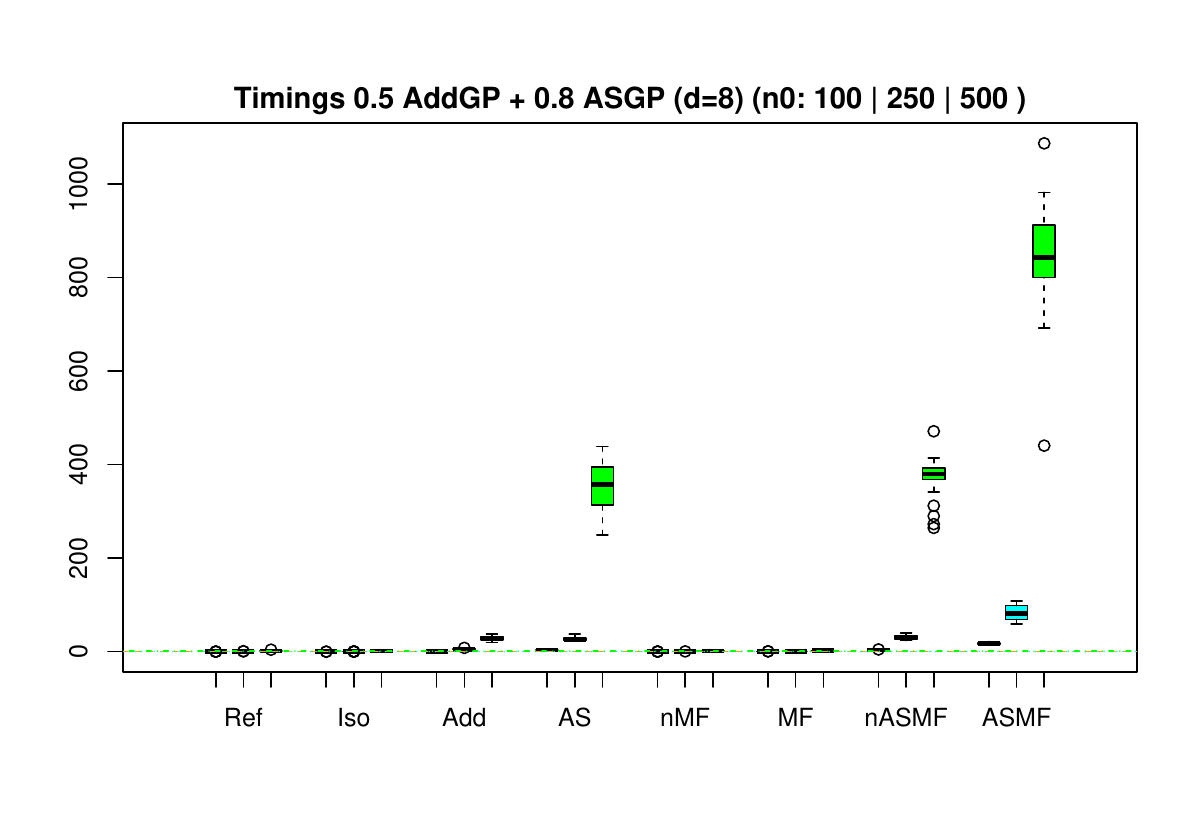}%
\includegraphics[width=0.3333\textwidth, trim= 20 45 30 40, clip]{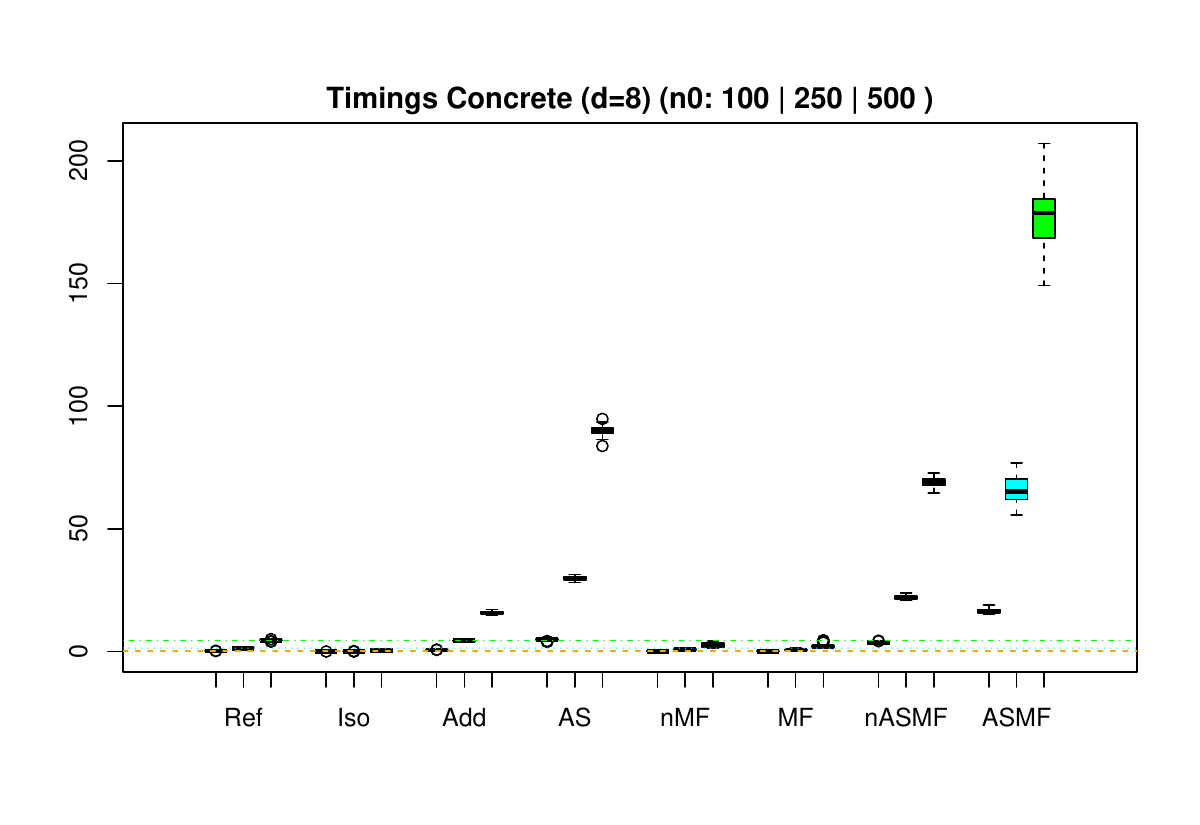}\\
\includegraphics[width=0.3333\textwidth, trim= 30 45 30 40, clip]{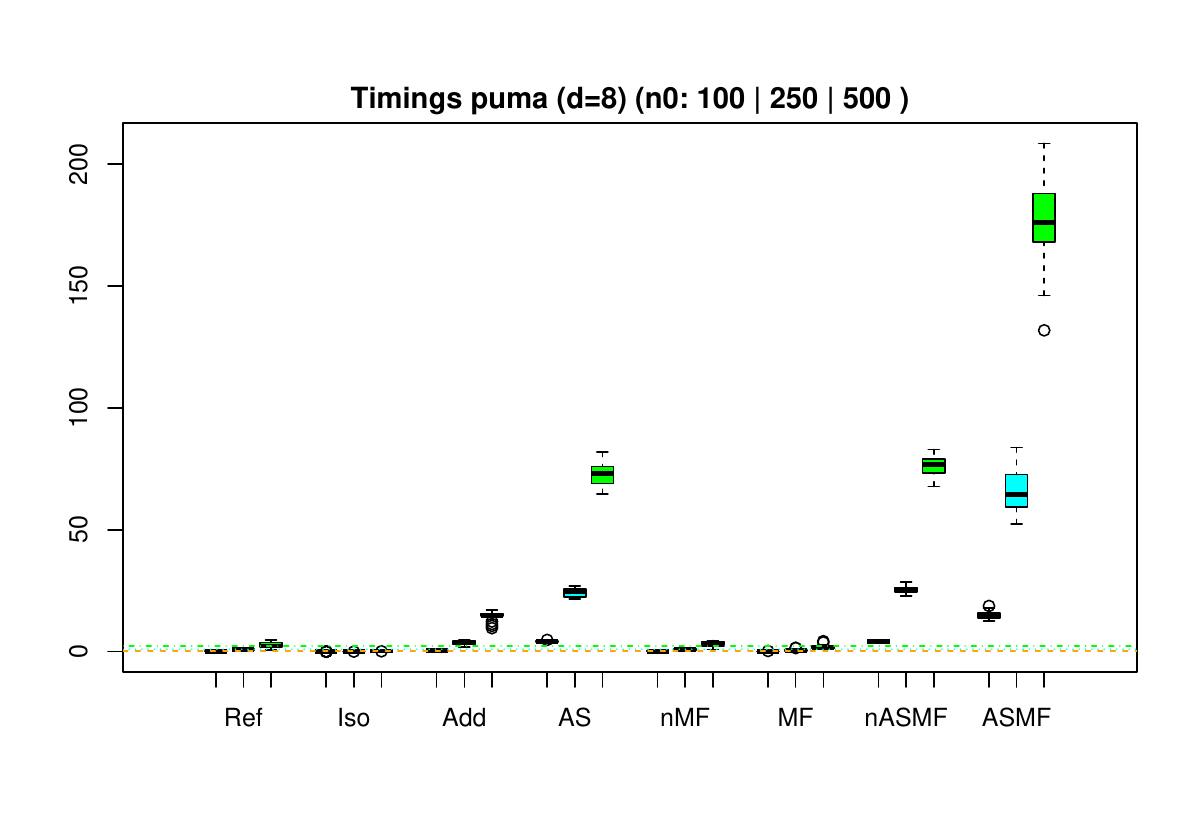}%
\includegraphics[width=0.3333\textwidth, trim= 20 45 30 40, clip]{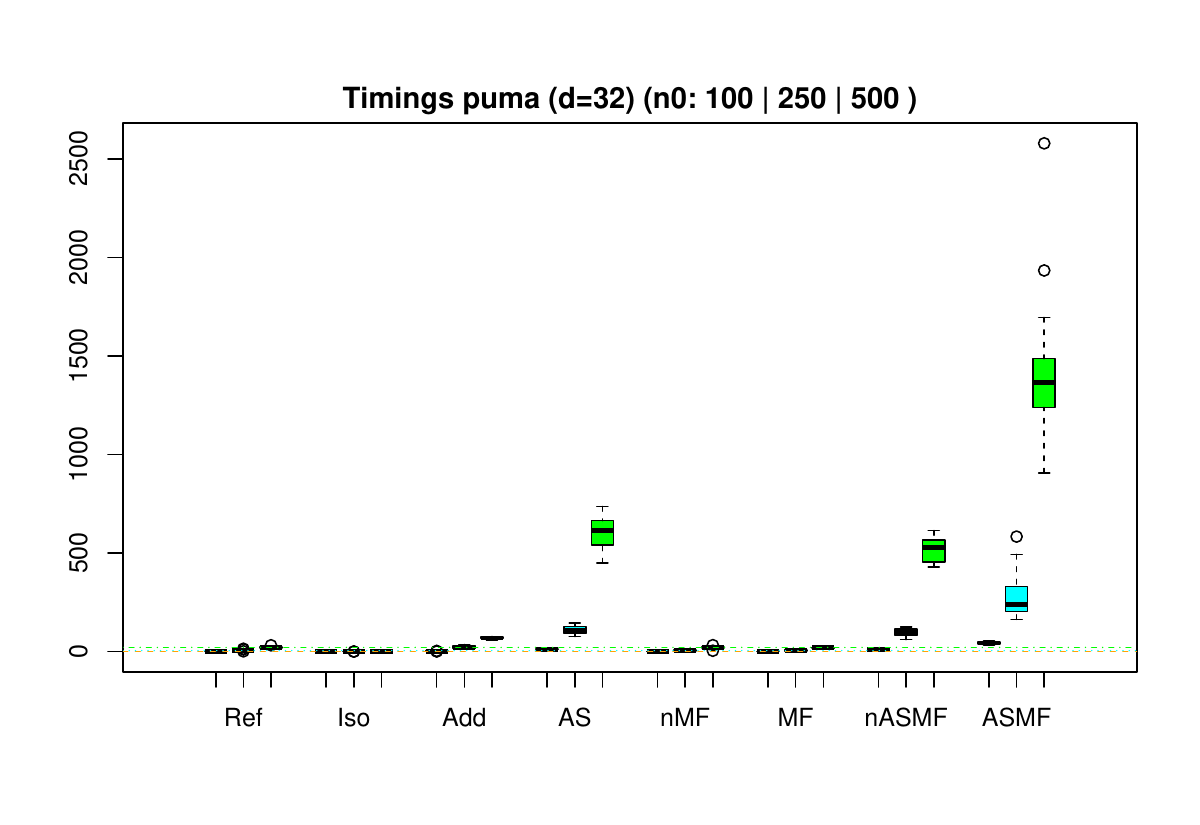}%
\includegraphics[width=0.3333\textwidth, trim= 20 45 30 40, clip]{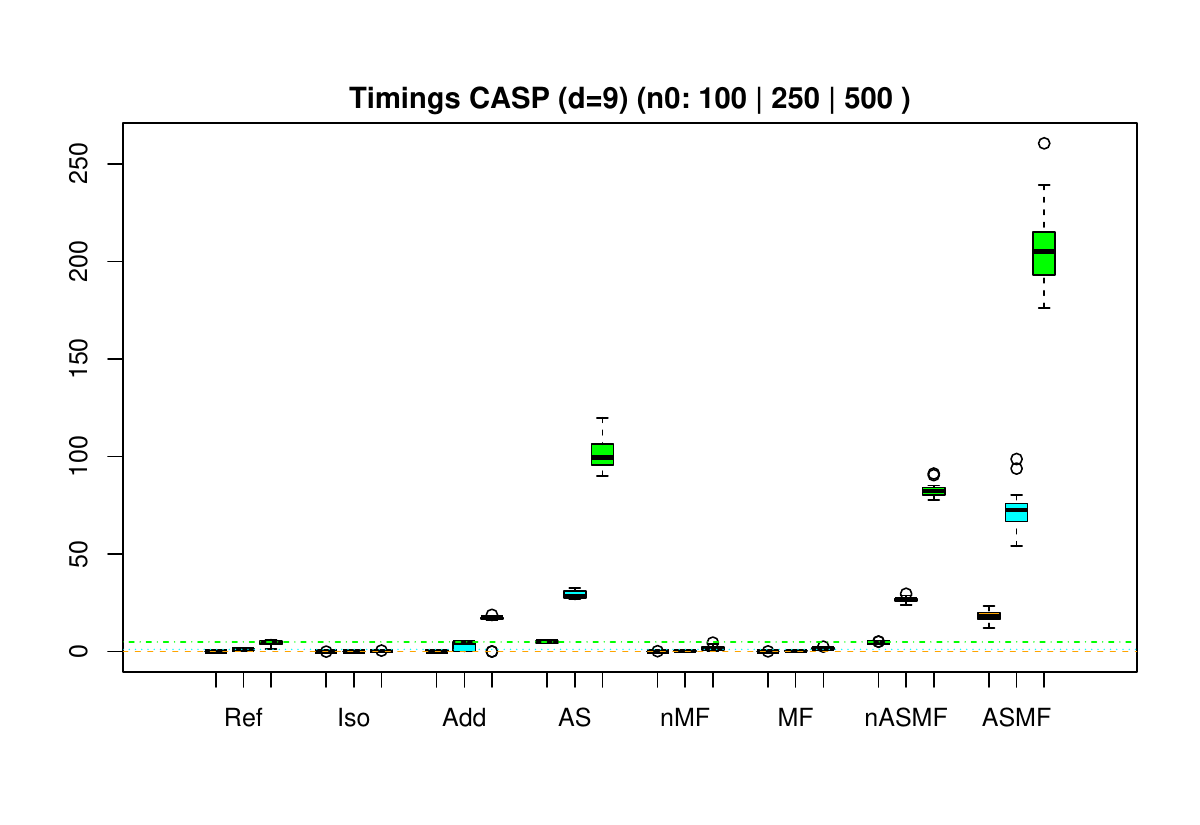}\\
\includegraphics[width=0.3333\textwidth, trim= 30 45 30 40, clip]{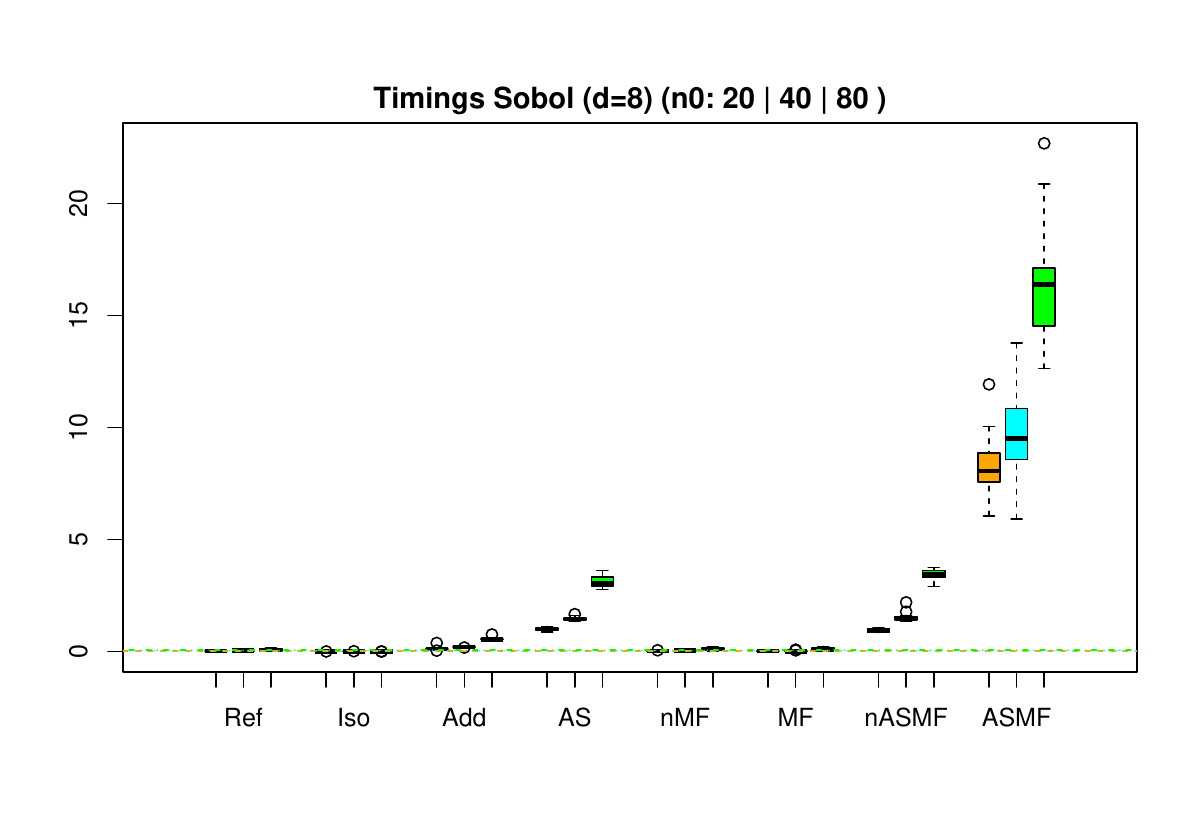}%
\includegraphics[width=0.3333\textwidth, trim= 20 45 30 40, clip]{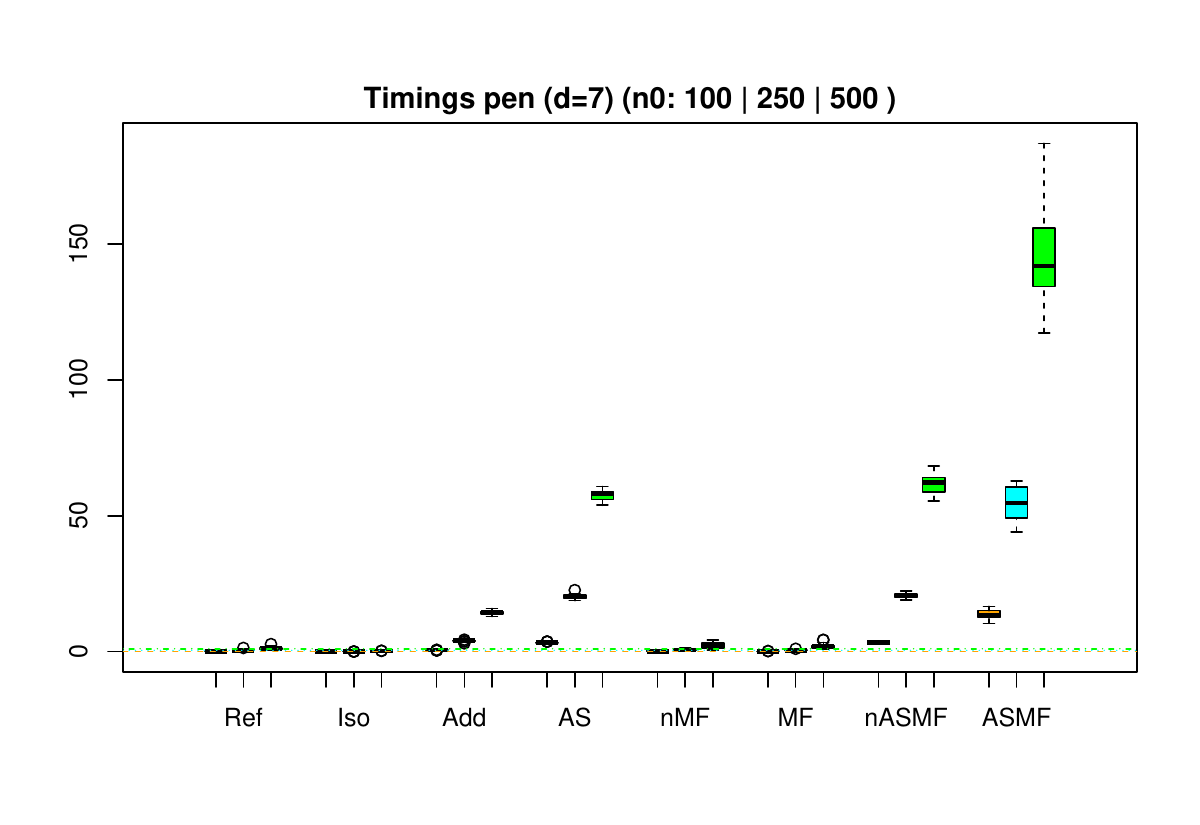}%
\includegraphics[width=0.3333\textwidth, trim= 20 45 30 40, clip]{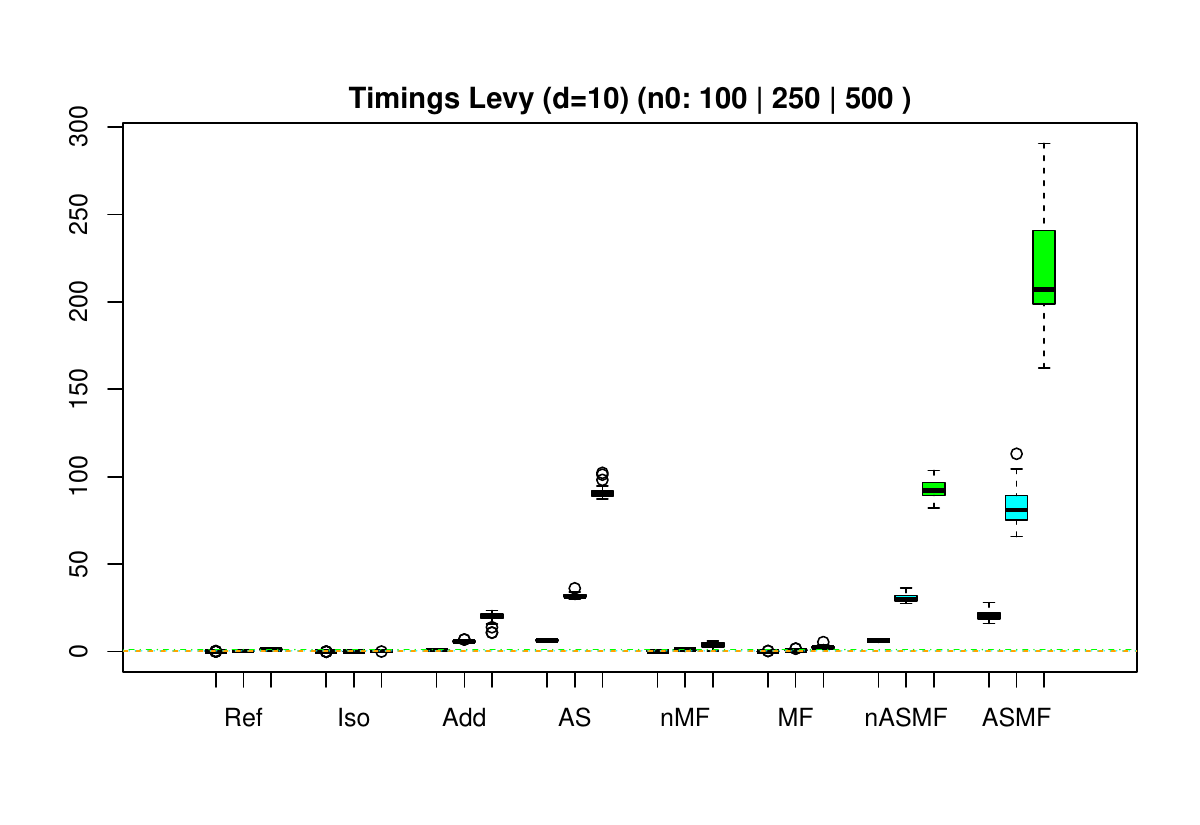}\\
\includegraphics[width=0.3333\textwidth, trim= 30 45 30 40, clip]{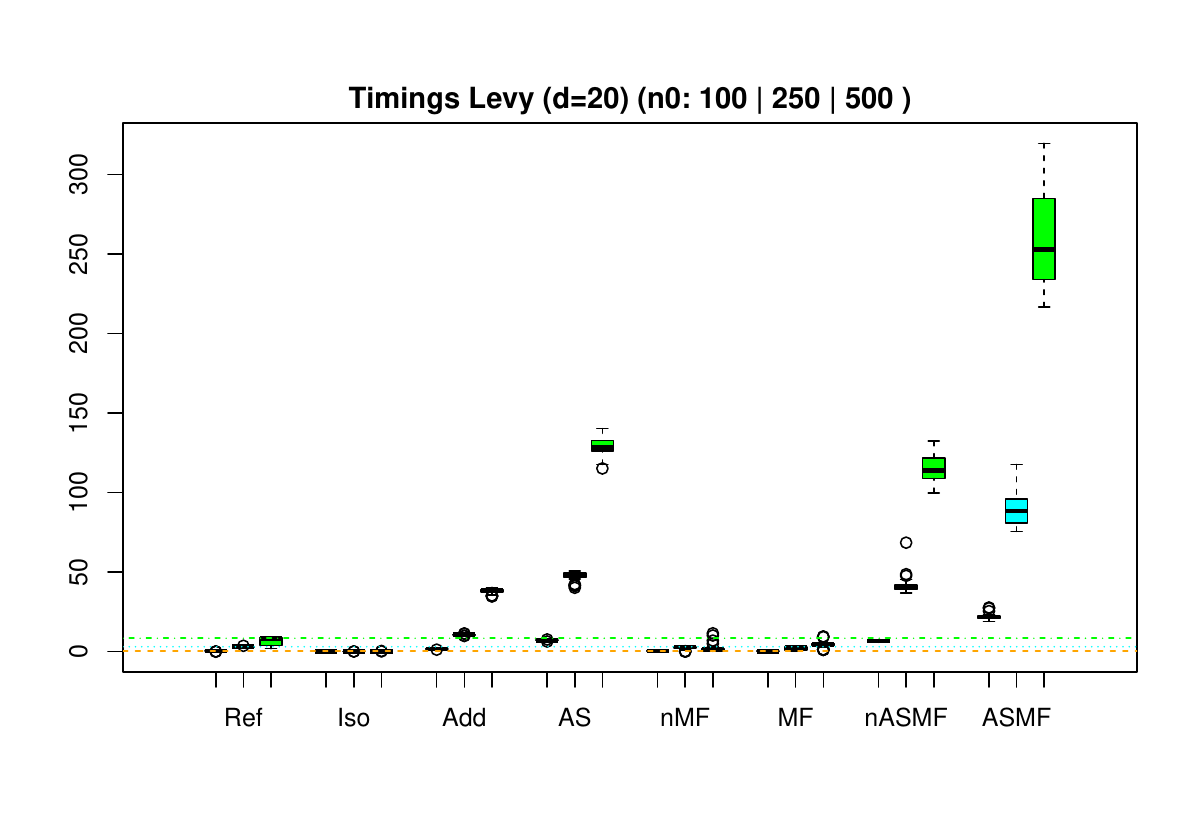}%
\includegraphics[width=0.3333\textwidth, trim= 20 45 30 40, clip]{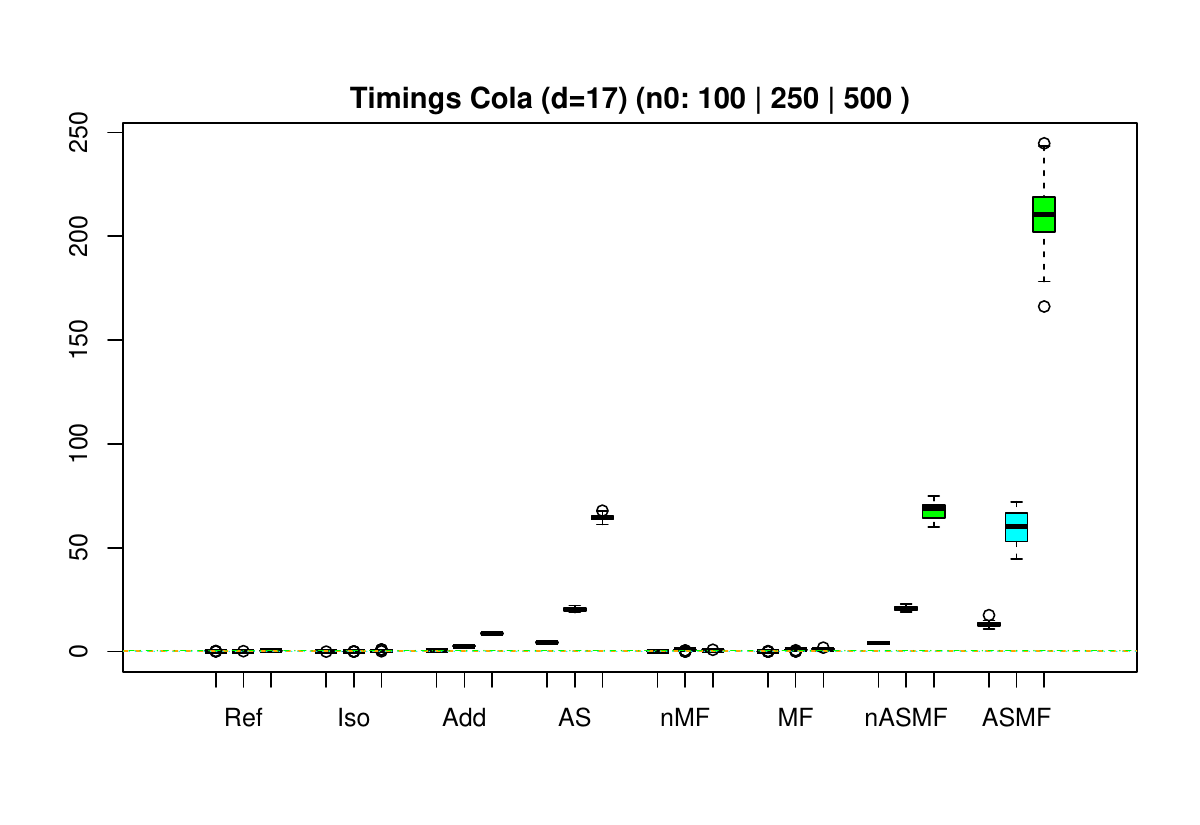}%
\includegraphics[width=0.3333\textwidth, trim= 20 45 30 40, clip]{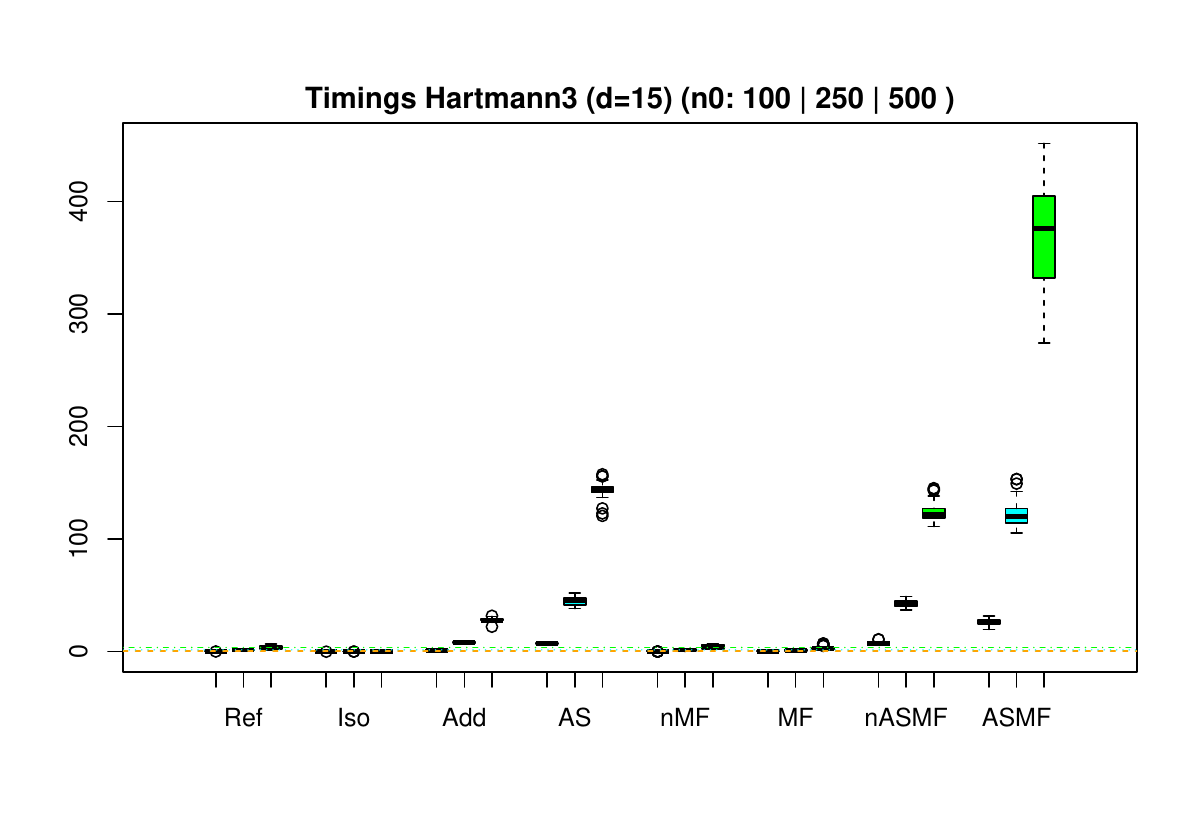}\\
\includegraphics[width=0.3333\textwidth, trim= 30 45 30 40, clip]{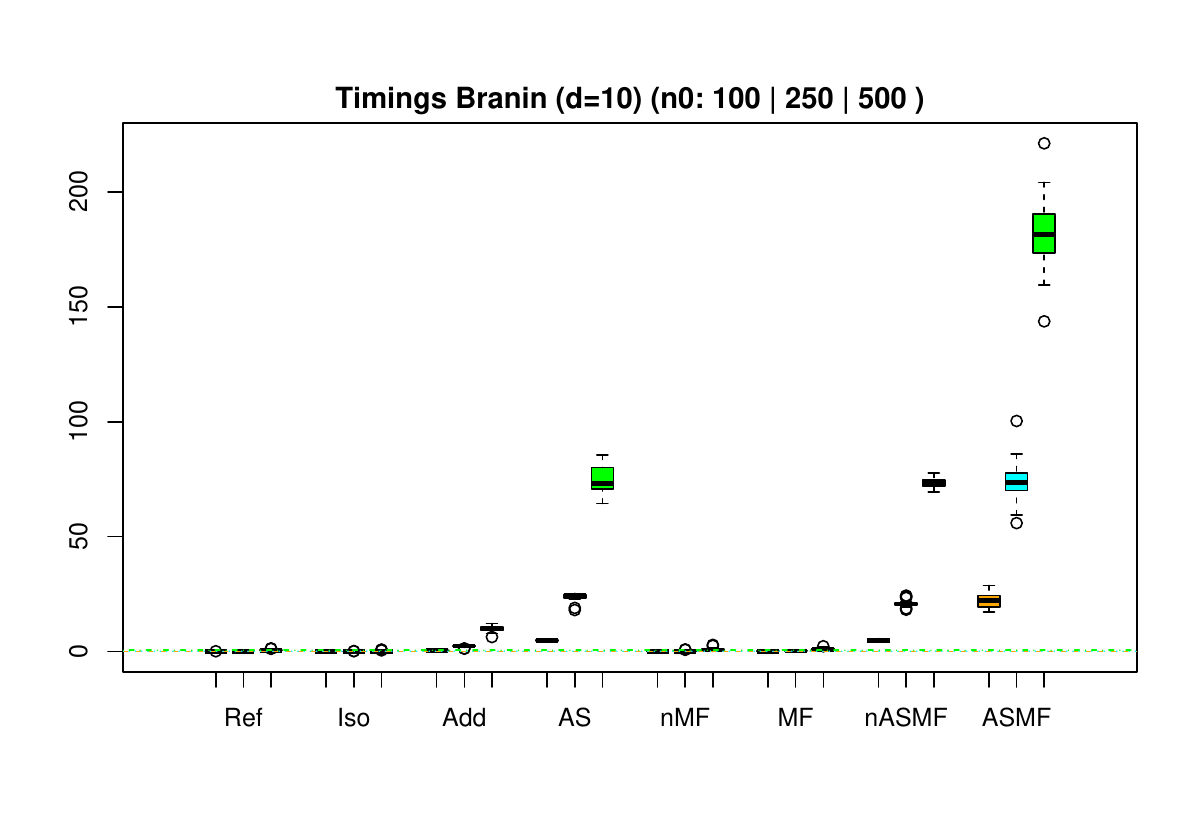}%
\includegraphics[width=0.3333\textwidth, trim= 20 45 30 40, clip]{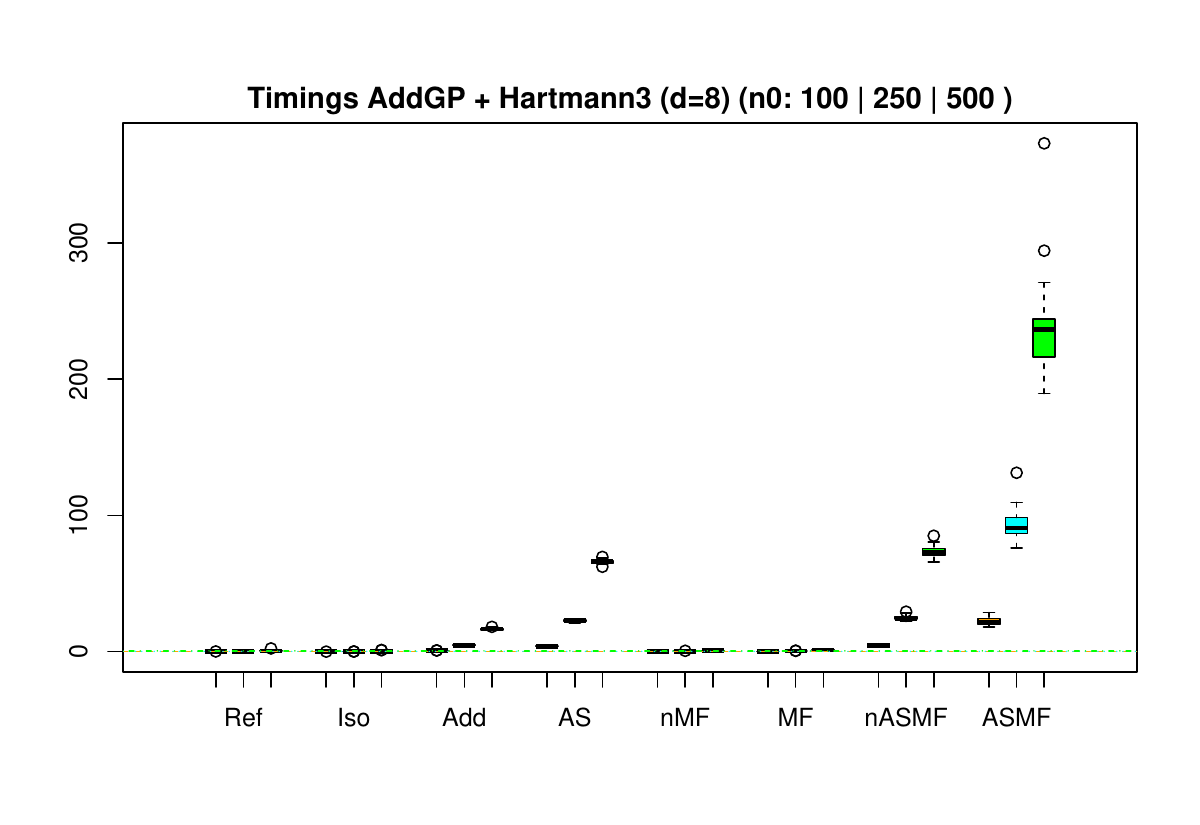}%
\includegraphics[width=0.3333\textwidth, trim= 20 45 30 40, clip]{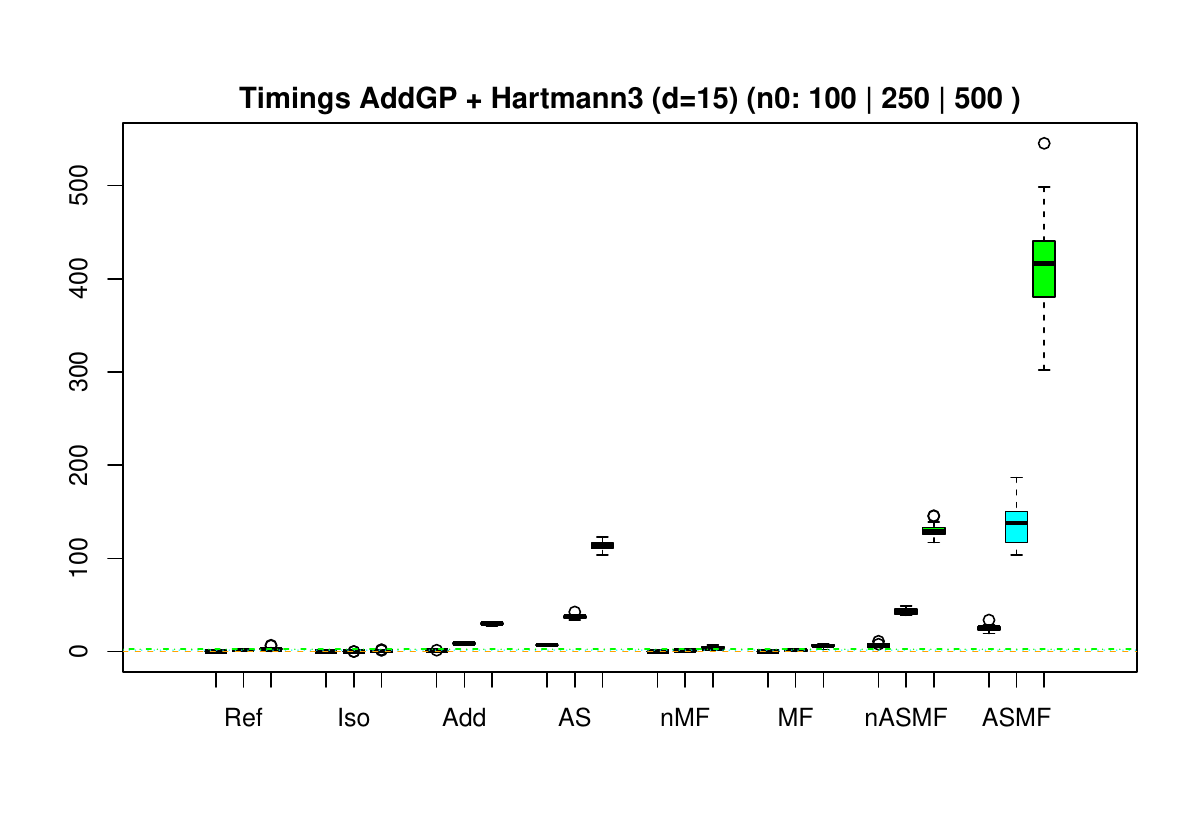}
\caption{Timings in seconds. The different $n_0$ values appear with different
 colors, the corresponding lines indicate the baseline result from standard
 GP models.}
\label{fig:timings}
\end{figure*}

To complement the results provided in the main part, Figure \ref{fig:nbactive}
focuses on the estimation of the low intrinsic dimensionality. Then a
comparison on the RMSE for the OAK model by \cite{Lu2022a} is given in Figure
\ref{fig:supprmseres}, before general timing results in Figure
\ref{fig:timings}. The experiments have been performed on four 2.40GHz Intel cores.

\end{document}